\documentclass[a4paper,twosided, 11.8pt]{amsart}

\usepackage{amsfonts}
\usepackage{amssymb}
\usepackage{latexsym}
\usepackage{graphics}
\usepackage[2emode]{psfrag}
\usepackage{mathrsfs}
\usepackage{amsthm}
\usepackage{amsmath}
\usepackage[all]{xy}
\usepackage{enumerate}
\usepackage[latin1]{inputenc}
\usepackage{lscape}
\usepackage{a4wide}

\xyoption{2cell} \xyoption{curve} \xyoption{rotate} \xyoption{v2}

\begin{document}

\newcommand{\thlabel}[1]{\label{th:#1}}
\newcommand{\thref}[1]{Theorem~\ref{th:#1}}
\newcommand{\selabel}[1]{\label{se:#1}}
\newcommand{\seref}[1]{Section~\ref{se:#1}}
\newcommand{\lelabel}[1]{\label{le:#1}}
\newcommand{\leref}[1]{Lemma~\ref{le:#1}}
\newcommand{\prlabel}[1]{\label{pr:#1}}
\newcommand{\prref}[1]{Proposition~\ref{pr:#1}}
\newcommand{\colabel}[1]{\label{co:#1}}
\newcommand{\coref}[1]{Corollary~\ref{co:#1}}
\newcommand{\relabel}[1]{\label{re:#1}}
\newcommand{\reref}[1]{Remark~\ref{re:#1}}
\newcommand{\exlabel}[1]{\label{ex:#1}}
\newcommand{\exref}[1]{Example~\ref{ex:#1}}
\newcommand{\delabel}[1]{\label{de:#1}}
\newcommand{\deref}[1]{Definition~\ref{de:#1}}
\newcommand{\eqlabel}[1]{\label{eq:#1}}
\newcommand{\equref}[1]{(\ref{eq:#1})}

\newcommand{\norm}[1]{\| #1 \|}
\def\N{\mathbb N}
\def\Z{\mathbb Z}
\def\Q{\mathbb Q}
\def\mod{\textit{\emph{~mod~}}}
\def\R{\mathcal R}
\def\S{\mathcal S}
\def\*C{{^*\mathcal C}}
\def\C{\mathcal C}
\def\D{\mathcal D}
\def\J{\mathcal J}
\def\M{\mathcal M}
\def\T{\mathcal T}

\newcommand{\Hom}{{\rm Hom}}
\newcommand{\End}{{\rm End}}
\newcommand{\Fun}{{\rm Fun}}
\newcommand{\Mor}{{\rm Mor}\,}
\newcommand{\Aut}{{\rm Aut}\,}
\newcommand{\Hopf}{{\rm Hopf}\,}
\newcommand{\Ann}{{\rm Ann}\,}
\newcommand{\Ker}{{\rm Ker}\,}
\newcommand{\Coker}{{\rm Coker}\,}
\newcommand{\im}{{\rm Im}\,}
\newcommand{\coim}{{\rm Coim}\,}
\newcommand{\Trace}{{\rm Trace}\,}
\newcommand{\Char}{{\rm Char}\,}
\newcommand{\Mod}{{\bf mod}}
\newcommand{\Spec}{{\rm Spec}\,}
\newcommand{\Span}{{\rm Span}\,}
\newcommand{\sgn}{{\rm sgn}\,}
\newcommand{\Id}{{\rm Id}\,}
\newcommand{\Com}{{\rm Com}\,}
\newcommand{\codim}{{\rm codim}}
\newcommand{\Mat}{{\rm Mat}}
\newcommand{\can}{{\rm can}}
\newcommand{\sign}{{\rm sign}}
\newcommand{\kar}{{\rm kar}}
\newcommand{\rad}{{\rm rad}}

\def\Ab{\underline{\underline{\rm Ab}}}
\def\lan{\langle}
\def\ran{\rangle}
\def\ot{\otimes}

\def\id{\textrm{{\small 1}\normalsize\!\!1}}
\def\To{{\multimap\!\to}}
\def\bigperp{{\LARGE\textrm{$\perp$}}}
\newcommand{\QED}{\hspace{\stretch{1}}
\makebox[0mm][r]{$\Box$}\\}

\def\AA{{\mathbb A}}
\def\BB{{\mathbb B}}
\def\CC{{\mathbb C}}
\def\DD{{\mathbb D}}
\def\EE{{\mathbb E}}
\def\FF{{\mathbb F}}
\def\GG{{\mathbb G}}
\def\HH{{\mathbb H}}
\def\II{{\mathbb I}}
\def\JJ{{\mathbb J}}
\def\KK{{\mathbb K}}
\def\LL{{\mathbb L}}
\def\MM{{\mathbb M}}
\def\NN{{\mathbb N}}
\def\OO{{\mathbb O}}
\def\PP{{\mathbb P}}
\def\QQ{{\mathbb Q}}
\def\RR{{\mathbb R}}
\def\SS{{\mathbb S}}
\def\TT{{\mathbb T}}
\def\UU{{\mathbb U}}
\def\VV{{\mathbb V}}
\def\WW{{\mathbb W}}
\def\XX{{\mathbb X}}
\def\YY{{\mathbb Y}}
\def\ZZ{{\mathbb Z}}

\def\aa{{\mathfrak A}}
\def\bb{{\mathfrak B}}
\def\cc{{\mathfrak C}}
\def\dd{{\mathfrak D}}
\def\ee{{\mathfrak E}}
\def\ff{{\mathfrak F}}
\def\gg{{\mathfrak G}}
\def\hh{{\mathfrak H}}
\def\ii{{\mathfrak I}}
\def\jj{{\mathfrak J}}
\def\kk{{\mathfrak K}}
\def\ll{{\mathfrak L}}
\def\mm{{\mathfrak M}}
\def\nn{{\mathfrak N}}
\def\oo{{\mathfrak O}}
\def\pp{{\mathfrak P}}
\def\qq{{\mathfrak Q}}
\def\rr{{\mathfrak R}}
\def\ss{{\mathfrak S}}
\def\tt{{\mathfrak T}}
\def\uu{{\mathfrak U}}
\def\vv{{\mathfrak V}}
\def\ww{{\mathfrak W}}
\def\xx{{\mathfrak X}}
\def\yy{{\mathfrak Y}}
\def\zz{{\mathfrak Z}}

\def\aaa{{\mathfrak a}}
\def\bbb{{\mathfrak b}}
\def\ccc{{\mathfrak c}}
\def\ddd{{\mathfrak d}}
\def\eee{{\mathfrak e}}
\def\fff{{\mathfrak f}}
\def\ggg{{\mathfrak g}}
\def\hhh{{\mathfrak h}}
\def\iii{{\mathfrak i}}
\def\jjj{{\mathfrak j}}
\def\kkk{{\mathfrak k}}
\def\lll{{\mathfrak l}}
\def\mmm{{\mathfrak m}}
\def\nnn{{\mathfrak n}}
\def\ooo{{\mathfrak o}}
\def\ppp{{\mathfrak p}}
\def\qqq{{\mathfrak q}}
\def\rrr{{\mathfrak r}}
\def\sss{{\mathfrak s}}
\def\ttt{{\mathfrak t}}
\def\uuu{{\mathfrak u}}
\def\vvv{{\mathfrak v}}
\def\www{{\mathfrak w}}
\def\xxx{{\mathfrak x}}
\def\yyy{{\mathfrak y}}
\def\zzz{{\mathfrak z}}

\newcommand{\aA}{\mathscr{A}}
\newcommand{\bB}{\mathscr{B}}
\newcommand{\cC}{\mathscr{C}}
\newcommand{\dD}{\mathscr{D}}
\newcommand{\eE}{\mathscr{E}}
\newcommand{\fF}{\mathscr{F}}
\newcommand{\gG}{\mathscr{G}}
\newcommand{\hH}{\mathscr{H}}
\newcommand{\iI}{\mathscr{I}}
\newcommand{\jJ}{\mathscr{J}}
\newcommand{\kK}{\mathscr{K}}
\newcommand{\lL}{\mathscr{L}}
\newcommand{\mM}{\mathscr{M}}
\newcommand{\nN}{\mathscr{N}}
\newcommand{\oO}{\mathscr{O}}
\newcommand{\pP}{\mathscr{P}}
\newcommand{\qQ}{\mathscr{Q}}
\newcommand{\rR}{\mathscr{R}}
\newcommand{\sS}{\mathscr{S}}
\newcommand{\tT}{\mathscr{T}}
\newcommand{\uU}{\mathscr{U}}
\newcommand{\vV}{\mathscr{V}}
\newcommand{\wW}{\mathscr{W}}
\newcommand{\xX}{\mathscr{X}}
\newcommand{\yY}{\mathscr{Y}}
\newcommand{\zZ}{\mathscr{Z}}

\newcommand{\Aa}{\mathcal{A}}
\newcommand{\Bb}{\mathcal{B}}
\newcommand{\Cc}{\mathcal{C}}
\newcommand{\Dd}{\mathcal{D}}
\newcommand{\Ee}{\mathcal{E}}
\newcommand{\Ff}{\mathcal{F}}
\newcommand{\Gg}{\mathcal{G}}
\newcommand{\Hh}{\mathcal{H}}
\newcommand{\Ii}{\mathcal{I}}
\newcommand{\Jj}{\mathcal{J}}
\newcommand{\Kk}{\mathcal{K}}
\newcommand{\Ll}{\mathcal{L}}
\newcommand{\Mm}{\mathcal{M}}
\newcommand{\Nn}{\mathcal{N}}
\newcommand{\Oo}{\mathcal{O}}
\newcommand{\Pp}{\mathcal{P}}
\newcommand{\Qq}{\mathcal{Q}}
\newcommand{\Rr}{\mathcal{R}}
\newcommand{\Ss}{\mathcal{S}}
\newcommand{\Tt}{\mathcal{T}}
\newcommand{\Uu}{\mathcal{U}}
\newcommand{\Vv}{\mathcal{V}}
\newcommand{\Ww}{\mathcal{W}}
\newcommand{\Xx}{\mathcal{X}}
\newcommand{\Yy}{\mathcal{Y}}
\newcommand{\Zz}{\mathcal{Z}}

\numberwithin{equation}{section}
\renewcommand{\theequation}{\thesection.\arabic{equation}}
\newcommand{\bara}[1]{\overline{#1}}
\newcommand{\ContFunt}[2]{\bara{\mathrm{Funt}}(#1,\,#2)}
\newcommand{\lBicomod}[2]{{}_{#1}\mM^{#2}}
\newcommand{\rBicomod}[2]{{}^{#1}\mM_{#2}}
\newcommand{\Bicomod}[2]{{}^{#1}\mM^{#2}}
\newcommand{\lrBicomod}[2]{{}^{#1}\mM^{#2}}
\newcommand{\rcomod}[1]{\mM^{#1}}
\newcommand{\rmod}[1]{\mM_{#1}}
\newcommand{\Bimod}[2]{{}_{#1}\mM_{#2}}
\newcommand{\Sf}[1]{\mathsf{#1}}
\newcommand{\Sof}[1]{S^{{\bf #1}}}
\newcommand{\Tof}[1]{T^{{\bf #1}}}
\renewcommand{\hom}[3]{\mathrm{Hom}_{#1}\left(\underset{}{}#2,\,#3\right)}
\newcommand{\Coint}[2]{\mathrm{Coint}(#1,#2)}
\newcommand{\IntCoint}[2]{\mathrm{InCoint}(#1,#2)}
\newcommand{\Coder}[2]{\mathrm{Coder}(#1,#2)}
\newcommand{\IntCoder}[2]{\mathrm{InCoder}(#1,#2)}
\newcommand{\lr}[1]{\left(\underset{}{} #1 \right)}
\newcommand{\Ext}[3]{\mathrm{Ext}_{\eE}^{#1}\lr{#2,\,#3}}
\newcommand{\equalizerk}[2]{\mathfrak{eq}_{#1,\,#2}^k}
\newcommand{\equalizer}[2]{\mathfrak{eq}_{#1,\,#2}}
\newcommand{\coring}[1]{\mathfrak{#1}}
\newcommand{\tensor}[1]{\otimes_{#1}}
\newcommand{\Mono}[1]{\rR_{(\coring{#1}\,:\,A)}}
\newcommand{\Rtensor}[1]{\underset{(\coring{#1}\,:\,A)}{\otimes}}
\newcommand{\tensork}[2]{#1\otimes_{\KK} #2}
\newcommand{\cotensor}[1]{\square_{#1}}
\newcommand{\tensorbajo}[1]{\underset{#1}{\otimes}}
\newcommand{\td}[1]{\widetilde{#1}}
\newcommand{\tha}[1]{\widehat{#1}}
\newcommand{\wrcomod}[2]{ \mM^{#1}_{#2}}
\newcommand{\bd}[1]{\boldsymbol{#1}}


\def\units{{\mathbb G}_m}
\def\rightact{\hbox{$\leftharpoonup$}}
\def\leftact{\hbox{$\rightharpoonup$}}

\def\*C{{}^*\hspace*{-1pt}{\Cc}}

\def\text#1{{\rm {\rm #1}}}

\def\smashco{\mathrel>\joinrel\mathrel\triangleleft}
\def\cosmash{\mathrel\triangleright\joinrel\mathrel<}

\def\ol{\overline}
\def\ul{\underline}
\def\dul#1{\underline{\underline{#1}}}
\def\Nat{\dul{\rm Nat}}
\def\Set{\dul{\rm Set}}

\renewcommand{\subjclassname}{\textup{2000} Mathematics Subject
     Classification}

\newtheorem{proposition}{Proposition}[section]
  \newtheorem{lemma}[proposition]{Lemma}
  \newtheorem{aussage}[proposition]{}
  \newtheorem{corollary}[proposition]{Corollary}
  \newtheorem{theorem}[proposition]{Theorem}

  \theoremstyle{definition}
  \newtheorem{definition}[proposition]{Definition}
  \newtheorem{example}[proposition]{Example}

  \theoremstyle{remark}
  \newtheorem{remark}[proposition]{Remark}

\title[Cowreath over Corings.]
{Extended distributive law: Co-Wreath over Co-Rings.}
\date{\today}
\author{L. El Kaoutit}
\address{Departamento de \'Algebra. Facultad de Educaci\'on y Humanidades de Ceuta.
Universidad de Granada. El Greco Nº 10. E-51002 Ceuta, Spain}
\email{kaoutit@ugr.es}

\keywords{ Bimodules. Corings. Comodules. Monoidal categories. Extended distributive law: Wreaths and co-Wreaths.\\
Research supported by grant MTM2004-01406 from the Ministerio de
Educaci\'{o}n y Ciencia of Spain. } 

\baselineskip 16pt

\begin{abstract}
A basic theory of cowreath or extended distributive laws in the
bicategory of unital bimodules, is deciphered. Precisely, we give in
terms of tensor product over a scalar base ring, a simplest and
equivalent definition for cowreath over coring and for comodule over
cowreath. An adjunction connecting the category of comodules over
the factor coring and the category of comodules over the coring
arising from the cowreath products is also given. The dual notions
i.e. wreaths over rings extension and their modules are included.
\end{abstract}

\maketitle

\section*{Introduction}
The notion of wreath in bicategory was introduced by S. Lack and R.
Street in \cite{Lack/Street:2002}. The notion of co-wreath is in
some sense dual to that of wreath. Explicitly, given a bicategory
$\Sf{B}$ (see \cite{Benabou:1967}), using \cite{Lack/Street:2002}
one can construct its (right) Eilenberg-Moore bicategory of comonads
denoted by $\Sf{REM(B)}$. A $0$-cell of $\Sf{REM(B)}$ is a pair
$(C,A)$ where $A$ is a $0$-cell of $\Sf{B}$ and $C$ is a comonad on
$A$ in $\Sf{B}$. A $1$-cell $(P,\ppp): (C,A) \to (D,B)$ consists of
a $1$-cell $P:A\to B$ and $2$-cell $\ppp: D \circ P \Rightarrow P
\circ C$ in $\Sf{B}$ satisfying two conditions with the counits and
comultiplications of the comonads $C$ and $D$. A $2$-cell (in
reduced form) $\varphi: (P,\ppp) \Rightarrow (Q,\qqq)$ is a $2$-cell
$\varphi: D \circ P \Rightarrow Q$ in $\Sf{B}$ satisfying one
condition with $\ppp$, $\qqq$ and the comultiplication of the
comonad $D$. A (right) \emph{cowreath} in $\Sf{B}$ is defined as in
\cite{Lack/Street:2002} to be a comonad in $\Sf{REM(B)}$.
Explicitly, a cowreath consists of a $0$-cell $A$ of $\Sf{B}$,
comonad $C$ on $A$ in $\Sf{B}$, a $1$-cell $R:A\to A$, and $2$-cells
$\rrr: C\circ R \Rightarrow R \circ C$, $\xi: C \circ R \Rightarrow
\II_A$ and $\delta: C \circ (R\circ R) \Rightarrow R$ satisfying
appropriate conditions. Notice that $R$ need not itself be a
comonad, but it could be while $\xi$ and $\delta$ could be obtained
from the counit and comultiplication of $R$: in this case $\rrr$ is
called a \emph{distributive law} between the comonad $C$ and $R$,
and due to J. Beck \cite{Beck:1969}.

In this paper we study co-wreath in the bicategory of bimodules
$\mathsf{Bim}$ ($0$-cells are unital rings, $1$-cells are unital
bimodules, $2$-cells are bilinear maps). For a given comonad in
$\Sf{Bim}$ i.e. coring, we review in Section \ref{Sec-1} its
Eilenberg-Moore monoidal category using $2$-cells in their unreduced
form. In Section \ref{Sect-DEf}, we give in terms of tensor product
over the base ring, a simplest and equivalent definition of
co-wreath (Proposition \ref{Equivalent-Def}). If the given coring
arises from some an entwining structure \cite{Brzezinski/Majid:1998}
we then establish a monoidal functor from the Eilenberg-Moore
monoidal category of the factor coalgebra to its Eilenberg-Moore
monoidal category. This gives a procedure to construct from the
commutative case, an examples of co-wreaths  with non commutative
base ring (Proposition \ref{ent-corona}). Comodules over co-wreath
and their morphisms are studied in Section \ref{Sect-Mod}, here
again we give an equivalent definitions in terms of tensor product
over base ring (Proposition \ref{comod} and \ref{comod1}). The
cowreath products is a coring with underlying bimodule a tensor
product of a base coring and a cowreath (Proposition
\ref{cow-prod}). Various functors and adjunctions are offered in
Section \ref{Sect-Mod} (diagram \eqref{digrama}). Section
\ref{Sect-Dual} contains a analogous results for wreath over rings
extension. A typical example of wreath is a distributive law between
two rings whose wreath products is known in the literature as
twisted tensor product algebra, see \cite{Tambara:1990, Van
Daele/Van Keer:1994, Cap/Schichl/Vanzura:1995,
Caenepeel/Ion/Militaru/Zhu:2000, LopezPena:2006} (subsection
\ref{TWALG}).

\bigskip
{\textsc{Notations and Basic Notions}:} Given any Hom-set category
$\cC$, the notation $X \in \aA$ means that $X$ is an object of
$\cC$. The identity morphism of $X$ will be denoted by $X$ itself.
The set of all morphisms $f:X \to X'$ in $\cC$, is denoted by
$\hom{\cC}{X}{X'}$. We work over a commutative ground ring with $1$
denoted by $\KK$. All rings are assumed to be associative
$\KK$--algebras. Modules are unital modules, and bimodules are left
and right unital modules and are assumed to be central
$\KK$--modules. Given $A$ and $B$ two rings, the category of $(A,
B)$-bimodules (left $A$-modules and right $B$-modules) is denoted as
usual by $\Bimod{A}{B}$. The $\KK$--module morphisms in this
category will be denoted by ${\rm Hom}_{A-B}\lr{-,-}$. The symbol
$-\tensor{A}-$ (or sometimes $-\otimes-$ depending on the context)
between bimodules and bilinear map denotes the tensor product over
$A$. Let $A$ be a ring an \emph{$A$-coring} \cite{Sweedler:1975} is
a three-tuple $(\coring{C}, \Delta, \varepsilon)$ consisting of an
$A$-bimodule $\cc$ and two $A$-bilinear maps
$$\Delta: \coring{C} \longrightarrow \coring{C} \tensor{A} \coring{C}\quad
\textrm{and}\quad \varepsilon: \coring{C} \longrightarrow A,$$ known
as the comultiplication and the counit of $\cc$, which satisfy
$$(\coring{C}\tensor{A}\Delta) \circ \Delta \,\, =\,\, (\Delta
\tensor{A} \coring{C}) \circ \Delta,\quad (\coring{C} \tensor{A}
\varepsilon) \circ \Delta \,\,=\,\, \coring{C} \,\,=\,\,
(\varepsilon\tensor{A}\coring{C}) \circ \Delta.$$ A Sweedler's
notation for comultiplication will be used i.e. $\Delta(c)\,=\,
c_{(1)}\tensor{A}c_{(2)}$ (summation is understood), for every $c
\in \cc$. We use the notation $(\cc:A)$ for an $A$--coring $\cc$. An
$A$--bilinear map $\phi: \dd \to \cc$ is a \emph{morphism of an
$A$-corings} if it satisfies $\varepsilon_{\cc} \circ \phi\,=\,
\varepsilon_{\dd}$ and $(\phi\tensor{A}\phi) \circ
\Delta_{\dd}\,=\,\Delta_{\cc} \circ \phi$.  A \emph{right
$\coring{C}$-comodule} is a pair $(M,\rho^M)$ with $M$ a right
$A$-module and $\rho^M: M \to M\tensor{A}\coring{C}$ a right
$A$-linear map (called right $\coring{C}$-coaction) satisfying two
equalities $(\rho^M \tensor{A} \coring{C}) \circ \rho^M \,=\,
(M\otimes_A \Delta) \circ \rho^M$, and $(M\otimes_A \varepsilon)
\circ \rho^M \,=\,M$. The Sweedler's notation for right
$\cc$-coactions is as usual: $\rho^M(m)\,=\,
m_{(0)}\tensor{A}m_{(1)}$ (summation is understood), for every $m
\in M$. A \emph{morphism} of right $\coring{C}$-comodules
$f:(M,\rho^M) \to (M',\rho^{M'})$ is a right $A$-linear map $f:M\to
M'$ which is compatible with coactions: $f \circ \varrho^{M'} \,=\,
(f\tensor{A} \coring{C}) \circ \varrho^{M}$ ($f$ is right
$\coring{C}$-colinear). The $\KK$--module of all colinear maps will
be denoted by ${\rm Hom}_{\cc}\lr{M,M'}$. We denote by
$\rcomod{\coring{C}}$ the category of all right
$\coring{C}$-comodules. Left $\coring{C}$-comodules are
symmetrically defined, we use the Greek letter $\lambda^{-}$ to
denote theirs coactions. The category of (right) $\cc$-comodules is
not in general an abelian category, its has cokernels and arbitrary
direct sums which can be already computed in the category of
$A$-modules. However, if ${}_A\cc$ is a flat module, then
$\rcomod{\cc}$ becomes a Grothendieck category (see
\cite{Kaoutit/Gomez/Lobillo:2004c, Brzezinski/Wisbauer:2003}). Let
$\dd$ be a $B$-coring, a $(\cc,\dd)$-\emph{bicomodule} is a
three-tuple $(M,\rho^M,\lambda^M)$ consisting of an $(A,B)$-bimodule
($M \in \Bimod{A}{B}$) and $A-B$-bilinear maps $\rho^M: M \to M
\otimes_B\coring{D}$ and $\lambda^M: M \to \coring{C}\otimes_A M$
such that $(M, \rho^M)$ is right $\coring{D}$-comodule and
$(M,\lambda^M)$ is left $\coring{C}$-comodule with compatibility
condition $(\coring{C}\otimes_A \rho^M) \circ \lambda^M  \,=\,
(\lambda^M \otimes_B \coring{D}) \circ \rho^M$. A \emph{morphism of
bicomodules} is a left and right colinear map (or \emph{bicolinear
map}); we use the notation ${\rm Hom}_{\cc-\dd}\lr{M,M'}$ for the
$\KK$--module of all bicolinear maps. The category of all
$(\coring{C},\coring{D})$-bicomodule is denoted by
$\Bicomod{\coring{C}}{\coring{D}}$. Obviously any ring $A$ can be
endowed with a trivial structure of an $A$-coring with
comultiplication the isomorphism $A\cong A\tensor{A}A$ and counit
the identity $A$. In this way an $(A,\dd)$-bicomodule is just a
right $\dd$-comodule $(M,\rho^M)$ whose underlying module $M$ is an
$A-B$-bimodule and whose coaction $\rho^M$ is an $A-B$-bilinear map.
The category of $(A,\dd)$--bicomodules is denoted by
${}_A\mM^{\dd}$. For more details on comodules, definitions and
basic properties of bicomodules and the cotensor product, the reader
is referred to monograph \cite{Brzezinski/Wisbauer:2003}.

\section{Eilenberg-Moore monoidal category associated to a
coring.}\label{Sec-1}

In all this section the symbole $-\tensor{}-$ stands by
$-\tensor{A}-$ the tensor product bi-funtor over a ring $A$. For a
fixed coring $(\cc:A)$  with comultiplication $\Delta$ and counit
$\varepsilon$. Consider, as in \cite{Brzezinski/Kaoutit/Gomez:2006}
(see \cite{Lack/Street:2002} for general notions), the pre-additive
category $\Mono{C}$ defined by the following data:
\begin{enumerate}[$\bullet$]
\item \emph{Objects}: Are pairs $(M,\mmm)$ consisting of an
$A$-bimodule $M$ and $A$-bilinear map $\mmm: \cc\tensor{}M \to
M\tensor{}\cc$ such that
\begin{eqnarray}
  (M\tensor{}\Delta) \circ \mmm  &=& (\mmm\tensor{}\cc) \circ (\cc\tensor{}\mmm)
  \circ (\Delta\tensor{}M) \label{1-cell} \\
  (M\tensor{}\varepsilon) \circ \mmm &=& \varepsilon \tensor{}M, \label{1-cell'}
\end{eqnarray}
where in the second equality $M$ was identified with  $A \tensor{}M$
and with $M\tensor{}A$ via the obvious isomorphism.

\item \emph{Morphisms}: Given any object $(M,\mmm)$ one can easily
check that $\cc\tensor{}M$ is a $\cc$-bicomodule with left
$\cc$-coaction $\lambda^{\cc\tensor{}M}\,=\,\Delta\tensor{}M$ and
right $\cc$-coaction $\varrho^{\cc\tensor{}M}\,=\,
(\cc\tensor{}\mmm) \circ (\Delta\tensor{}M)$. By \cite[Proposition
2.2]{Brzezinski/Kaoutit/Gomez:2006}, the $\KK$--module of
homomorphisms in $\Mono{C}$ are then defined (in unreduced form) by
$$\hom{\Mono{C}}{(M,\mmm)}{(M',\mmm')} \,\, :=\,\,
\hom{\cc-\cc}{\cc\tensor{}M}{\cc\tensor{}M'}.$$
\end{enumerate}
That is a  morphism $\varphi:(M,\mmm) \to (M',\mmm')$ in
$\Mono{C}$ is an $A$-bilinear map $\varphi: \cc\tensor{}M \to
\cc\tensor{}M'$ which satisfies
\begin{eqnarray}
  (\Delta\tensor{}M') \circ \varphi &=&  (\cc\tensor{}\varphi) \circ (\Delta\tensor{}M) \label{2-cell}\\
  (\cc\tensor{}\mmm') \circ (\Delta\tensor{}M') \circ \varphi &=&
  (\varphi\tensor{}\cc) \circ (\cc\tensor{}\mmm) \circ
  (\Delta\tensor{}M). \label{2-cell'}
\end{eqnarray}

\begin{remark}\label{i-ii}
Let $(\cc:A)$ be any coring and $M$ an $A$--bimodule. Then one can
easily check that the following statements are equivalent.
\begin{enumerate}[(i)]
\item $\cc\tensor{A}M$ is a $\cc$-bicomodule with left
$\cc$-coaction $\lambda^{\cc\tensor{}M}\,=\, \Delta\tensor{}M$;

\item there is an $A$-bilinear map $\mmm: \cc\tensor{}M \to
M\tensor{}\cc$ such that $(M,\mmm)$ is an object of the category
$\Mono{C}$.
\end{enumerate}
\end{remark}

The category $\Mono{C}$ is in fact a multiplicative additive
category (i.e. an additive Monoidal category). Its multiplication is
defined as follows.

\begin{proposition}\label{monoidal}
For any coring $(\cc:A)$, the category $\Mono{C}$ is a monoidal
category with multiplication defined as follows.
\begin{enumerate} [ {}]
\item Given two object $(M,\mmm)$ and $(M',\mmm')$ of $\Mono{C}$, we define a
new object of $\Mono{C}$
$$(M,\mmm) \Rtensor{C} (M',\mmm') \,\,:=\,\,
\left(\underset{}{}M\tensor{}M',(M\tensor{}\mmm') \circ
(\mmm\tensor{}M') \right)$$

\item If $\varphi:(M,\mmm) \to (M',\mmm')$ and $\psi:(N,\nnn) \to
(N',\nnn')$ are morphisms in $\Mono{C}$, then their
multiplication is defined by the composition\\
\xy *+{\cc\tensor{}M\tensor{}N}="p",
p+<12cm,0pt>*+{\cc\tensor{}M'\tensor{}N'}="1",
p+<0pt,-2cm>*+{\cc\tensor{}\cc\tensor{}M\tensor{}N}="2",
p+<12cm,-2cm>*+{\cc\tensor{}M'\tensor{}\cc\tensor{}N'}="3",
p+<3cm,-4cm>*+{\cc\tensor{}\cc\tensor{}M'\tensor{}N}="4",
p+<9cm,-4cm>*+{\cc\tensor{}M'\tensor{}\cc\tensor{}N}="5", {"p"
\ar@{-->}^-{\varphi \Rtensor{C} \psi} "1"}, {"p"
\ar@{->}_-{\Delta\tensor{}M\tensor{}N} "2"}, {"2"
\ar@{->}_-{\cc\tensor{}\varphi \tensor{}N} "4"}, {"4"
\ar@{->}_-{\cc\tensor{}\mmm'\tensor{}N} "5"}, {"5"
\ar@{->}_-{\cc\tensor{}M'\tensor{}\psi} "3"}, {"3"
\ar@{->}_-{\cc\tensor{}M'\tensor{}\varepsilon\tensor{}N'} "1"}
\endxy
\end{enumerate}
Equivalently $\varphi \Rtensor{C} \psi\,=\,
(\cc\tensor{}M'\tensor{}\varepsilon\tensor{}N) \circ
(\varphi\tensor{}\psi) \circ (\cc\tensor{}\mmm\tensor{}N) \circ
(\Delta\tensor{}M\tensor{}N)$. The identity object for this
multiplication is proportioned by the pair $(A,\cc)$ where $\cc$ was
identified with $\cc\tensor{}A \cong A\tensor{}\cc$.
\end{proposition}
\begin{proof}
It is clearly seen that $(M,\mmm) \Rtensor{C}(M',\mmm')$ belongs to
$\Mono{C}$, whenever $(M,\mmm)$ and $(M',\mmm')$ are objects of
$\Mono{C}$. By definition $\varphi \Rtensor{C}\psi$ is an
$A$-bilinear map. So we need to check its compatibility with left
and right $\cc$-coactions of both $\cc$--bicomodules
$\cc\tensor{}M\tensor{}N$ and $\cc\tensor{}M'\tensor{}N'$. We know
that $\varrho^{\cc\tensor{}M'\tensor{}N'}\,=\,
(\cc\tensor{}M'\tensor{}\nnn') \circ (\cc\tensor{}\mmm'\tensor{}N')
\circ (\Delta\tensor{}M'\tensor{}N')$. Applying equation
\eqref{1-cell} to $\mmm'$ and \eqref{2-cell} to both $\varphi$ and
$\psi$, we get
$$\varrho^{\cc\tensor{}M'\tensor{}N'} \circ
(\varphi\Rtensor{C}\psi)\,=\, (\cc\tensor{}M'\tensor{}\nnn') \circ
(\cc\tensor{}M'\tensor{}\psi)   \circ (\cc\tensor{}\mmm'\tensor{}N)
\circ (\cc\tensor{}\varphi\tensor{}N) \circ
(\Delta\tensor{}M\tensor{}N)$$ Analogously, we obtain
$$\lr{(\varphi\Rtensor{C}\psi) \tensor{}\cc} \circ
\varrho^{\cc\tensor{}M\tensor{}N} \,=\,
(\cc\tensor{}M'\tensor{}\nnn') \circ (\cc\tensor{}M'\tensor{}\psi)
\circ (\cc\tensor{} \mmm'\tensor{}N) \circ
(\cc\tensor{}\varphi\tensor{}N) \circ  (\Delta\tensor{}M
\tensor{}N)$$ Therefore, $ \varrho^{\cc\tensor{}M'\tensor{}N'} \circ
(\varphi\Rtensor{C}\psi)=((\varphi\Rtensor{C}\psi) \tensor{}\cc)
\circ \varrho^{\cc\tensor{}M\tensor{}N}$. That is
$(\varphi\Rtensor{C}\psi)$ is right $\cc$-colinear. This map is
clearly left $\cc$--colinear.
\end{proof}

\begin{example}
Let $(\cc:A)$ be any coring. Define the map
$$\xymatrix{ \ccc: \cc\tensor{A}\cc \ar@{->}[r] &
\cc\tensor{A}\cc & & \lr{c\tensor{}c' \longmapsto \,\,c_{(1)}
\tensor{}c_{(2)} \varepsilon(c') +
\varepsilon(c)c'_{(1)}\tensor{}c'_{(2)} - c\tensor{}c'}. }$$ Then
$(\cc,\ccc)$ is an object of $\Mono{C}$. By definition $\ccc$ is
an $A$--bilinear map. It is clear that $(\cc\tensor{}\varepsilon)
\circ \ccc = \varepsilon \tensor{}\cc$. Now, we have
\begin{multline*}
(\ccc\tensor{}\cc) \circ (\cc\tensor{}\ccc) \circ
(\Delta\tensor{}\cc) (c\tensor{}c') \,\,=\, (\ccc\tensor{}\cc)
\circ (\cc\tensor{}\ccc) (c_{(1)}\tensor{}c_{(2)} \tensor{}c')  \\
\,=\, \ccc\tensor{}\cc \lr{ c_{(1)}
\tensor{}\lr{c_{(2)}\tensor{}c_{(3)} \varepsilon(c') +
\varepsilon(c_{(2)})c'_{(1)}\tensor{}c'_{(2)} -
c_{(2)}\tensor{}c'} }  \\ \,=\, \ccc\tensor{}\cc \lr{ c_{(1)}
\tensor{}c_{(2)}\tensor{}c_{(3)}
\varepsilon(c') + c\tensor{}c'_{(1)}\tensor{}c'_{(2)} - c_{(1)}\tensor{}c_{(2)}\tensor{}c' }  \\
\,=\,
c_{(1)}\tensor{}c_{(2)}\varepsilon(c_{(3)})\tensor{}c_{(4)}\varepsilon(c')+
\varepsilon(c_{(1)})c_{(2)}\tensor{}c_{(3)}\tensor{}c_{(4)}\varepsilon(c')
- c_{(1)}\tensor{}c_{(2)}\tensor{}c_{(3)}\varepsilon(c') +
c_{(1)}\tensor{}c_{(2)}\varepsilon(c'_{(1)})\tensor{}c'_{(2)} \\
\,\,\,\, + \varepsilon(c)c'_{(1)}\tensor{}c'_{(2)}\tensor{}c'_{(3)}
- c\tensor{}c'_{(1)}\tensor{}c'_{(2)}
-c_{(1)}\tensor{}c_{(2)}\varepsilon(c_{(3)})\tensor{}c'-
\varepsilon(c_{(1)})c_{(2)}\tensor{}c_{(3)}\tensor{}c' +
c_{(1)}\tensor{}c_{(2)}\tensor{}c' \\ \,=\,
c_{(1)}\tensor{}c_{(2)}\tensor{}c_{(3)}\varepsilon(c')+
c_{(1)}\tensor{}c_{(2)}\tensor{}c_{(3)}\varepsilon(c') -
c_{(1)}\tensor{}c_{(2)}\tensor{}c_{(3)}\varepsilon(c') +
c_{(1)}\tensor{}c_{(2)}\tensor{}c'  +
\varepsilon(c)c'_{(1)}\tensor{}c'_{(2)}\tensor{}c'_{(3)} \\
\,\,\,\, - c\tensor{}c'_{(1)}\tensor{}c'_{(2)}
-c_{(1)}\tensor{}c_{(2)}\tensor{}c'-
c_{(1)}\tensor{}c_{(2)}\tensor{}c' +
c_{(1)}\tensor{}c_{(2)}\tensor{}c'
\\ \,=\,
c_{(1)}\tensor{}c_{(2)}\tensor{}c_{(3)}\varepsilon(c')  +
\varepsilon(c)c'_{(1)}\tensor{}c'_{(2)}\tensor{}c'_{(3)} -
c\tensor{}c'_{(1)}\tensor{}c'_{(2)} \,=\, (\cc\tensor{}\Delta) \circ
\ccc(c\tensor{}c'),
\end{multline*} for every pair $(c,c') \in \cc \times \cc$.
Therefore $(\cc,\ccc)$ satisfies equalities \eqref{1-cell} and
\eqref{1-cell'}, and so $(\cc,\ccc)$ is an object of $\Mono{C}$
(this is dual to \cite[proposition 1.7]{Menini/Dragos:2003}).\\
Given $(\dd:A)$ another coring and suppose that there is an
$A$-coring morphism $\phi: \dd \to \cc$. Define the map
$$\xymatrix{ \ddd: \cc\tensor{A}\dd \ar@{->}[r] &
\dd\tensor{A}\cc & & \lr{c\tensor{}d \longmapsto \,\,\varepsilon(c)
d_{(1)}\tensor{}\phi(d_{(2)})}. }$$ It is clear that $\ddd$ is an
$A$-bilinear, and that $(\dd\tensor{}\varepsilon) \circ \ddd =
\varepsilon \tensor{}\dd$. Furthermore, for every pair of elements
$(c,d) \in \cc \times \dd$, we have
\begin{eqnarray*}
  (\ddd\tensor{}\cc) \circ (\cc\tensor{}\ddd) \circ (\Delta\tensor{}\dd) (c\tensor{}d)  &=&
  (\ddd\tensor{}\cc) \circ (\cc\tensor{}\ddd)(c_{(1)} \tensor{}c_{(2)} \tensor{}d) \\
   &=& \ddd\tensor{}\cc\lr{c_{(1)}\tensor{}\varepsilon(c_{(2)})d_{(1)}\tensor{}\phi(d_{(2)}) }\\
   &=& \varepsilon(c)d_{(1)}\tensor{}\phi(d_{(2)})\tensor{}\phi(d_{(2)}) \\
   &=& (\dd\tensor{}\Delta) \circ \ddd(c\tensor{}d).
\end{eqnarray*} That is $(\dd,\ddd)$ satisfies equalities \eqref{1-cell} and
\eqref{1-cell'}, and hence $(\dd,\ddd)$ is an object of
$\Mono{C}$. This in fact comes from a general setting. Namely, if
we have any $(A',A)$--corings morphism $(\phi,\varphi): (\cc':A')
\to (\cc:A)$ in the sense of \cite{Gomez:2002}. That is  $\varphi:
A' \to A$ is a rings morphism and $\phi:\cc' \to \cc$ is by scalar
restriction an $A'$--bilinear map satisfying $\varepsilon_{\cc}
\circ \phi \,=\, \varphi \circ \varepsilon_{\cc'}$ and
$\Delta_{\cc} \circ \phi \,=\, \omega_{A',A} \circ
(\phi\tensor{A'}\phi) \circ \Delta_{\cc'}$, where $\omega_{A',A}$
is the obvious map. Then one can prove that
$(A\tensor{A'}\cc'\tensor{A'}A, \mmm)$ is an object of $\Mono{C}$,
where the map $\mmm: \cc\tensor{A'}\cc'\tensor{A'}A \to
A\tensor{A'}\cc'\tensor{A'}\cc$ sends $c\tensor{A'}c'\tensor{A'}a
\mapsto
\varepsilon(c)\tensor{A'}c'_{(1)}\tensor{A'}\phi(c'_{(2)})a$, for
every  element $a \in A$, $c \in \cc$ and $c' \in \cc'$.
\end{example}

Recall from \cite{Brzezinski/Majid:1998} that an entwining
structure over $\KK$ is a three-tuple $(A,C)_{\aaa}$ consisting of
$\KK$--algebra $A$ with multiplication $\mu$ and unit $1$,
$\KK$--coalgebra $C$ with comultiplication $\Delta$ and counit
$\varepsilon$, and a $\KK$--module map $\aaa: C\tensor{\KK}A
\rightarrow A\tensor{\KK}C$ satisfying
\begin{eqnarray}
  \aaa \circ (\tensork{C}{\mu}) &=& (\tensork{\mu}{C}) \circ
(\tensork{A}{\aaa}) \circ (\tensork{\aaa}{A}), \label{ent-10} \\
  \aaa \circ
(\tensork{C}{1}) &=& \tensork{1}{C}; \label{ent-11} \\
(\tensork{A}{\Delta}) \circ \aaa &=& (\tensork{\aaa}{C}) \circ
(\tensork{C}{\aaa}) \circ (\tensork{\Delta}{A}), \label{ent-20}\\
(\tensork{A}{\varepsilon}) \circ \aaa &=& \tensork{\varepsilon}{A}
\label{ent-21}.
\end{eqnarray} By \cite[Proposition 2.2]{Brzezinski:2002} the corresponding
$A$--coring is the $A$-bimodule $\coring{C} = A\tensor{\KK} C$ with
obvious left $A$--action and its right $A$--action is given by
$(a'\tensor{\KK}c).a=a'\aaa(c\tensor{\KK}a)$, for every $a,a' \in
A$, $c \in C$. The comultiplication map is $\Delta_{\coring{C}}=
A\tensor{\KK}\Delta$ and its counit is $\varepsilon_{\coring{C}}=
A\tensor{\KK}\varepsilon$. For instance, assume a $\KK$--bialgebra
$\Hh$ is given together with a right $\Hh$--comodule algebra $A$ and
right $\Hh$--module coalgebra. That is the right coaction $\rho^A: A
\to A\tensor{\KK}\Hh$ is a $\KK$-algebras morphism, while the right
action $\centerdot: C\tensor{\KK}\Hh \to C$ is a $\KK$--coalgebra
morphism. It is clearly seen that the map $\aaa: C\tensor{\KK}A \to
A\tensor{\KK}C$ sending $c\tensor{}a \mapsto a_{(1)}
\tensor{}(c\centerdot a_{(2)})$ (summation is understood), for every
$c \in C$ and $a \in A$, satisfies all equalities
\eqref{ent-10}-\eqref{ent-21}. Thus $(A,C)_{\aaa}$ is an entwining
structure over $\KK$.

Given any entwining structure $(A,C)_{\aaa}$ over $\KK$, one can
immediately check that $(A,\aaa)$ is an object of $\rR_{(C:K)}$.
Furthermore, we have

\begin{lemma}\label{ent-R}
Let $(A,C)_{\aaa}$ be an entwining structure over $\KK$. There is
a functor $A\tensor{}-\tensor{}A: \rR_{(C:\,\KK)} \to \Mono{C}$,
defined over objects by $(N,\nnn) \to
\lr{A\tensor{\KK}N\tensor{\KK}A,\,
(A\tensor{\KK}N\tensor{\KK}\aaa) \circ
(A\tensor{\KK}\nnn\tensor{\KK}A)}$ (up to canonical isomorphisms),
and over morphisms by  $f \to A\tensor{\KK}f\tensor{\KK}A$
\end{lemma}
\begin{proof}
Straightforward.
\end{proof}

\begin{remark}\label{left-right}
Reversing the twist maps, one can construct, for any coring
$(\cc:A)$, another monoidal category denoted by $\lL_{(\cc:\,A)}$.
The objects of $\lL_{(\cc:\,A)}$ are pairs $(\lll,L)$ consisting
of an $A$--bimodule $L$ and $A$--bilinear map $\lll:
L\tensor{A}\cc \to \cc\tensor{A}L$ compatible with the
comultiplication and the counit, i.e. satisfies the equalities
\begin{eqnarray}
  (\varepsilon\tensor{}L) \circ \lll &=& L\tensor{}\cc \label{lr-1} \\
  (\Delta\tensor{}L) \circ \lll &=& (\cc\tensor{}\lll) \circ
  (\lll\tensor{}\cc) \circ (L\tensor{}\Delta). \label{lr-2}
\end{eqnarray}
Notice, that here the $\cc$--bicomodule structure over
$L\tensor{}\cc$ is given by
$\varrho^{L\tensor{}\cc}=L\tensor{}\Delta$ and
$\lambda^{L\tensor{}\cc}= (\lll\tensor{}\cc) \circ
(L\tensor{}\Delta)$. The $\KK$--modules of morphisms in this
category are defined by
$$\hom{\lL_{(\cc:\,A)}}{(\lll,L)}{(\lll',L')} \,\, :=\,\,
\hom{\cc-\cc}{L\tensor{}\cc}{L'\tensor{}\cc}.$$ The
multiplications of this monoidal category are defined as follows:
Given $\gamma: (\lll,L) \to (\lll',L')$ and $\sigma: (\kkk,K) \to
(\kkk',K')$ two morphisms in $\lL_{(\cc:\,A)}$, the vertical
multiplication is defined by $$ (\lll,L) \tensorbajo{(\cc:\,A)}
(\kkk,K)\,\,=\,\, \lr{(\lll\tensor{}K) \circ (L\tensor{}\kkk),\,
L\tensor{}K}$$ and the vertical multiplication is defined by
\begin{equation}
\gamma\tensorbajo{(\cc:\,A)}\sigma\,\,=\,\,
(L'\tensor{}\varepsilon\tensor{}K'\tensor{}\cc) \circ
(\gamma\tensor{}K'\tensor{}\cc) \circ
(L\tensor{}\kkk'\tensor{}\cc) \circ (L\tensor{}\sigma\tensor{}\cc)
\circ (L\tensor{}K\tensor{}\Delta).
\end{equation}
\end{remark}

\section{Equivalent definitions of cowreath in
$\mathsf{Bim}$.}\label{Sect-DEf}

In this section we give in terms of the tensor product over the base
ring $A$, an equivalent definition of cowreath over a given coring
$(\cc:A)$. If our coring arises from entwining structures
\cite{Brzezinski/Majid:1998}, we then prove a procedure to construct
a new cowreath from a given cowreath over the factor coalgebra.

\begin{definition}\label{def}
Let $(\cc:A)$ be a coring. A \emph{cowreath over} $(\cc:A)$ (or
$\cc$-\emph{cowreath}) is coalgebra in the monoidal additive
category $\Mono{C}$ defined in Section \ref{Sec-1}. A \emph{wreath
over} $\cc$ (or $\cc$--\emph{wreath}) is an algebra in $\Mono{C}$.
Notice, that here in fact we are defining a \emph{right wreath} and
\emph{right cowreath}. The left notions are defined in the monoidal
category $\lL_{(\cc\,:A)}$ of Remark \ref{left-right}.
\end{definition}

\begin{proposition}\label{Equivalent-Def}
Let $(\cc:A)$ be a coring, and $(M,\mmm)$ an object of $\Mono{C}$.
The following statements are equivalent
\begin{enumerate}[(i)]
\item $(M,\mmm)$ is a $\cc$--cowreath.

\item There is a $\cc$--bicolinear maps $\xi:\cc\tensor{A}M \to
\cc$ and $\delta: \cc\tensor{A}M \to \cc\tensor{A}M\tensor{A}M$
rendering commutative the following diagrams
$$
\xy *+{\cc\tensor{}M}="p",
p+<3cm,0pt>*+{\cc\tensor{}M\tensor{}M}="1",
p+<3cm,-2cm>*+{\cc\tensor{}M}="2",{"p" \ar@{=} "2"}, {"p"
\ar@{->}^-{\delta} "1"}, {"1" \ar@{->}^-{\xi\tensor{}M} "2"}
\endxy \qquad \xy *+{\cc\tensor{}M}="p",
p+<4cm,0pt>*+{\cc\tensor{}M\tensor{}M}="1",
p+<4cm,-2cm>*+{M\tensor{}\cc\tensor{}M}="2",
p+<0pt,-2cm>*+{M\tensor{}\cc}="3", {"p" \ar@{->}_-{\mmm} "3"}, {"p"
\ar@{->}^-{\delta} "1"}, {"1" \ar@{->}^-{\mmm\tensor{}M} "2"}, {"2"
\ar@{->}^-{M\tensor{}\xi} "3"}
\endxy $$
$$
\xy *+{\cc\tensor{}M}="p",
p+<4cm,0pt>*+{\cc\tensor{}M\tensor{}M}="1",
p+<8cm,0pt>*+{\cc\tensor{}M\tensor{}M\tensor{}M}="2",
p+<0pt,-2cm>*+{\cc\tensor{}M\tensor{}M}="3",
p+<4cm,-2cm>*+{M\tensor{}\cc\tensor{}M}="4",
p+<8cm,-2cm>*+{M\tensor{}\cc\tensor{}M\tensor{}M}="5",  {"p"
\ar@{->}^-{\delta} "1"}, {"p" \ar@{->}_-{\delta} "3"}, {"1"
\ar@{->}^-{\delta\tensor{}M} "2"}, {"2"
\ar@{->}^{\mmm\tensor{}M\tensor{}M} "5"}, {"3"
\ar@{->}_-{\mmm\tensor{}M} "4"}, {"4" \ar@{->}_-{M\tensor{}\delta}
"5"}
\endxy
$$
\end{enumerate}
\end{proposition}
\begin{proof}
Given a coalgebra $(M,\mmm)$ in $\rR_{(\cc\,:A)}$ is equivalent to
given two morphisms $\xi: (M,\mmm) \to (A,\cc)$ and $\delta:
(M,\mmm) \to (M,\mmm) \Rtensor{C} (M,\mmm)$ in $\rR_{(\cc:\,A)}$
satisfying the usual coassociativity and counitary properties. That
is, up to natural isomorphisms, the following equalities are
verified
\begin{eqnarray*}
 \lr{(M,\,\mmm)\Rtensor{C}\xi} \circ \delta  &=& (M,\mmm) \,\, =\,\, \lr{ \xi\Rtensor{C}(M,\,\mmm)} \circ \delta   \\
  \lr{(M,\,\mmm) \Rtensor{C}\delta} \circ \delta &=&
  \lr{\delta\Rtensor{C}(M,\,\mmm)} \circ \delta.
\end{eqnarray*}
Of course $\xi:\cc\tensor{}M \to \cc$ and $\delta: \cc\tensor{}M
\to \cc\tensor{}M\tensor{}M$ are two $\cc$--bicolinear maps. Using
the definition of the multiplication  $-\Rtensor{C}-$ given in
Proposition \ref{monoidal}, the previous equalities are equivalent
to the following ones: the counitary properties are
\begin{eqnarray}
(\cc\tensor{}M\tensor{}\varepsilon) \circ
(\cc\tensor{}M\tensor{}\xi) \circ
  (\cc\tensor{}\mmm\tensor{}M) \circ (\Delta\tensor{}M\tensor{}M)
  \circ \delta
  &=& \cc\tensor{}M \label{counit-pr1} \\
(\cc\tensor{}\varepsilon\tensor{}M) \circ
(\cc\tensor{}\xi\tensor{}M) \circ (\Delta\tensor{}M\tensor{}M) \circ
\delta &=& \cc\tensor{}M \label{counit-pr2}
\end{eqnarray}
and the coassociativity is
\begin{multline}\label{coass-pr}
  (\cc\tensor{}M\tensor{}M\tensor{}\varepsilon\tensor{}M) \circ
  (\cc\tensor{}M\tensor{}\mmm\tensor{}M) \circ
  (\cc\tensor{}\mmm\tensor{}M\tensor{}M) \circ
  (\cc\tensor{}\delta\tensor{}M) \circ (\Delta\tensor{}M\tensor{}M)
  \circ \delta \\ =\,\, (\cc\tensor{}M \tensor{}
  \varepsilon\tensor{}M \tensor{} M) \circ
  (\cc\tensor{}M\tensor{}\delta) \circ (\cc\tensor{}\mmm\tensor{}M)
  \circ (\Delta\tensor{}M\tensor{}M) \circ \delta
\end{multline}
In conclusion $(M,\mmm)$ is a coalgebras in $\Mono{C}$ if and only
if there exist $\xi$ and $\delta$ satisfying equalities
\eqref{counit-pr1}-\eqref{coass-pr}.

$(i) \Rightarrow (ii)$. By equation \eqref{counit-pr2}, we have
\begin{eqnarray*}
  \cc\tensor{}M &=& (\cc\tensor{}\varepsilon\tensor{}M) \circ (\cc\tensor{}\xi\tensor{}M)
  \circ (\Delta\tensor{}M\tensor{}M) \circ \delta \\ &=&
  (\cc\tensor{}\varepsilon\tensor{}M) \circ \lr{ \lr{(\cc\tensor{}\xi) \circ (\Delta\tensor{}M)}
  \tensor{}M} \circ \delta \\ & \overset{\eqref{2-cell}}{=} &
  (\cc\tensor{}\varepsilon\tensor{}M) \circ \lr{ \lr{\Delta \circ \xi}
  \tensor{}M} \circ \delta \\ &=&
  (\cc\tensor{}\varepsilon\tensor{}M) \circ (\Delta\tensor{}M) \circ (\xi\tensor{}M)
   \circ \delta \\ &=& (\xi\tensor{}M)\circ \delta,
\end{eqnarray*}
which gives the commutativity of the first diagram in $(ii)$.
Multiplied on the left by $\mmm$, equation \eqref{counit-pr1}
becomes
\begin{eqnarray*}
  \mmm &=& (\mmm\tensor{}A) \circ (\cc\tensor{}M\tensor{}\varepsilon) \circ
  (\cc\tensor{}M\tensor{}\xi) \circ (\cc\tensor{}\mmm\tensor{}M)
  \circ (\Delta\tensor{}M\tensor{}M) \circ \delta \\
   &=& (M\tensor{}\cc\tensor{}\varepsilon) \circ (\mmm\tensor{}\cc) \circ
  (\cc\tensor{}M\tensor{}\xi) \circ (\cc\tensor{}\mmm\tensor{}M)
  \circ (\Delta\tensor{}M\tensor{}M) \circ \delta \\
   &=& (M\tensor{}\cc\tensor{}\varepsilon) \circ (M\tensor{}\cc\tensor{}\xi)
   \circ (\mmm\tensor{}\cc\tensor{}M)\circ (\cc\tensor{}\mmm\tensor{}M)
  \circ (\Delta\tensor{}M\tensor{}M) \circ \delta  \\
   &=& (M\tensor{}\cc\tensor{}\varepsilon) \circ (M\tensor{}\cc\tensor{}\xi)
   \circ \lr{\lr{(\mmm\tensor{}\cc) \circ (\cc\tensor{}\mmm)\circ (\Delta\tensor{}M) }\tensor{}M } \circ \delta \\
   & \overset{\eqref{1-cell}}{=}& (M\tensor{}\cc\tensor{}\varepsilon) \circ (M\tensor{}\cc\tensor{}\xi)
   \circ (M\tensor{}\Delta\tensor{}M) \circ (\mmm\tensor{}M) \circ \delta \\
   &=& (M\tensor{}\cc\tensor{}\varepsilon) \circ \lr{ M\tensor{} \lr{(\cc\tensor{}\xi)
   \circ (\Delta\tensor{}M) }} \circ (\mmm\tensor{}M) \circ \delta \\
   &\overset{\eqref{2-cell}}{=}& (M\tensor{}\cc\tensor{}\varepsilon) \circ (M\tensor{} \Delta)
   \circ (M\tensor{}\xi)  \circ (\mmm\tensor{}M) \circ \delta \\
   &=& (M\tensor{}\xi)  \circ (\mmm\tensor{}M) \circ \delta.
\end{eqnarray*}
Whence the commutativity of the second diagram of item $(ii)$.
Composing equation \eqref{coass-pr} with
$(\mmm\tensor{}M\tensor{}M)$, we obtain from the left hand term
\begin{multline*}
(\mmm\tensor{}M\tensor{}M) \circ
(\cc\tensor{}M\tensor{}M\tensor{}\varepsilon\tensor{}M) \circ
(\cc\tensor{}M\tensor{}\mmm\tensor{}M) \circ
(\cc\tensor{}\mmm\tensor{}M\tensor{}M) \circ
(\cc\tensor{}\delta\tensor{}M) \circ (\Delta\tensor{}M\tensor{}M)
\circ \delta  \\ \,=\,(\mmm\tensor{}M\tensor{}M) \circ
(\cc\tensor{}M\tensor{}M\tensor{}\varepsilon\tensor{}M) \circ
(\cc\tensor{}M\tensor{}\mmm\tensor{}M) \circ
(\cc\tensor{}\mmm\tensor{}M\tensor{}M) \circ
\lr{\lr{(\cc\tensor{}\delta) \circ (\Delta\tensor{}M)}\tensor{}M}
\circ \delta
\\ \overset{\eqref{2-cell}}{=}(\mmm\tensor{}M\tensor{}M) \circ
(\cc\tensor{}M\tensor{}M\tensor{}\varepsilon\tensor{}M) \circ
(\cc\tensor{}M\tensor{}\mmm\tensor{}M) \circ
(\cc\tensor{}\mmm\tensor{}M\tensor{}M) \circ
(\Delta\tensor{}M\tensor{}M\tensor{}M) \circ (\delta\tensor{}M)
\circ \delta
\\ \,=\,(M\tensor{}\cc\tensor{}M\tensor{}\varepsilon\tensor{}M) \circ
(\mmm\tensor{}M\tensor{}\cc\tensor{}M) \circ
(\cc\tensor{}M\tensor{}\mmm\tensor{}M) \circ
(\cc\tensor{}\mmm\tensor{}M\tensor{}M) \circ
(\Delta\tensor{}M\tensor{}M\tensor{}M) \circ (\delta\tensor{}M)
\circ \delta
\\ \,=\,(M\tensor{}\cc\tensor{}M\tensor{}\varepsilon\tensor{}M) \circ
(M\tensor{}\cc\tensor{}\mmm\tensor{}M) \circ
(\mmm\tensor{}\cc\tensor{}M\tensor{}M) \circ
(\cc\tensor{}\mmm\tensor{}M\tensor{}M) \circ
(\Delta\tensor{}M\tensor{}M\tensor{}M) \circ (\delta\tensor{}M)
\circ \delta
\\ \,=\,(M\tensor{}\cc\tensor{}M\tensor{}\varepsilon\tensor{}M) \circ
(M\tensor{}\cc\tensor{}\mmm\tensor{}M) \circ \lr{
\lr{(\mmm\tensor{}\cc) \circ (\cc\tensor{}\mmm) \circ
(\Delta\tensor{}M)}\tensor{}M\tensor{}M} \circ (\delta\tensor{}M)
\circ \delta
\\ \overset{\eqref{1-cell}}{=} (M\tensor{}\cc\tensor{}M\tensor{}\varepsilon\tensor{}M) \circ
(M\tensor{}\cc\tensor{}\mmm\tensor{}M) \circ
(M\tensor{}\Delta\tensor{}M\tensor{}M) \circ
(\mmm\tensor{}M\tensor{}M)  \circ (\delta\tensor{}M) \circ \delta
\\ \,=\, \lr{ M\tensor{}\cc\tensor{}\lr{ (M\tensor{}\varepsilon) \circ \mmm }\tensor{}M} \circ
(M\tensor{}\Delta\tensor{}M\tensor{}M) \circ
(\mmm\tensor{}M\tensor{}M)  \circ (\delta\tensor{}M) \circ \delta \\
\overset{\eqref{1-cell'}}{=}
(M\tensor{}\cc\tensor{}\varepsilon\tensor{}M\tensor{}M) \circ
(M\tensor{}\Delta\tensor{}M\tensor{}M) \circ
(\mmm\tensor{}M\tensor{}M)  \circ (\delta\tensor{}M) \circ \delta \\
\,=\, (\mmm\tensor{}M\tensor{}M)  \circ (\delta\tensor{}M) \circ
\delta.
\end{multline*}
From the right hand term, we get
\begin{multline*}
(\mmm\tensor{}M\tensor{}M) \circ
(\cc\tensor{}M\tensor{}\varepsilon\tensor{}M\tensor{}M) \circ
(\cc\tensor{}M\tensor{}\delta) \circ (\cc\tensor{}\mmm\tensor{}M)
\circ (\Delta\tensor{}M\tensor{}M) \circ \delta \\
\,=\, (M\tensor{}\cc\tensor{}\varepsilon\tensor{}M\tensor{}M) \circ
(\mmm\tensor{}\cc\tensor{}M\tensor{}M)\circ
(\cc\tensor{}M\tensor{}\delta) \circ (\cc\tensor{}\mmm\tensor{}M)
\circ (\Delta\tensor{}M\tensor{}M) \circ \delta \\
\,=\, (M\tensor{}\cc\tensor{}\varepsilon\tensor{}M\tensor{}M) \circ
(M\tensor{}\cc\tensor{}\delta) \circ (\mmm\tensor{}\cc\tensor{}M)
\circ (\cc\tensor{}\mmm\tensor{}M) \circ
(\Delta\tensor{}M\tensor{}M) \circ \delta \\ \,=\,
(M\tensor{}\cc\tensor{}\varepsilon\tensor{}M\tensor{}M) \circ
(M\tensor{}\cc\tensor{}\delta) \circ \lr{ \lr{(\mmm\tensor{}\cc)
\circ (\cc\tensor{}\mmm) \circ (\Delta\tensor{}M)}\tensor{}M} \circ
\delta \\ \overset{\eqref{1-cell}}{=}
(M\tensor{}\cc\tensor{}\varepsilon\tensor{}M\tensor{}M) \circ
(M\tensor{}\cc\tensor{}\delta) \circ (M\tensor{}\Delta\tensor{}M)
\circ (\mmm\tensor{}M) \circ \delta\\ \,=\,
(M\tensor{}\cc\tensor{}\varepsilon\tensor{}M\tensor{}M) \circ
\lr{M\tensor{}\lr{ (\cc\tensor{}\delta) \circ (\Delta\tensor{}M)}}
\circ (\mmm\tensor{}M) \circ \delta \\ \overset{\eqref{2-cell}}{=}
(M\tensor{}\cc\tensor{}\varepsilon\tensor{}M\tensor{}M) \circ
(M\tensor{}\Delta\tensor{}M\tensor{}M) \circ (M\tensor{}\delta)
\circ (\mmm\tensor{}M) \circ \delta \\ \,=\, (M\tensor{}\delta)
\circ (\mmm\tensor{}M) \circ \delta.
\end{multline*}
Therefore, $(\mmm\tensor{}M\tensor{}M)  \circ (\delta\tensor{}M)
\circ \delta\,=\,(M\tensor{}\delta) \circ (\mmm\tensor{}M) \circ
\delta$, which gives the commutativity of the last diagram of item
$(ii)$.\\ $(ii) \Rightarrow (i)$. We need to show that $\xi$ and
$\delta$ satisfy equations \eqref{counit-pr1}, \eqref{counit-pr2}
and \eqref{coass-pr}. By hypothesis $\xi$ and $\delta$ are
$\cc$--colinear maps satisfying
\begin{eqnarray}
  (\xi\tensor{}M) \circ \delta &=& \cc \tensor{}M \label{ii-1}\\
  (M\tensor{}\xi) \circ (\mmm\tensor{}M) \circ \delta &=& \mmm \label{ii-2}\\
  (M\tensor{}\delta) \circ(\mmm\tensor{}M) \circ \delta &=&
  (\mmm\tensor{}M\tensor{}M) \circ (\delta\tensor{}M) \circ \delta
  \label{ii-3}.
\end{eqnarray}
We then compute
\begin{multline*}
(\cc\tensor{}M \tensor{}\varepsilon) \circ
(\cc\tensor{}M\tensor{}\xi) \circ (\cc\tensor{}\mmm\tensor{}M) \circ
(\Delta\tensor{}M\tensor{}M) \circ \delta \\
\overset{\eqref{2-cell}}{=} (\cc\tensor{}M \tensor{}\varepsilon)
\circ (\cc\tensor{}M\tensor{}\xi) \circ (\cc\tensor{}\mmm\tensor{}M)
\circ (\cc\tensor{}\delta) \circ (\Delta\tensor{}M) \\ \,=\,
(\cc\tensor{}M \tensor{}\varepsilon) \circ \lr{
\cc\tensor{}\lr{(M\tensor{}\xi) \circ (\mmm\tensor{}M) \circ
\delta}} \circ (\Delta\tensor{}M) \\  \overset{\eqref{ii-2}}{=}
(\cc\tensor{}M \tensor{}\varepsilon) \circ (\cc\tensor{}\mmm) \circ
(\Delta\tensor{}M)  \overset{\eqref{1-cell'}}{=}
(\cc\tensor{}\varepsilon\tensor{}M) \circ (\Delta\tensor{}M) \,=\,
\cc\tensor{}M,
\end{multline*}
which gives equality \eqref{counit-pr1}. The equality
\eqref{counit-pr2} is obtained as follows:
\begin{eqnarray*}
    (\cc\tensor{}\varepsilon\tensor{}M) \circ
   (\cc\tensor{}\xi\tensor{}M) \circ (\Delta\tensor{}M\tensor{}M) \circ \delta
   &\overset{\eqref{2-cell}}{=}& (\cc\tensor{}\varepsilon\tensor{}M) \circ
   (\cc\tensor{}\xi\tensor{}M) \circ (\cc\tensor{}\delta) \circ (\Delta\tensor{}M)  \\
   &=& (\cc\tensor{}\varepsilon\tensor{}M) \circ
   \lr{\cc\tensor{}\lr{ (\xi\tensor{}M) \circ \delta}} \circ (\Delta\tensor{}M) \\
  &\overset{\eqref{ii-1}}{=}& (\cc\tensor{}\varepsilon\tensor{}M)
  \circ (\Delta\tensor{}M)\,\,=\,\, \cc\tensor{}M
\end{eqnarray*}
Lastly, the coassociativity that is equality \eqref{coass-pr} is
derived from the following computation:
\begin{multline*}
(\cc\tensor{}M\tensor{}M\tensor{}\varepsilon \tensor{}M) \circ
(\cc\tensor{}M\tensor{}\mmm\tensor{}M) \circ
(\cc\tensor{}\mmm\tensor{}M\tensor{}M) \circ
(\cc\tensor{}\delta\tensor{}M) \circ (\Delta\tensor{}M\tensor{}M)
\circ \delta \\ \overset{\eqref{2-cell}}{=}
(\cc\tensor{}M\tensor{}M\tensor{}\varepsilon \tensor{}M) \circ
(\cc\tensor{}M\tensor{}\mmm\tensor{}M) \circ
(\cc\tensor{}\mmm\tensor{}M\tensor{}M) \circ
(\cc\tensor{}\delta\tensor{}M) \circ (\cc\tensor{}\delta) \circ
(\Delta\tensor{}M) \\ \,=\,
(\cc\tensor{}M\tensor{}M\tensor{}\varepsilon \tensor{}M) \circ
(\cc\tensor{}M\tensor{}\mmm\tensor{}M) \circ
\lr{\cc\tensor{}\lr{(\mmm\tensor{}M\tensor{}M) \circ
(\delta\tensor{}M) \circ \delta}} \circ (\Delta\tensor{}M) \\
\overset{\eqref{ii-3}}{=}
(\cc\tensor{}M\tensor{}M\tensor{}\varepsilon \tensor{}M) \circ
(\cc\tensor{}M\tensor{}\mmm\tensor{}M) \circ
(\cc\tensor{}M\tensor{}\delta) \circ (\cc\tensor{}\mmm\tensor{}M)
\circ(\cc\tensor{}\delta) \circ (\Delta\tensor{}M) \\
\overset{\eqref{2-cell}}{=}
(\cc\tensor{}M\tensor{}M\tensor{}\varepsilon \tensor{}M) \circ
(\cc\tensor{}M\tensor{}\mmm\tensor{}M) \circ
(\cc\tensor{}M\tensor{}\delta) \circ (\cc\tensor{}\mmm\tensor{}M)
\circ (\Delta\tensor{}M\tensor{}M) \circ \delta \\  \,=\,
\lr{\cc\tensor{}M\tensor{}\lr{(M\tensor{} \varepsilon) \circ \mmm}
\tensor{} M} \circ (\cc\tensor{}M\tensor{}\delta) \circ
(\cc\tensor{}\mmm\tensor{}M) \circ (\Delta\tensor{}M\tensor{}M)
\circ \delta \\ \overset{\eqref{1-cell'}}{=} (\cc\tensor{}M\tensor{}
\varepsilon\tensor{}M\tensor{}M) \circ
(\cc\tensor{}M\tensor{}\delta) \circ (\cc\tensor{}\mmm\tensor{}M)
\circ (\Delta\tensor{}M\tensor{}M) \circ \delta.
\end{multline*}
\end{proof}

\begin{example}
Of course any $A$-coring $\cc$ can be seen as a cowreath over the
trivial coring $(A:A)$.\\
Let $C$ and $D$ two $\KK$-coalgebras. It is clear that $(D,\tau)$
belongs to $\rR_{(C:\,\KK)}$, where $\tau:C\tensor{\KK}D \to
D\tensor{\KK}C$ is the usual flip. Consider the maps
$\xi=C\tensor{}\varepsilon_D: C\tensor{\KK}D \to C$ and
$\delta=C\tensor{}\Delta_D:C\tensor{\KK}D \to
C\tensor{\KK}D\tensor{\KK}D$. One can easily prove that those maps
define in fact a morphisms in the category $\rR_{(C:\,\KK)}$; that
is $\xi: (D,\tau) \to (\KK,C)$ and $\delta: (D,\tau) \to
(D,\tau)\underset{(C:\,\KK)}{\otimes}(D,\tau)$. Moreover, $\xi$ and
$\delta$ respect the commutativity of the diagrams stated in
Proposition \ref{Equivalent-Def}(ii) with $(\cc: A)\,=\, (C:\KK)$.
Whence, $(D,\tau)$ is a $C$--cowreath.
\end{example}

\begin{example}
Let $(\cc:A)$ and $(\dd:A)$ two corings. Assume that there is an
$A$-bilinear map $\ddd: \cc\tensor{A}\dd \to \dd\tensor{A}\cc$
satisfying
\begin{eqnarray}
  (\dd\tensor{}\varepsilon_{\cc}) \circ \ddd &=& \varepsilon_{\cc}\tensor{}\dd \label{CD-1} \\
   (\dd \tensor{}\Delta_{\cc}) \circ \ddd &=& (\ddd\tensor{}\cc) \circ
   (\cc\tensor{}\dd) \circ (\Delta_{\cc}\tensor{}\dd) \label{CD-2}  \\
  (\varepsilon_{\dd}\tensor{}\cc) \circ \ddd &=& \cc\tensor{}\varepsilon_{\dd} \label{CD-3}  \\
  (\Delta_{\dd}\tensor{}\cc) \circ \dd &=& (\dd\tensor{}\dd) \circ
  (\ddd\tensor{}\dd) \circ (\cc\tensor{}\Delta_{\dd}). \label{CD-4}
\end{eqnarray}
Equations \eqref{CD-1} and \eqref{CD-2} say that $(\dd,\ddd)$ is
an object of the category $\rR_{(\cc:\,A)}$. While equations
\eqref{CD-3} and \eqref{CD-3} say that $(\ddd,\cc)$ is an object
of the category $\lL_{(\dd:\,A)}$ of Remark \ref{left-right}. One
can check that $(\dd,\ddd)$ is a right $\cc$-cowreath with
structure maps $\cc\tensor{}\varepsilon_{\dd}$ and
$\cc\tensor{}\Delta_{\dd}$, and similarly $(\ddd,\cc)$ is a left
$\dd$-cowreath with structure maps $\varepsilon_{\cc}\tensor{}\dd$
and $\Delta_{\cc}\tensor{}\dd$.
\end{example}

The following proposition gives, using an entwining structures, a
method to construct a cowreath from the commutative case to the non
commutative one. This proposition can be deduced from Lemma
\ref{ent-R} if we success to show that the functor occurring there
is monoidal. We prefer here to include a direct proof without
recalling more general theory.

\begin{proposition}\label{ent-corona}
Let $(A,C)_{\aaa}$ be an entwining structure over $\KK$ with
$\aaa: C\tensor{\KK}A \to A\tensor{\KK}C$. Consider its associated
coring $(A\tensor{\KK}C:A)$. If $(N,\nnn)$ is a $C$--cowreath,
then $(A\tensor{\KK}N\tensor{\KK}A,
(A\tensor{\KK}N\tensor{\KK}\aaa) \circ
(A\tensor{\KK}\nnn\tensor{\KK}A))$ is an
$(A\tensor{\KK}C)$--cowreath.
\end{proposition}
\begin{proof}
During this proof the tensor product $-\tensor{\KK}-$ will be
denoted by $-\tensor{}-$. Set $\cc:=A\tensor{}C$,
$M:=A\tensor{}N\tensor{}A$, and $\mmm:= (A\tensor{}N\tensor{}\aaa)
\circ (A\tensor{}\nnn\tensor{}A)$. By Lemma \ref{ent-R}, the pair
$(M,\mmm)$ belongs to $\rR_{(\cc:\,A)}$. Let us denote by
$\xi:C\tensor{}N \to C$ and $\delta: C\tensor{}N \to
C\tensor{}N\tensor{}N$ the structure of the $C$-cowreath $(N,\nnn)$.
Define the following $A$-bilinear maps
$$\xymatrix@C=60pt{ \td{\xi}: \cc\tensor{A}M \ar@{->}^-{A\tensor{}\xi\tensor{}A}[r] & A\tensor{}C\tensor{}A
\ar@{->}^-{A\tensor{}\aaa}[r] & A\tensor{}A\tensor{}C
\ar@{->}^-{\mu\tensor{}C}[r] & A\tensor{}C=\cc}$$ $$
\xymatrix@C=80pt{\td{\delta}: \cc\tensor{A}M
\ar@{->}^-{A\tensor{}\delta\tensor{}A}[r] &
A\tensor{}C\tensor{}N\tensor{}N\tensor{}A
\ar@{->}^-{A\tensor{}C\tensor{}N\tensor{}1\tensor{}N\tensor{}A}[r] &
\cc\tensor{A}M\tensor{A}M},$$ where we have used the isomorphisms
$\cc\tensor{A}M\cong A\tensor{}C\tensor{}N\tensor{}A$ and
$\cc\tensor{A}M\tensor{A}M\cong
A\tensor{}C\tensor{}N\tensor{}A\tensor{}N \tensor{}A$. We claim that
$(M,\mmm)$ is a $\cc$--cowreath with structure maps $\td{\xi}$ and
$\td{\delta}$. First, we need to show that those maps are
$\cc$--bicolinear. Starting with $\td{\xi}$, we know that
\begin{multline*}
(\cc\tensor{A}\td{\xi}) \circ (\Delta_{\cc}\tensor{A}M) \,=\,
\lr{(A\tensor{}C)\tensor{A}(\mu\tensor{}C)} \circ \lr{(A\tensor{}C)
\tensor{A}(A\tensor{}\aaa)} \circ
\lr{(A\tensor{}C)\tensor{A}(A\tensor{}\xi\tensor{}A)} \\ \,\,\,
\circ \lr{(A\tensor{}\Delta)\tensor{A}(A\tensor{}N\tensor{}A)}.
\end{multline*}
Since the the right $A$--action on $A\tensor{}C$ is given by the
twist map $\aaa$, the term
$(A\tensor{}C)\tensor{A}(\mu\tensor{}C)$ is identified with
$(\mu\tensor{}C\tensor{}C) \circ (A\tensor{}\aaa\tensor{}C)$,
which implies that
\begin{eqnarray*}
  (\cc\tensor{A}\td{\xi}) \circ (\Delta_{\cc}\tensor{A}M) &=&
  (\mu\tensor{}C\tensor{}C) \circ (A\tensor{}\aaa\tensor{}C) \circ
(A\tensor{}C\tensor{}\aaa) \circ
(A\tensor{}C\tensor{}\xi\tensor{}A) \circ
(A\tensor{}\Delta\tensor{}N\tensor{}A) \\
   &=& (\mu\tensor{}C\tensor{}C) \circ (A\tensor{}\aaa\tensor{}C)
\circ (A\tensor{}C\tensor{}\aaa) \circ
\lr{A\tensor{}\lr{(C\tensor{}\xi) \circ (\Delta\tensor{}N)}\tensor{}A} \\
   &\overset{\eqref{2-cell}}{=}& (\mu\tensor{}C\tensor{}C) \circ (A\tensor{}\aaa\tensor{}C)
\circ (A\tensor{}C\tensor{}\aaa) \circ
(A\tensor{}\Delta\tensor{}A)
\circ (A\tensor{}\xi\tensor{}A) \\
   &=& (\mu\tensor{}C\tensor{}C)
\circ \lr{A\tensor{}\lr{(\aaa\tensor{}C) \circ (C\tensor{}\aaa)
\circ (\Delta\tensor{}A)}} \circ (A\tensor{}\xi\tensor{}A) \\
   &\overset{\eqref{ent-20}}{=}& (\mu\tensor{}C\tensor{}C) \circ (A\tensor{}A\tensor{}\Delta)
\circ (A\tensor{}\aaa) \circ (A\tensor{}\xi\tensor{}A) \\
   &=& (A\tensor{}\Delta) \circ (\mu\tensor{}C) \circ (A\tensor{}\aaa)
\circ (A\tensor{}\xi\tensor{}A) \,\,=\,\, \Delta_{\cc} \circ
\td{\xi}.
\end{eqnarray*}
This proves that $\td{\xi}$ is left $\cc$--colinear. Its right
$\cc$--colinearity is obtained as follows: We know that
\begin{multline*}
(\td{\xi}\tensor{A}\cc) \circ (\cc\tensor{A}\mmm) \circ
(\Delta_{\cc}\tensor{A}M) \,=\,
\lr{(\mu\tensor{}C)\tensor{A}(A\tensor{}C)}  \circ
\lr{(A\tensor{}\aaa)\tensor{A}(A\tensor{}C)} \circ
\lr{(A\tensor{}\xi\tensor{}A)\tensor{A}(A\tensor{}C)} \\ \, \,
\circ\, \lr{(A\tensor{}C)\tensor{A}(A\tensor{}N\tensor{}\aaa)}
\circ \lr{(A\tensor{}C)\tensor{A}(A\tensor{}\nnn\tensor{}A)} \circ
\lr{(A\tensor{}\Delta)\tensor{A}(A\tensor{}N\tensor{}A}.
\end{multline*}
So using the obvious isomorphisms, we compute
\begin{eqnarray*}
  (\td{\xi}\tensor{A}\cc) \circ (\cc\tensor{A}\mmm) \circ
(\Delta_{\cc}\tensor{A}M) &=& (\mu\tensor{}C\tensor{}C) \circ
(A\tensor{}\aaa\tensor{}C) \circ
(A\tensor{}\xi\tensor{}A\tensor{}C) \circ
(A\tensor{}C\tensor{}N\tensor{}\aaa) \\
   &\,\,& \circ \,\, (A\tensor{}C\tensor{}\nnn\tensor{}A) \circ (A\tensor{}\Delta\tensor{}N\tensor{}A)  \\
   &=& (\mu\tensor{}C\tensor{}C) \circ (A\tensor{}\aaa\tensor{}C)
\circ (A\tensor{}C\tensor{}\aaa) \circ
(A\tensor{}\xi\tensor{}C\tensor{}A) \\ &\,\,&
\circ \,\, (A\tensor{}C\tensor{}\nnn\tensor{}A) \circ (A\tensor{}\Delta\tensor{}N\tensor{}A) \\
   &\overset{\eqref{2-cell'}}{=} & (\mu\tensor{}C\tensor{}C) \circ
(A\tensor{}\aaa\tensor{}C) \circ (A\tensor{}C\tensor{}\aaa)
\circ (A\tensor{}\Delta\tensor{}A) \circ (A\tensor{}\xi\tensor{}A)  \\
   &=& (\mu\tensor{}C\tensor{}C) \circ
\lr{A\tensor{}\lr{(\aaa\tensor{}C)
\circ (C\tensor{}\aaa) \circ (\Delta\tensor{}A)}} \circ (A\tensor{}\xi\tensor{}A) \\
   &\overset{\eqref{ent-20}}{=} & (\mu\tensor{}C\tensor{}C) \circ
(A\tensor{}A\tensor{}\Delta)
\circ (A\tensor{}\aaa) \circ (A\tensor{}\xi\tensor{}A)  \\
   &=& (A\tensor{}\Delta) \circ (\mu\tensor{}C) \circ
(A\tensor{}\aaa) \circ (A\tensor{}\xi\tensor{}A) \,\,=\,\,
\Delta_{\cc} \circ \td{\xi}.
\end{eqnarray*}
Concerning the colinearity of $\td{\delta}$, we have
\begin{eqnarray*}
   (\cc\tensor{A}\td{\delta}) \circ (\Delta_{\cc}\tensor{A}M) &=&
   \lr{(A\tensor{}C)\tensor{A}(A\tensor{}N\tensor{}1\tensor{}N\tensor{}A)}
   \circ \lr{(A\tensor{}C)\tensor{A}(A\tensor{}\delta\tensor{}A)}
   \\ &\,\,&\, \circ \lr{(A\tensor{}\Delta)\tensor{A}(A\tensor{}N\tensor{}A)} \\
   &=& (A\tensor{}C\tensor{}C\tensor{}N\tensor{}1\tensor{}N\tensor{}A)
   \circ (A\tensor{}C\tensor{}\delta\tensor{}A)
   \circ (A\tensor{}\Delta\tensor{}N\tensor{}A)  \\
   &=& (A\tensor{}C\tensor{}C\tensor{}N\tensor{}1\tensor{}N\tensor{}A)
   \circ \lr{A\tensor{}\lr{(C\tensor{}\delta)
   \circ (\Delta\tensor{}N)}\tensor{}A} \\
   &\overset{\eqref{2-cell}}{=}& (A\tensor{}C\tensor{}C\tensor{}N\tensor{}1\tensor{}N\tensor{}A)
   \circ (A\tensor{}\Delta\tensor{}N\tensor{}N\tensor{}A)
   \circ (A\tensor{}\delta\tensor{}A) \\
   &=& (A\tensor{}\Delta\tensor{}N\tensor{}A\tensor{}N\tensor{}A)
   \circ (A\tensor{}C\tensor{}N\tensor{}1\tensor{}N\tensor{}A)
   \circ (A\tensor{}\delta\tensor{}A) \\
   &=& (\Delta_{\cc}\tensor{A}M\tensor{A}M) \circ \td{\delta},
\end{eqnarray*}which implies that
 $\td{\delta}$ is left $\cc$--colinear. $\td{\delta}$ is right
$\cc$--colinear by the following computations: First we know that
\begin{multline*}
(\td{\delta}\tensor{A}\cc) \circ (\cc\tensor{A}\mmm) \circ
(\Delta_{\cc}\tensor{A}M) \,=\,
\lr{(A\tensor{}C\tensor{}N\tensor{}1\tensor{}N\tensor{}A)\tensor{A}(A\tensor{}C)}
\circ \lr{(A\tensor{}\delta\tensor{}A) \tensor{A}(A\tensor{}C)} \\
\,\, \circ\,
\lr{(A\tensor{}C)\tensor{A}(A\tensor{}N\tensor{}\aaa)} \circ
\lr{(A\tensor{}C)\tensor{A}(A\tensor{}\nnn\tensor{}A)} \circ
\lr{(A\tensor{}\Delta)\tensor{A}(A\tensor{}N\tensor{}A)}.
\end{multline*}
Using as before the obvious isomorphisms, we compute
\begin{multline*}
(\td{\delta}\tensor{A}\cc) \circ (\cc\tensor{A}\mmm) \circ
(\Delta_{\cc}\tensor{A}M) \,=\,
(A\tensor{}C\tensor{}N\tensor{}1\tensor{}N\tensor{}A\tensor{}C)
\circ (A\tensor{}\delta\tensor{}A\tensor{}C) \circ
(A\tensor{}C\tensor{}N\tensor{}\aaa) \\ \,\,  \,\circ\,\,
(A\tensor{}C\tensor{}\nnn\tensor{}A) \circ
(A\tensor{}\Delta\tensor{}N\tensor{}A) \\
\,=\,
(A\tensor{}C\tensor{}N\tensor{}1\tensor{}N\tensor{}A\tensor{}C)
\circ (A\tensor{}C\tensor{}N\tensor{}N\tensor{}\aaa) \circ
(A\tensor{}\delta\tensor{}C\tensor{}A)  \circ
(A\tensor{}C\tensor{}\nnn\tensor{}A) \circ
(A\tensor{}\Delta\tensor{}N\tensor{}A) \\
\,=\,
(A\tensor{}C\tensor{}N\tensor{}1\tensor{}N\tensor{}A\tensor{}C)
\circ (A\tensor{}C\tensor{}N\tensor{}N\tensor{}\aaa) \circ
\lr{A\tensor{}\lr{(\delta\tensor{}C)  \circ (C\tensor{}\nnn) \circ
(\Delta\tensor{}N)}\tensor{}A}\\
\overset{\eqref{2-cell'}}{=}
(A\tensor{}C\tensor{}N\tensor{}1\tensor{}N\tensor{}A\tensor{}C)
\circ (A\tensor{}C\tensor{}N\tensor{}N\tensor{}\aaa) \circ
(A\tensor{}C\tensor{}N\tensor{}\nnn\tensor{}A)  \circ
(A\tensor{}C\tensor{}\nnn\tensor{}N\tensor{}A) \\ \,\, \circ
(A\tensor{}\Delta\tensor{}N\tensor{}N\tensor{}A) \circ
(A\tensor{}\delta\tensor{}A) \\ \,=\,
(A\tensor{}C\tensor{}N\tensor{}A\tensor{}N\tensor{}\aaa) \circ
(A\tensor{}C\tensor{}N\tensor{}1\tensor{}N\tensor{}C\tensor{}A)
\circ (A\tensor{}C\tensor{}N\tensor{}\nnn\tensor{}A)  \circ
(A\tensor{}C\tensor{}\nnn\tensor{}N\tensor{}A) \\ \,\, \circ
(A\tensor{}\Delta\tensor{}N\tensor{}N\tensor{}A) \circ
(A\tensor{}\delta\tensor{}A) \\ \,=\,
(A\tensor{}C\tensor{}N\tensor{}A\tensor{}N\tensor{}\aaa) \circ
(A\tensor{}C\tensor{}N\tensor{}A\tensor{}\nnn\tensor{}A) \circ
(A\tensor{}C\tensor{}N\tensor{}1\tensor{}C \tensor{}N\tensor{}A)
\circ (A\tensor{}C\tensor{}\nnn\tensor{}N\tensor{}A) \\ \,\, \circ
(A\tensor{}\Delta\tensor{}N\tensor{}N\tensor{}A) \circ
(A\tensor{}\delta\tensor{}A),
\end{multline*}
by equation \eqref{ent-11}, we then get
$$ (\td{\delta}\tensor{A}\cc) \circ (\cc\tensor{A}\mmm) \circ
(\Delta_{\cc}\tensor{A}M) \,=\, (\cc\tensor{A}M\tensor{A}\mmm) \circ
(\cc\tensor{A}\mmm\tensor{A}M) \circ
(\Delta_{\cc}\tensor{A}M\tensor{A}M) \circ \td{\delta}.$$ To finish
the proof we need to check that $\td{\xi}$ and $\td{\delta}$ satisfy
the commutativity of the sated  diagrams in Proposition
\ref{Equivalent-Def}(ii). By assumptions, we know that $\xi$ and
$\delta$ satisfy the same property with base coring the coalgebra
$(C:\KK)$. So, we have
\begin{multline*}
(\td{\xi}\tensor{A}M) \circ \td{\delta} \,=\,
\lr{\lr{(\mu\tensor{}C) \circ (A\tensor{}\aaa) \circ
(A\tensor{}\xi\tensor{}A)}\tensor{A}(A\tensor{}N\tensor{}A)} \circ
(A\tensor{}C\tensor{}N\tensor{}1\tensor{}N\tensor{}A) \circ
(A\tensor{}\delta\tensor{}A) \\ \,=\,
(\mu\tensor{}C\tensor{}N\tensor{}A) \circ
(A\tensor{}\aaa\tensor{}N\tensor{}A) \circ
(A\tensor{}\xi\tensor{}A\tensor{}N\tensor{}A) \circ
(A\tensor{}C\tensor{}N\tensor{}1\tensor{}N\tensor{}A) \circ
(A\tensor{}\delta\tensor{}A) \\ \,=\,
(\mu\tensor{}C\tensor{}N\tensor{}A) \circ
(A\tensor{}\aaa\tensor{}N\tensor{}A) \circ
(A\tensor{}C\tensor{}1\tensor{}N\tensor{}A) \circ
(A\tensor{}\xi\tensor{}N\tensor{}A) \circ
(A\tensor{}\delta\tensor{}A) \\ \overset{\eqref{ent-11}}{=}
(\mu\tensor{}C\tensor{}N\tensor{}A) \circ
(A\tensor{}1\tensor{}C\tensor{}N\tensor{}A) \,=\,
A\tensor{}C\tensor{}N\tensor{}A \,=\, \cc\tensor{A}M,
\end{multline*}
and
\begin{multline*}
(M\tensor{A}\td{\xi}) \circ (\mmm\tensor{A}M) \circ \td{\delta}
\,=\, \lr{(A\tensor{}N\tensor{}A)\tensor{A}\lr{(\mu\tensor{}C) \circ
(A\tensor{}\aaa) \circ (A\tensor{}\xi\tensor{}A)}} \\
\,\,  \circ\, \lr{\lr{(A\tensor{}N\tensor{}\aaa) \circ
(A\tensor{}\nnn\tensor{}A)}\tensor{A}(A\tensor{}N\tensor{}A)}  \circ
(A\tensor{}C\tensor{}N\tensor{}1\tensor{}N\tensor{}A) \circ
(A\tensor{}\delta\tensor{}A) \\ \,=\,
(A\tensor{}N\tensor{}\mu\tensor{}C) \circ
(A\tensor{}N\tensor{}A\tensor{}\aaa) \circ
(A\tensor{}N\tensor{}A\tensor{}\xi\tensor{}A) \circ
(A\tensor{}N\tensor{}\aaa\tensor{}N\tensor{}A)  \circ
(A\tensor{}\nnn\tensor{}A\tensor{}N\tensor{}A) \\ \,\, \circ
(A\tensor{}C\tensor{}N\tensor{}1\tensor{}N\tensor{}A) \circ
(A\tensor{}\delta\tensor{}A) \\ \,=\,
(A\tensor{}N\tensor{}\mu\tensor{}C) \circ
(A\tensor{}N\tensor{}A\tensor{}\aaa) \circ
(A\tensor{}N\tensor{}A\tensor{}\xi\tensor{}A) \circ
(A\tensor{}N\tensor{}\aaa\tensor{}N\tensor{}A)  \circ
(A\tensor{}N\tensor{}C\tensor{}1\tensor{}N\tensor{}A) \\
\,\, \circ (A\tensor{}\nnn\tensor{}N\tensor{}A) \circ
(A\tensor{}\delta\tensor{}A) \\ \overset{\eqref{ent-11}}{=}
(A\tensor{}N\tensor{}\mu\tensor{}C) \circ
(A\tensor{}N\tensor{}A\tensor{}\aaa) \circ
(A\tensor{}N\tensor{}A\tensor{}\xi\tensor{}A) \circ
(A\tensor{}N\tensor{}1\tensor{}C\tensor{}N\tensor{}A)  \circ
\lr{A\tensor{}\lr{(\nnn\tensor{}N) \circ \delta}\tensor{}A} \\ \,=\,
(A\tensor{}N\tensor{}\mu\tensor{}C) \circ
(A\tensor{}N\tensor{}A\tensor{}\aaa) \circ
(A\tensor{}N\tensor{}1\tensor{}C\tensor{}A) \circ
(A\tensor{}N\tensor{}\xi\tensor{}A)  \circ
\lr{A\tensor{}\lr{(\nnn\tensor{}N) \circ \delta}\tensor{}A} \\ \,=\,
(A\tensor{}N\tensor{}\mu\tensor{}C) \circ
(A\tensor{}N\tensor{}A\tensor{}\aaa) \circ
(A\tensor{}N\tensor{}1\tensor{}C\tensor{}A) \circ
\lr{A\tensor{}\lr{(N\tensor{}\xi) \circ (\nnn\tensor{}N) \circ
\delta}\tensor{}A} \\ \,=\, (A\tensor{}N\tensor{}\mu\tensor{}C)
\circ (A\tensor{}N\tensor{}1\tensor{}A\tensor{}C) \circ
(A\tensor{}N\tensor{}\aaa) \circ (A\tensor{}\nnn\tensor{}A) \\ \,=\,
(A\tensor{}N\tensor{}\aaa) \circ (A\tensor{}\nnn\tensor{}A)\,\,=\,\,
\mmm,
\end{multline*}
These give the commutativity of the two first diagrams in
Proposition \eqref{Equivalent-Def}(ii). The commutativity of the
third one is derived from the following computations: From one hand,
we have
\begin{multline*}
(\mmm\tensor{A}M\tensor{A}M) \circ (\td{\delta}\tensor{A}M) \circ
\td{\delta}  \,=\,
(A\tensor{}N\tensor{}\aaa\tensor{}N\tensor{}A\tensor{}N\tensor{}A)
\circ
(A\tensor{}\nnn\tensor{}A\tensor{}N\tensor{}A\tensor{}N\tensor{}A)
\\ \,\, \circ
(A\tensor{}C\tensor{}N\tensor{}1\tensor{}N\tensor{}A\tensor{}N\tensor{}A)
\circ (A\tensor{}\delta\tensor{}A\tensor{}N\tensor{}A)
 \circ (A\tensor{}C\tensor{}N\tensor{}1\tensor{}N\tensor{}A)
\circ (A\tensor{}\delta\tensor{}A)  \\ \,=\,
(A\tensor{}N\tensor{}\aaa\tensor{}N\tensor{}A\tensor{}N\tensor{}A)
\circ
(A\tensor{}N\tensor{}C\tensor{}1\tensor{}N\tensor{}A\tensor{}N\tensor{}A)
\circ (A\tensor{}\nnn\tensor{}N\tensor{}A\tensor{}N\tensor{}A)
 \circ (A\tensor{}\delta\tensor{}A\tensor{}N\tensor{}A)
\\ \,\, \circ (A\tensor{}C\tensor{}N\tensor{}1\tensor{}N\tensor{}A)
\circ (A\tensor{}\delta\tensor{}A) \\ \overset{\eqref{ent-11}}{=}
(A\tensor{}N\tensor{}1\tensor{}C\tensor{}N\tensor{}A\tensor{}N\tensor{}A)
\circ (A\tensor{}\nnn\tensor{}N\tensor{}A\tensor{}N\tensor{}A) \circ
(A\tensor{}\delta\tensor{}A\tensor{}N\tensor{}A) \circ
(A\tensor{}C\tensor{}N\tensor{}1\tensor{}N\tensor{}A)  \circ
(A\tensor{}\delta\tensor{}A)\\ \,=\,
(A\tensor{}N\tensor{}1\tensor{}C\tensor{}N\tensor{}A\tensor{}N\tensor{}A)
\circ (A\tensor{}\nnn\tensor{}N\tensor{}A\tensor{}N\tensor{}A) \circ
(A\tensor{}C\tensor{}N\tensor{}N\tensor{}1\tensor{}N\tensor{}A)
\circ (A\tensor{}\delta\tensor{}N\tensor{}A) \circ
(A\tensor{}\delta\tensor{}A)
\\ \,=\,
(A\tensor{}N\tensor{}1\tensor{}C\tensor{}N\tensor{}A\tensor{}N\tensor{}A)
\circ
(A\tensor{}N\tensor{}C\tensor{}N\tensor{}1\tensor{}N\tensor{}A)
\circ (A\tensor{}\nnn\tensor{}N\tensor{}N\tensor{}A) \circ
(A\tensor{}\delta\tensor{}N\tensor{}A) \circ
(A\tensor{}\delta\tensor{}A) \\ \,=\,
(A\tensor{}N\tensor{}1\tensor{}C\tensor{}N\tensor{}A\tensor{}N\tensor{}A)
\circ
(A\tensor{}N\tensor{}C\tensor{}N\tensor{}1\tensor{}N\tensor{}A)
\circ \lr{A\tensor{}\lr{(\nnn\tensor{}N\tensor{}N)
 \circ (\delta\tensor{}N)
\circ \delta}\tensor{}A} \\ \,=\,
(A\tensor{}N\tensor{}1\tensor{}C\tensor{}N\tensor{}A\tensor{}N\tensor{}A)
\circ
(A\tensor{}N\tensor{}C\tensor{}N\tensor{}1\tensor{}N\tensor{}A)
\circ (A\tensor{}N\tensor{}\delta\tensor{}A) \circ
(A\tensor{}\nnn\tensor{}N\tensor{}A) \circ
(A\tensor{}\delta\tensor{}A),
\end{multline*}
and from an other hand, we obtain
\begin{multline*}
(M\tensor{A}\td{\delta}) \circ (\mmm\tensor{A}M) \circ \td{\delta}
\,=\,
(A\tensor{}N\tensor{}C\tensor{}N\tensor{}1\tensor{}N\tensor{}A)
\circ (A\tensor{}N\tensor{}A\tensor{}\delta\tensor{}A) \circ
(A\tensor{}N\tensor{}\aaa\tensor{}N\tensor{}A) \circ
(A\tensor{}\nnn\tensor{}A\tensor{}N\tensor{}A) \\ \,\, \circ
(A\tensor{}C\tensor{}N\tensor{}1\tensor{}N\tensor{}A) \circ
(A\tensor{}\delta\tensor{}A) \\ \,=\,
(A\tensor{}N\tensor{}C\tensor{}N\tensor{}1\tensor{}N\tensor{}A)
\circ (A\tensor{}N\tensor{}A\tensor{}\delta\tensor{}A) \circ
(A\tensor{}N\tensor{}\aaa\tensor{}N\tensor{}A) \circ
(A\tensor{}N\tensor{}C\tensor{}1\tensor{}N\tensor{}A) \\ \,\, \circ
(A\tensor{}\nnn\tensor{}N\tensor{}A) \circ
(A\tensor{}\delta\tensor{}A) \\ \overset{\eqref{ent-11}}{=}
(A\tensor{}N\tensor{}C\tensor{}N\tensor{}1\tensor{}N\tensor{}A)
\circ (A\tensor{}N\tensor{}A\tensor{}\delta\tensor{}A) \circ
(A\tensor{}N\tensor{}1\tensor{}C\tensor{}N\tensor{}A)  \circ
(A\tensor{}\nnn\tensor{}N\tensor{}A) \circ
(A\tensor{}\delta\tensor{}A) \\ \,=\,
(A\tensor{}N\tensor{}C\tensor{}N\tensor{}1\tensor{}N\tensor{}A)
\circ
(A\tensor{}N\tensor{}1\tensor{}C\tensor{}N\tensor{}N\tensor{}A)
\circ (A\tensor{}N\tensor{}\delta\tensor{}A) \circ
(A\tensor{}\nnn\tensor{}N\tensor{}A) \circ
(A\tensor{}\delta\tensor{}A) \\ \,=\,
(A\tensor{}N\tensor{}1\tensor{}C\tensor{}N\tensor{}A\tensor{}N\tensor{}A)
\circ
(A\tensor{}N\tensor{}C\tensor{}N\tensor{}1\tensor{}N\tensor{}A)
\circ (A\tensor{}N\tensor{}\delta\tensor{}A) \circ
(A\tensor{}\nnn\tensor{}N\tensor{}A) \circ
(A\tensor{}\delta\tensor{}A)
\end{multline*}
Therefore, $(M\tensor{A}\td{\delta}) \circ (\mmm\tensor{A}M) \circ
\td{\delta}\,=\, (\mmm\tensor{A}M\tensor{A}M) \circ
(\td{\delta}\tensor{A}M) \circ \td{\delta}$.
\end{proof}

\begin{remark}
Entwining structures give an example of a wreath over coalgebras.
Explicitly, given any entwining structure $(A,C)_{\aaa}$ over $\KK$
with $\aaa: C\tensor{\KK}A \to A\tensor{\KK}C$. As we have already
observe $(A,\aaa)$ is an object of the monoidal category
$\rR_{(C:\,\KK)}$. Taking $\boldsymbol{\eta}=C\tensor{\KK}1: C \to
C\tensor{\KK}A$ and $\boldsymbol{\mu}=C\tensor{\KK}\mu:
C\tensor{\KK}A\tensor{\KK}A \to C\tensor{\KK}A$. One can easily
checks that $\boldsymbol{\eta}$ and $\boldsymbol{\mu}$ are in fact a
$C$--bicolinear morphisms, that is $\boldsymbol{\eta}: (\KK,C) \to
(A,\aaa)$ and $\boldsymbol{\mu}: (A,\aaa)
\underset{(C:\,\KK)}{\otimes} (A,\aaa) \to (A,\aaa)$ are morphisms
in $\rR_{(C:\,\KK)}$. Furthermore, $\boldsymbol{\eta}$ and
$\boldsymbol{\mu}$ endow $(A,\aaa)$ with a structure of an algebra
in the monoidal category $\rR_{(C:\,\KK)}$. That is $(A,\aaa)$ is in
our terminology a $C$--wreath.
\end{remark}

\section{The category of comodules over cowreath, cowreath products and
functors.}\label{Sect-Mod}

This section presents a simplest and equivalent definitions of the
objects and morphisms of the category of (right) comodules over a
cowreath. The definition of cowreath product and the construction of
some functors connecting categories are presented as well. For sake
of completeness a detailed proofs are included.

Fix a coring $(\cc:A)$ and let us use the symbol $-\tensor{}-$ to
denote the tensor product $-\tensor{A}-$ over the base ring $A$. Let
$(M,\mmm)$ be a $\cc$--cowreath with structure maps $\xi:
\cc\tensor{}M \to \cc$ and $\delta: \cc\tensor{}M \to
\cc\tensor{}M\tensor{}M$. Since $(M,\mmm)$ is a coalgebra in the
monoidal category $\Mono{C}$, it is natural to ask for the category
of (right) $(M,\mmm)$--comodules. Thus, an object $(X,\xxx)$ of
$\Mono{C}$ is said to be a \emph{right $(M,\mmm)$--comodule} if
there exists a morphisms $\varrho^{(X,\xxx)}: (X,\xxx) \to
(X,\xxx)\Rtensor{C} (M,\mmm)$ in $\Mono{C}$ which satisfies
\begin{equation}\label{comod-cow}
\lr{(X,\xxx)\Rtensor{C}\xi} \circ \varrho^{(X,\xxx)} \,=\, (X,\xxx),
\quad \lr{(X,\xxx)\Rtensor{C}\delta} \circ \varrho^{(X,\xxx)} \,=\,
\lr{\varrho^{(X,\xxx)}\Rtensor{C}(M,\mmm)} \circ \varrho^{(X,\xxx)}.
\end{equation}

\begin{proposition}\label{comod}
Let $(M,\mmm)$ be a $\cc$--cowreath with structure maps $\xi:
\cc\tensor{A}M \to \cc$ and $\delta: \cc\tensor{A}M \to
\cc\tensor{A}M\tensor{A}M$.
\begin{enumerate}[(a)]
\item Consider $(X,\xxx)$ an object of $\Mono{C}$. The following conditions are
equivalent
\begin{enumerate}[(i)]
\item $(X,\xxx)$ is right $(M,\mmm)$--comodule;

\item There is a $\cc$--bicomodules morphism
$\varrho^{(X,\xxx)}: \cc\tensor{A}X \to \cc\tensor{A}X \tensor{A}M$
such that $$\xymatrix@R=30pt{\cc\tensor{}X
\ar@{->}^-{\varrho^{(X,\xxx)}}[r]  \ar@{->}_-{\xxx}[rrd] &
\cc\tensor{}X\tensor{}M \ar@{->}^-{\xxx\tensor{}M}[r] &
X\tensor{}\cc\tensor{}M \ar@{->}^-{X\tensor{}\xi}[d]\\ & &
X\tensor{}\cc,} \,\,\, \xymatrix@C=40pt@R=30pt{\cc\tensor{}X
\ar@{->}^-{\varrho^{(X,\xxx)}}[r] \ar@{->}_-{\varrho^{(X,\xxx)}}[d]
& \cc\tensor{}X\tensor{}M \ar@{->}^-{\xxx\tensor{}M}[r] &
X\tensor{}\cc\tensor{}M
\ar@{->}^-{X\tensor{}\delta}[d]  \\
\cc\tensor{}X\tensor{}M \ar@{->}^-{\varrho^{(X,\xxx)}\tensor{}M}[r]
& \cc\tensor{}X\tensor{}M\tensor{}M
\ar@{->}^-{\xxx\tensor{}M\tensor{}M}[r]&
X\tensor{}\cc\tensor{}M\tensor{}M }$$ are commutative diagrams.
\end{enumerate}

\item Given two right $(M,\mmm)$--comodules $(X,\xxx)$ and $(X',\xxx')$.
A morphism $\Sf{f}: (X,\xxx) \to (X',\xxx')$ in $\Mono{C}$ is a
morphism of right $(M,\mmm)$--comodules if and only if $\,\Sf{f}$
turns commutative the following diagram
$$\xymatrix@C=50pt@R=30pt{ \cc\tensor{}X
\ar@{->}^-{\Sf{f}}[r] \ar@{->}_-{\varrho^{(X,\xxx)}}[d] & \cc\tensor{}X' \ar@{->}^-{\varrho^{(X',\xxx')}}[d] \\
\cc\tensor{}X\tensor{}M \ar@{->}^-{\Sf{f}\tensor{}M}[r] &
\cc\tensor{}X'\tensor{}M}$$
\end{enumerate}
\end{proposition}
\begin{proof}
$(a)$. The proof of Proposition \ref{Equivalent-Def} can be adapted
to the setting of this item.\\ $(b)$. The map $\Sf{f}$ is a morphism
of right $(M,\mmm)$--comodule if and only if
$$\varrho^{(X,\xxx)} \circ \Sf{f} \,=\,
\lr{\Sf{f}\Rtensor{C}(M,\mmm)} \circ \varrho^{(X',\xxx')},$$ if and
only if
\begin{eqnarray*}
  \varrho^{(X,\xxx)} \circ \Sf{f} &=& (\cc\tensor{}X'\tensor{}\varepsilon\tensor{}M)
  \circ (\cc\tensor{}\xxx'\tensor{}M) \circ (\cc\tensor{}\Sf{f}\tensor{}M)
  \circ (\Delta\tensor{}X\tensor{}M) \circ \varrho^{(X',\xxx')} \\
   &\overset{\eqref{1-cell'}}{=}& (\cc\tensor{}\varepsilon\tensor{}X'\tensor{}M)
  \circ (\cc\tensor{}\Sf{f}\tensor{}M)
  \circ (\Delta\tensor{}X\tensor{}M) \circ \varrho^{(X',\xxx')} \\
   &=& (\cc\tensor{}\varepsilon\tensor{}X'\tensor{}M)
  \circ \lr{\lr{(\cc\tensor{}\Sf{f})
  \circ (\Delta\tensor{}X)}\tensor{}M} \circ \varrho^{(X',\xxx')} \\
   &\overset{\eqref{2-cell}}{=}& (\cc\tensor{}\varepsilon\tensor{}X'\tensor{}M)
  \circ (\Delta\tensor{}X'\tensor{}M) \circ (\Sf{f}\tensor{}M) \circ
  \varrho^{(X',\xxx')} \\ &=& (\Sf{f}\tensor{}M) \circ
  \varrho^{(X',\xxx')}.
\end{eqnarray*}
\end{proof}
Clearly $(M,\mmm)$ is right $(M,\mmm)$-comodule with coaction
$\rho^{(M,\mmm)}\,=\, \delta$, and $(M\tensor{}M,(M\tensor{}\mmm)
\circ (\mmm\tensor{}M))$ is right $(M,\mmm)$--comodule with coaction
$$\rho^{\lr{M\tensor{}M,(M\tensor{}\mmm) \circ
(\mmm\tensor{}M)}}\,\,=\,\,
(\cc\tensor{}\varepsilon\tensor{}M\tensor{}M) \circ
(\cc\tensor{}M\tensor{}\delta) \circ (\cc\tensor{}\mmm\tensor{}M)
\circ (\Delta\tensor{}M\tensor{}M).$$

The expressions of objects and morphisms in the category of left
$(M,\mmm)$--comodules are given by the following

\begin{proposition}\label{comod1}
Let $(M,\mmm)$ be a $\cc$--cowreath with structure maps $\xi:
\cc\tensor{A}M \to \cc$ and $\delta: \cc\tensor{A}M \to
\cc\tensor{A}M\tensor{A}M$.
\begin{enumerate}[(a)]
\item Consider $(X,\xxx)$ an object of $\Mono{C}$. The following
conditions are equivalent
\begin{enumerate}[(i)]
\item $(X,\xxx)$ is left $(M,\mmm)$--comodule;

\item There is a $\cc$--bicomodules morphism
$\lambda^{(X,\xxx)}: \cc\tensor{A}X \to \cc\tensor{A}M
\tensor{A}X$ such that
$$\xymatrix@C=40pt@R=30pt{\cc\tensor{}X \ar@{->}^-{\lambda^{(X,\xxx)}}[r]
\ar@{=}[dr] & \cc\tensor{}M\tensor{}X \ar@{->}^-{\xi\tensor{}X}[d] \\
& \cc\tensor{}X ,} \,\,\, \xymatrix@C=40pt@R=30pt{\cc\tensor{}X
\ar@{->}^-{\lambda^{(X,\xxx)}}[r] \ar@{->}_-{\lambda^{(X,\xxx)}}[d]
& \cc\tensor{}M\tensor{}X \ar@{->}^-{\delta\tensor{}X}[r]
&\cc\tensor{}M\tensor{}M\tensor{}X
\ar@{->}^-{\mmm\tensor{}M\tensor{}X}[d]  \\
\cc\tensor{}M\tensor{}X \ar@{->}^-{\mmm\tensor{}X}[r] &
M\tensor{}\cc\tensor{}X
\ar@{->}^-{M\tensor{}\lambda^{(X,\,\xxx)}}[r]&
M\tensor{}\cc\tensor{}M\tensor{}X }$$ are commutative diagrams.
\end{enumerate}

\item Given two left $(M,\mmm)$--comodules $(X,\xxx)$ and
$(X',\xxx')$. A morphism $\Sf{f}: (X,\xxx) \to (X',\xxx')$ in
$\Mono{C}$ is a morphism of right $(M,\mmm)$--comodules if and
only if $\,\Sf{f}$ turns commutative the following diagram
$$\xymatrix@C=50pt@R=20pt{ & \cc\tensor{}X
\ar@{->}^-{\Sf{f}}[r] \ar@{->}_-{\lambda^{(X,\xxx)}}[dl] &
\cc\tensor{}X' \ar@{->}^-{\lambda^{(X',\xxx')}}[dr] & \\
\cc\tensor{}M\tensor{}X \ar@{->}_-{\mmm\tensor{}X}[rd] & & & \cc\tensor{}M\tensor{}X' \ar@{->}^-{\mmm\tensor{}X'}[ld] \\
& M\tensor{}\cc\tensor{}X \ar@{->}^-{M\tensor{}\Sf{f}}[r] &
M\tensor{}\cc\tensor{}X' & }$$
\end{enumerate}
\end{proposition}

The cowreath product, as was defined in more general setting in
\cite{Lack/Street:2002}, is given in our case by the following.

\begin{proposition}\label{cow-prod}
Let $(\cc:A)$ be any coring and $(M,\mmm)$ a $\cc$--cowreath with
structure maps $\xi:\cc\tensor{A}M \to \cc$ and
$\delta:\cc\tensor{A}M \to \cc\tensor{A}M\tensor{A}M$. Then
$(\cc\tensor{A}M:A)$ is a coring whose comultiplication and counit
are defined by
$$ \Delta' \,\,=\,\, (\cc\tensor{A}\mmm\tensor{A}M) \circ (\cc\tensor{A}\delta) \circ
(\Delta\tensor{A}M), \qquad \varepsilon' \,\,=\,\, \varepsilon
\circ \xi.$$ Moreover, $\xi: \cc\tensor{A}M \to \cc$ is a morphism
of $A$--corings.
\end{proposition}
\begin{proof}
By definition $\varepsilon'$ and $\Delta'$ are $A$--bilinear maps.
We need to show the counitary and the coassociativity properties.
The counitary property is derived from the following two
computations
\begin{eqnarray*}
  (\cc\tensor{}M\tensor{}\varepsilon') \circ \Delta' &=& (\cc\tensor{}M\tensor{}\varepsilon)
   \circ (\cc\tensor{}M\tensor{}\xi) \circ (\cc\tensor{}\mmm\tensor{}M) \circ (\cc\tensor{}\delta) \circ
(\Delta\tensor{}M) \\
   &=& (\cc\tensor{}M\tensor{}\varepsilon)
   \circ \lr{\cc\tensor{}\lr{(M\tensor{}\xi) \circ (\mmm\tensor{}M) \circ \delta}} \circ
(\Delta\tensor{}M) \\
   &\overset{\ref{Equivalent-Def}(ii)}{=}& (\cc\tensor{}M\tensor{}\varepsilon)
   \circ (\cc\tensor{}\mmm) \circ (\Delta\tensor{}M) \\
   &\overset{\eqref{1-cell'}}{=}& (\cc\tensor{}\varepsilon\tensor{}M)
   \circ (\Delta\tensor{}M) \,\,=\,\, \cc\tensor{}M ,
\end{eqnarray*}
and
\begin{eqnarray*}
  (\varepsilon'\tensor{}\cc\tensor{}M) \circ \Delta' &=& (\varepsilon\tensor{}\cc\tensor{}M)
  \circ (\xi\tensor{}\cc\tensor{}M) \circ (\cc\tensor{}\mmm\tensor{}M) \circ (\cc\tensor{}\delta) \circ
(\Delta\tensor{}M) \\
   &\overset{\eqref{2-cell}}{=}& (\varepsilon\tensor{}\cc\tensor{}M)
  \circ (\xi\tensor{}\cc\tensor{}M) \circ (\cc\tensor{}\mmm\tensor{}M)
  \circ (\Delta\tensor{}M\tensor{}M) \circ \delta \\
   &=& (\varepsilon\tensor{}\cc\tensor{}M)
  \circ \lr{\lr{(\xi\tensor{}\cc) \circ (\cc\tensor{}\mmm)
  \circ (\Delta\tensor{}M)}\tensor{}M} \circ \delta \\
   &\overset{\eqref{2-cell'}}{=}& (\varepsilon\tensor{}\cc\tensor{}M)
  \circ (\Delta\tensor{}M) \circ (\xi\tensor{}M) \circ \delta \,\,
  \overset{\ref{Equivalent-Def}(ii)}{=}\,\, \cc\tensor{}M
\end{eqnarray*}
The coassociative property is obtained as follows:
\begin{multline*}
(\Delta'\tensor{}\cc\tensor{}M) \circ \Delta'  \,=\,
(\cc\tensor{}\mmm\tensor{}M\tensor{}\cc\tensor{}M) \circ
(\cc\tensor{}\delta\tensor{}\cc\tensor{}M) \circ
(\Delta\tensor{}M\tensor{}\cc\tensor{}M) \circ
(\cc\tensor{}\mmm\tensor{}M) \circ (\cc\tensor{}\delta) \circ
(\Delta\tensor{}M) \\ \,=\,
(\cc\tensor{}\mmm\tensor{}M\tensor{}\cc\tensor{}M) \circ
(\cc\tensor{}\delta\tensor{}\cc\tensor{}M) \circ
(\cc\tensor{}\cc\tensor{}\mmm\tensor{}M) \circ
(\Delta\tensor{}\cc\tensor{}M\tensor{}M) \circ (\cc\tensor{}\delta)
\circ (\Delta\tensor{}M) \\ \,=\,
(\cc\tensor{}\mmm\tensor{}M\tensor{}\cc\tensor{}M) \circ
(\cc\tensor{}\delta\tensor{}\cc\tensor{}M) \circ
(\cc\tensor{}\cc\tensor{}\mmm\tensor{}M) \circ
(\cc\tensor{}\cc\tensor{}\delta) \circ
(\Delta\tensor{}\cc\tensor{}M) \circ (\Delta\tensor{}M) \\ \,=\,
(\cc\tensor{}\mmm\tensor{}M\tensor{}\cc\tensor{}M) \circ
(\cc\tensor{}\delta\tensor{}\cc\tensor{}M) \circ
(\cc\tensor{}\cc\tensor{}\mmm\tensor{}M) \circ
(\cc\tensor{}\cc\tensor{}\delta) \circ
(\cc\tensor{}\Delta\tensor{}M) \circ (\Delta\tensor{}M) \\
\overset{\eqref{2-cell}}{=}
(\cc\tensor{}\mmm\tensor{}M\tensor{}\cc\tensor{}M) \circ
(\cc\tensor{}\delta\tensor{}\cc\tensor{}M) \circ
(\cc\tensor{}\cc\tensor{}\mmm\tensor{}M) \circ
(\cc\tensor{}\Delta\tensor{}M\tensor{}M) \circ
(\cc\tensor{}\delta) \circ (\Delta\tensor{}M) \\ \,=\,
(\cc\tensor{}\mmm\tensor{}M\tensor{}\cc\tensor{}M) \circ
\lr{\cc\tensor{}\lr{ (\delta\tensor{}\cc) \circ (\cc\tensor{}\mmm)
\circ (\Delta\tensor{}M)}\tensor{}M} \circ (\cc\tensor{}\delta)
\circ (\Delta\tensor{}M) \\ \overset{\eqref{2-cell'}}{=}
(\cc\tensor{}\mmm\tensor{}M\tensor{}\cc\tensor{}M) \circ
(\cc\tensor{}\cc\tensor{}M\tensor{}\mmm\tensor{}M) \circ
(\cc\tensor{}\cc\tensor{}\mmm\tensor{}M\tensor{}M) \circ
(\cc\tensor{}\Delta\tensor{}M\tensor{}M\tensor{}M) \circ
(\cc\tensor{}\delta\tensor{}M)\\ \,\, \circ  (\cc\tensor{}\delta)
\circ (\Delta\tensor{}M) \\ \,=\,
(\cc\tensor{}M\tensor{}\cc\tensor{}\mmm\tensor{}M) \circ
(\cc\tensor{}\mmm\tensor{}\cc\tensor{}M\tensor{}M) \circ
(\cc\tensor{}\cc\tensor{}\mmm\tensor{}M\tensor{}M) \circ
(\cc\tensor{}\Delta\tensor{}M\tensor{}M\tensor{}M) \circ
(\cc\tensor{}\delta\tensor{}M)\\ \,\, \circ  (\cc\tensor{}\delta)
\circ (\Delta\tensor{}M) \\ \,=\,
(\cc\tensor{}M\tensor{}\cc\tensor{}\mmm\tensor{}M) \circ
\lr{\cc\tensor{}\lr{(\mmm\tensor{}\cc) \circ (\cc\tensor{}\mmm)
\circ (\Delta\tensor{}M)}\tensor{}M\tensor{}M} \circ
(\cc\tensor{}\delta\tensor{}M) \circ  (\cc\tensor{}\delta) \circ
(\Delta\tensor{}M) \\ \overset{\eqref{1-cell}}{=}
(\cc\tensor{}M\tensor{}\cc\tensor{}\mmm\tensor{}M) \circ
(\cc\tensor{}M\tensor{}\Delta\tensor{}M\tensor{}M) \circ
(\cc\tensor{}\mmm\tensor{}M\tensor{}M)  \circ
(\cc\tensor{}\delta\tensor{}M) \circ (\cc\tensor{}\delta) \circ
(\Delta\tensor{}M) \\ \,=\,
(\cc\tensor{}M\tensor{}\cc\tensor{}\mmm\tensor{}M) \circ
(\cc\tensor{}M\tensor{}\Delta\tensor{}M\tensor{}M) \circ
\lr{\cc\tensor{}\lr{ (\mmm\tensor{}M\tensor{}M)  \circ
(\delta\tensor{}M) \circ \delta}} \circ (\Delta\tensor{}M),
\end{multline*}
applying consecutively Proposition \ref{Equivalent-Def}(ii) and
equation \eqref{2-cell} to the map $\delta$, we get $
(\Delta'\tensor{}\cc\tensor{}M) \circ \Delta' \,=\,
(\cc\tensor{}M\tensor{}\Delta') \circ \Delta'$.

By definition $\xi$ is an $A$--bilinear map compatible with both
counits $\varepsilon$ and $\varepsilon'$. Let us show that $\xi$
is compatible with comultiplications. We have
\begin{eqnarray*}
  (\xi\tensor{}\xi) \circ \Delta' &=& (\xi\tensor{}\xi)
  \circ (\cc\tensor{}\mmm\tensor{}M) \circ (\cc\tensor{}\delta) \circ
(\Delta\tensor{}M)  \\
   &=& (\cc\tensor{}\xi) \circ (\xi\tensor{}\cc\tensor{}M)
   \circ (\cc\tensor{}\mmm\tensor{}M) \circ (\cc\tensor{}\delta) \circ
(\Delta\tensor{}M) \\
   &\overset{\eqref{2-cell}}{=}& (\cc\tensor{}\xi) \circ (\xi\tensor{}\cc\tensor{}M)
   \circ (\cc\tensor{}\mmm\tensor{}M) \circ (\Delta\tensor{}M\tensor{}M) \circ \delta  \\
   &\overset{\eqref{2-cell'}}{=}&  (\cc\tensor{}\xi) \circ
   (\Delta\tensor{}M) \circ (\xi\tensor{}M) \circ \delta
   \,\overset{\ref{Equivalent-Def}(ii)}{=}\, (\cc\tensor{}\xi) \circ
   (\Delta\tensor{}M) \,\overset{\eqref{2-cell}}{=}\, \Delta \circ
   \xi.
\end{eqnarray*}
This proves that $\xi$ is a corings morphism, and finishes the
proof.
\end{proof}

\begin{lemma}\label{tensorM}
Let $(\cc:A)$ be any coring and $(M,\mmm)$ a $\cc$--cowreath with
structure maps $\xi:\cc\tensor{A}M \to \cc$ and
$\delta:\cc\tensor{A}M \to \cc\tensor{A}M\tensor{A}M$. Consider
$\cc\tensor{A}M$ as an $A$--coring with structure given by
proposition \ref{cow-prod}. There is a functor $-\tensor{A}M:
{}_A\mM^{\cc} \to {}_A\mM^{\cc\tensor{A}M}$ sending $(X,\rho^X) \to
\lr{X\tensor{A}M,(X\tensor{A}\mmm\tensor{A}M) \circ
(X\tensor{A}\delta) \circ (\rho^X\tensor{A}M)}$ and $f \to
f\tensor{A}M$.
\end{lemma}
\begin{proof}
We only prove that $-\tensor{}M$ is a well defined functor. Let
$(X,\rho^X)$ be an arbitrary object of the category
${}_A\mM^{\cc}$, i.e. $\rho^X$ is left $A$--linear and right
$\cc$--coaction. Put, $\varrho^{X\tensor{}M}:=
(X\tensor{}\mmm\tensor{}M) \circ (X\tensor{}\delta) \circ
(\rho^X\tensor{}M)$, so we have
\begin{eqnarray*}
  (X\tensor{}M\tensor{}\varepsilon') \circ \varrho^{X\tensor{}M}
  &=& (X\tensor{}M\tensor{}\varepsilon) \circ (X\tensor{}M\tensor{}\xi) \circ
  (X\tensor{}\mmm\tensor{}M) \circ
(X\tensor{}\delta) \circ (\rho^X\tensor{}M)  \\
   &=& (X\tensor{}M\tensor{}\varepsilon) \circ \lr{X\tensor{}\lr{(M\tensor{}\xi) \circ
  (\mmm\tensor{}M) \circ \delta}} \circ (\rho^X\tensor{}M) \\
   &\overset{\ref{Equivalent-Def}(ii)}{=}& (X\tensor{}M\tensor{}\varepsilon)
   \circ (X\tensor{}\mmm)  \circ (\rho^X\tensor{}M) \\
   &\overset{\eqref{1-cell'}}{=}& (X\tensor{}\varepsilon\tensor{}M) \circ (\rho^X\tensor{}M)
   \,\, =\,\, \cc\tensor{}M,
\end{eqnarray*}
which gives the counit property. Now, we have from one hand
\begin{multline*}
(X\tensor{}M\tensor{}\Delta') \circ \varrho^{X\tensor{}M}  \,=\,
(X\tensor{}M\tensor{}\cc\tensor{}\mmm\tensor{}M) \circ
(X\tensor{}M\tensor{}\cc\tensor{}\delta) \circ
(X\tensor{}M\tensor{}\Delta\tensor{}M) \circ
(X\tensor{}\mmm\tensor{}M) \circ (X\tensor{}\delta) \circ (\rho^X\tensor{}M)  \\
\overset{\eqref{2-cell}}{=}
(X\tensor{}M\tensor{}\cc\tensor{}\mmm\tensor{}M) \circ
(X\tensor{}M\tensor{}\Delta\tensor{}M\tensor{}M) \circ
(X\tensor{}M\tensor{}\delta) \circ (X\tensor{}\mmm\tensor{}M)
\circ (X\tensor{}\delta) \circ (\rho^X\tensor{}M) \\
\,=\, (X\tensor{}M\tensor{}\cc\tensor{}\mmm\tensor{}M) \circ
(X\tensor{}M\tensor{}\Delta\tensor{}M\tensor{}M) \circ
\lr{X\tensor{}\lr{(M\tensor{}\delta) \circ (\mmm\tensor{}M) \circ
\delta}} \circ  (\rho^X\tensor{}M)  \\
\overset{\ref{Equivalent-Def}(ii)}{=}
(X\tensor{}M\tensor{}\cc\tensor{}\mmm\tensor{}M) \circ
(X\tensor{}M\tensor{}\Delta\tensor{}M\tensor{}M) \circ
(X\tensor{}\mmm\tensor{}M\tensor{}M) \circ
(X\tensor{}\delta\tensor{}M) \circ (X\tensor{}\delta) \circ
(\rho^X\tensor{}M),
\end{multline*}
and from another
\begin{multline*}
(\varrho^{X\tensor{}M}\tensor{}\cc\tensor{}M) \circ
\varrho^{X\tensor{}M}  \,=\,
(X\tensor{}\mmm\tensor{}M\tensor{}\cc\tensor{}M) \circ
(X\tensor{}\delta\tensor{}\cc\tensor{}M) \circ
(\rho^X\tensor{}M\tensor{}\cc\tensor{}M) \circ
(X\tensor{}\mmm\tensor{}M) \circ (X\tensor{}\delta) \circ
(\rho^X\tensor{}M) \\ \,=\,
(X\tensor{}\mmm\tensor{}M\tensor{}\cc\tensor{}M) \circ
(X\tensor{}\delta\tensor{}\cc\tensor{}M) \circ
(X\tensor{}\cc\tensor{}\mmm\tensor{}M) \circ
(\rho^X\tensor{}\cc\tensor{}M\tensor{}M) \circ (X\tensor{}\delta)
\circ (\rho^X\tensor{}M) \\ \,=\,
(X\tensor{}\mmm\tensor{}M\tensor{}\cc\tensor{}M) \circ
(X\tensor{}\delta\tensor{}\cc\tensor{}M) \circ
(X\tensor{}\cc\tensor{}\mmm\tensor{}M) \circ
(X\tensor{}\cc\tensor{}\delta) \circ (\rho^X\tensor{}\cc\tensor{}M)
\circ (\rho^X\tensor{}M) \\ \,=\,
(X\tensor{}\mmm\tensor{}M\tensor{}\cc\tensor{}M) \circ
(X\tensor{}\delta\tensor{}\cc\tensor{}M) \circ
(X\tensor{}\cc\tensor{}\mmm\tensor{}M) \circ
(X\tensor{}\cc\tensor{}\delta) \circ (X\tensor{}\Delta\tensor{}M)
\circ (\rho^X\tensor{}M) \\ \overset{\eqref{2-cell}}{=}
(X\tensor{}\mmm\tensor{}M\tensor{}\cc\tensor{}M) \circ
(X\tensor{}\delta\tensor{}\cc\tensor{}M) \circ
(X\tensor{}\cc\tensor{}\mmm\tensor{}M) \circ
(X\tensor{}\Delta\tensor{}M\tensor{}M) \circ (X\tensor{}\delta)
\circ (\rho^X\tensor{}M) \\ \,=\,
(X\tensor{}\mmm\tensor{}M\tensor{}\cc\tensor{}M) \circ
\lr{X\tensor{}\lr{(\delta\tensor{}\cc) \circ (\cc\tensor{}\mmm)
\circ (\Delta\tensor{}M)}\tensor{}M} \circ
(X\tensor{}\delta) \circ (\rho^X\tensor{}M)  \\
\overset{\eqref{2-cell'}}{=}
(X\tensor{}\mmm\tensor{}M\tensor{}\cc\tensor{}M) \circ
(X\tensor{}\cc\tensor{}M\tensor{}\mmm\tensor{}M) \circ
(X\tensor{}\cc\tensor{}\mmm\tensor{}M\tensor{}M) \circ
(X\tensor{}\Delta\tensor{}M\tensor{}M\tensor{}M) \circ
(X\tensor{}\delta\tensor{}M) \circ (X\tensor{}\delta) \circ
(\rho^X\tensor{}M) \\ \,=\,
(X\tensor{}M\tensor{}\cc\tensor{}\mmm\tensor{}M) \circ
\lr{X\tensor{}\lr{ (\mmm\tensor{}\cc) \circ (\cc\tensor{}\mmm) \circ
(\Delta\tensor{}M)}\tensor{}M\tensor{}M} \circ
(X\tensor{}\delta\tensor{}M) \circ (X\tensor{}\delta) \circ
(\rho^X\tensor{}M) \\ \overset{\eqref{1-cell}}{=}
(X\tensor{}M\tensor{}\cc\tensor{}\mmm\tensor{}M) \circ
(X\tensor{}M\tensor{}\Delta\tensor{}M\tensor{}M) \circ
(X\tensor{}\mmm\tensor{}M\tensor{}M) \circ
(X\tensor{}\delta\tensor{}M) \circ (X\tensor{}\delta) \circ
(\rho^X\tensor{}M)
\end{multline*}
Therefore, $(X\tensor{}M\tensor{}\Delta') \circ
\varrho^{X\tensor{}M} \,=\,
(\varrho^{X\tensor{}M}\tensor{}\cc\tensor{}M) \circ
\varrho^{X\tensor{}M}$, thus $\varrho^{X\tensor{}M}$ is a right
$\cc\tensor{}M$--coaction. Let $f:(X,\rho^X) \to (X',\rho^{X'})$ be
a right $\cc$--colinear map. It is easily check that
$\varrho^{X'\tensor{}M} \circ (f\tensor{}M) \,=\,
(f\tensor{}\cc\tensor{}M) \circ \varrho^{X\tensor{}M}$, that is
$f\tensor{}M$ is a right $\cc\tensor{}M$--colinear maps.
\end{proof}

By Proposition \ref{cow-prod}, $\xi: \cc\tensor{A}M \to \cc$ is a
morphism of $A$-corings. So one can associated to $\xi$ it induction
functor $(-)_{\xi}: {}_A\mM^{\cc\tensor{A}M} \to {}_A\mM^{\cc}$
which sends $(Y,\rho^Y) \to (Y,\rho^Y)_{\xi}\,=\, \lr{Y,
(Y\tensor{A}\xi) \circ \rho^Y}$ and acts by identity on morphisms,
see \cite{Gomez:2002} or \cite{Brzezinski/Wisbauer:2003}. We prove
by the following, that the induction functor $(-)_{\xi}$ admits
$-\tensor{A}M$ as right adjoint functor.

\begin{proposition}\label{cow-adjoint}
Let $(\cc:A)$ be any coring and $(M,\mmm)$ a $\cc$--cowreath with
structure maps $\xi:\cc\tensor{A}M \to \cc$ and
$\delta:\cc\tensor{A}M \to \cc\tensor{A}M\tensor{A}M$. Consider
$\cc\tensor{A}M$ as an $A$--coring with structure given by
proposition \ref{cow-prod}. There is a natural isomorphism
$$\xymatrix@R=0pt{ {\rm Hom}_{A-\cc}\lr{(Y, \rho^Y)_{\xi},\, (X,\rho^X)}
\ar@{->}[rr] & & {\rm Hom}_{A-(\cc\tensor{A}M)}\lr{(Y,\rho^Y),\,
(X,\rho^X)\tensor{A}M} \\ f \ar@{|->}[rr] &&
\left[(X\tensor{A}M\tensor{A}\varepsilon) \circ
(f\tensor{A}M\tensor{A}\cc) \circ (Y\tensor{A}\mmm) \circ \rho^Y\right] \\
\left[(X\tensor{A}\varepsilon) \circ (X\tensor{A}\xi) \circ
(\rho^X\tensor{A}M) \circ g\right] & & \ar@{|->}[ll] g, }$$ for
every pair of bicomodules $(X,\rho^X) \in {}_A\mM^{\cc}$ and
$(Y,\rho^Y) \in {}_A\mM^{\cc\tensor{A}M}$. That is the induction
functor $(-)_{\xi}$ is left adjoint to the functor $-\tensor{A}M$
defined in lemma \ref{tensorM}.
\end{proposition}
\begin{proof}
We first show that the stated maps are well defined. Let $f:
(Y,\rho^Y)_{\xi} \to (X,\rho^X)$ be any morphism in the category
${}_A\mM^{\cc}$, and denote by $\tha{f}=
(X\tensor{}M\tensor{}\varepsilon) \circ (f\tensor{}M\tensor{}\cc)
\circ (Y\tensor{}\mmm) \circ \rho^Y$ its image. We know that $\rho^X
\circ f = (f\tensor{}\cc) \circ (Y\tensor{}\xi) \circ \rho^Y$, so
compute $\varrho^{X\tensor{}M} \circ \tha{f}$
\begin{multline*}
\varrho^{X\tensor{}M} \circ \tha{f} \,=\, (X\tensor{}\mmm\tensor{}M)
\circ (X\tensor{}\delta) \circ (\rho^X\tensor{}M) \circ
(X\tensor{}M\tensor{}\varepsilon) \circ
(f\tensor{}M\tensor{}\cc) \circ (Y\tensor{}\mmm) \circ \rho^Y \\
\,=\, (X\tensor{}\mmm\tensor{}M) \circ (X\tensor{}\delta) \circ
(X\tensor{}\cc\tensor{}M\tensor{}\varepsilon) \circ
(\rho^X\tensor{}M\tensor{}\cc) \circ (f\tensor{}M\tensor{}\cc) \circ (Y\tensor{}\mmm) \circ \rho^Y \\
\,=\, (X\tensor{}\mmm\tensor{}M) \circ (X\tensor{}\delta) \circ
(X\tensor{}\cc\tensor{}M\tensor{}\varepsilon) \circ
(f\tensor{}\cc\tensor{}M\tensor{}\cc) \circ
(Y\tensor{}\xi\tensor{}M\tensor{}\cc) \circ
(\rho^Y\tensor{}M\tensor{}\cc) \circ (Y\tensor{}\mmm) \circ \rho^Y
\\ \,=\, (X\tensor{}\mmm\tensor{}M) \circ (X\tensor{}\delta) \circ
(X\tensor{}\cc\tensor{}M\tensor{}\varepsilon) \circ
(f\tensor{}\cc\tensor{}M\tensor{}\cc) \circ
(Y\tensor{}\xi\tensor{}M\tensor{}\cc) \circ
(Y\tensor{}\cc\tensor{}M\tensor{}\mmm) \circ
(\rho^Y\tensor{}\cc\tensor{}M) \circ \rho^Y \\ \,=\,
(X\tensor{}\mmm\tensor{}M) \circ (X\tensor{}\delta) \circ
(X\tensor{}\cc\tensor{}M\tensor{}\varepsilon) \circ
(f\tensor{}\cc\tensor{}M\tensor{}\cc) \circ
(Y\tensor{}\xi\tensor{}M\tensor{}\cc) \circ
(Y\tensor{}\cc\tensor{}M\tensor{}\mmm) \circ (Y\tensor{}\Delta')
\circ \rho^Y \\ \,=\, (X\tensor{}\mmm\tensor{}M) \circ
(X\tensor{}\delta) \circ
(X\tensor{}\cc\tensor{}M\tensor{}\varepsilon) \circ
(f\tensor{}\cc\tensor{}M\tensor{}\cc) \circ
(Y\tensor{}\xi\tensor{}M\tensor{}\cc) \circ
(Y\tensor{}\cc\tensor{}M\tensor{}\mmm) \circ
(Y\tensor{}\cc\tensor{}\mmm\tensor{}M) \\ \,\, \circ
(Y\tensor{}\cc\tensor{}\delta) \circ (Y\tensor{}\Delta\tensor{}M)
\circ \rho^Y \\ \overset{\eqref{2-cell}}{=}
(X\tensor{}\mmm\tensor{}M) \circ (X\tensor{}\delta) \circ
(X\tensor{}\cc\tensor{}M\tensor{}\varepsilon) \circ
(f\tensor{}\cc\tensor{}M\tensor{}\cc) \circ
(Y\tensor{}\cc\tensor{}\mmm) \circ
(Y\tensor{}\xi\tensor{}\cc\tensor{}M) \circ
(Y\tensor{}\cc\tensor{}\mmm\tensor{}M) \\ \,\, \circ
(Y\tensor{}\Delta\tensor{}M\tensor{}M) \circ (Y\tensor{}\delta)
\circ \rho^Y \\ \overset{\eqref{2-cell'}}{=}
(X\tensor{}\mmm\tensor{}M) \circ (X\tensor{}\delta) \circ
(X\tensor{}\cc\tensor{}M\tensor{}\varepsilon) \circ
(f\tensor{}\cc\tensor{}M\tensor{}\cc) \circ
(Y\tensor{}\cc\tensor{}\mmm) \circ (Y\tensor{}\Delta\tensor{}M)
\circ (Y\tensor{}\xi\tensor{}M) \circ (Y\tensor{}\delta) \circ
\rho^Y \\ \overset{\ref{Equivalent-Def}(ii)}{=}
(X\tensor{}\mmm\tensor{}M) \circ (X\tensor{}\cc\tensor{}M\tensor{}M
\tensor{}\varepsilon) \circ (X\tensor{}\delta\tensor{}\cc) \circ
(f\tensor{}\cc\tensor{}M\tensor{}\cc) \circ
(Y\tensor{}\cc\tensor{}\mmm) \circ (Y\tensor{}\Delta\tensor{}M)
\circ \rho^Y \\ \,=\, (X\tensor{}\mmm\tensor{}M) \circ
(X\tensor{}\cc\tensor{}M\tensor{}M \tensor{}\varepsilon) \circ
(f\tensor{}\cc\tensor{}M\tensor{}M\tensor{}\cc) \circ
(Y\tensor{}\delta\tensor{}\cc) \circ (Y\tensor{}\cc\tensor{}\mmm)
\circ (Y\tensor{}\Delta\tensor{}M) \circ \rho^Y \\
\overset{\eqref{2-cell'}}{=} (X\tensor{}\mmm\tensor{}M) \circ
(X\tensor{}\cc\tensor{}M\tensor{}M \tensor{}\varepsilon) \circ
(f\tensor{}\cc\tensor{}M\tensor{}M\tensor{}\cc) \circ
(Y\tensor{}\cc\tensor{}M\tensor{}\mmm) \circ
(Y\tensor{}\cc\tensor{}\mmm\tensor{}M) \\ \,\, \circ
(Y\tensor{}\Delta\tensor{}M\tensor{}M) \circ (Y\tensor{}\delta)
\circ \rho^Y \\ \,=\,
(X\tensor{}M\tensor{}\cc\tensor{}M\tensor{}\varepsilon) \circ
(X\tensor{}\mmm\tensor{}M \tensor{}\cc) \circ
(f\tensor{}\cc\tensor{}M\tensor{}M\tensor{}\cc) \circ
(Y\tensor{}\cc\tensor{}M\tensor{}\mmm) \circ
(Y\tensor{}\cc\tensor{}\mmm\tensor{}M) \\ \,\, \circ
(Y\tensor{}\Delta\tensor{}M\tensor{}M) \circ (Y\tensor{}\delta)
\circ \rho^Y \\ \,=\,
(X\tensor{}M\tensor{}\cc\tensor{}M\tensor{}\varepsilon) \circ
(f\tensor{}M\tensor{}\cc\tensor{}M \tensor{}\cc) \circ
(Y\tensor{}\mmm\tensor{}M\tensor{}\cc) \circ
(Y\tensor{}\cc\tensor{}M\tensor{}\mmm) \circ
(Y\tensor{}\cc\tensor{}\mmm\tensor{}M) \\ \,\, \circ
(Y\tensor{}\Delta\tensor{}M\tensor{}M) \circ (Y\tensor{}\delta)
\circ \rho^Y \\ \,=\,
(X\tensor{}M\tensor{}\cc\tensor{}M\tensor{}\varepsilon) \circ
(f\tensor{}M\tensor{}\cc\tensor{}M \tensor{}\cc) \circ
(Y\tensor{}M\tensor{}\cc\tensor{}\mmm) \circ
(Y\tensor{}\mmm\tensor{}\cc\tensor{}M) \circ
(Y\tensor{}\cc\tensor{}\mmm\tensor{}M) \\ \,\, \circ
(Y\tensor{}\Delta\tensor{}M\tensor{}M) \circ (Y\tensor{}\delta)
\circ \rho^Y,
\end{multline*}
applying consecutively equations \eqref{1-cell} and
\eqref{1-cell'} to the twist map $\mmm$, we get
$$\varrho^{X\tensor{}M} \circ \tha{f} \,\,=\,\,
(f\tensor{}M\tensor{}\cc\tensor{}M)  \circ
(Y\tensor{}\mmm\tensor{}M) \circ (Y\tensor{}\delta) \circ \rho^Y
$$ Now, we compute $(\tha{f}\tensor{}\cc\tensor{}M)
\circ \rho^Y$, we find
\begin{eqnarray*}
  (\tha{f}\tensor{}\cc\tensor{}M) \circ \rho^Y
  &=& (X\tensor{}M\tensor{}\varepsilon\tensor{}\cc\tensor{}M) \circ
(f\tensor{}M\tensor{}\cc\tensor{}\cc\tensor{}M) \circ
(Y\tensor{}\mmm\tensor{}\cc\tensor{}M) \circ
(\rho^Y\tensor{}\cc\tensor{}M) \circ \rho^Y   \\
   &=& (X\tensor{}M\tensor{}\varepsilon\tensor{}\cc\tensor{}M) \circ
(f\tensor{}M\tensor{}\cc\tensor{}\cc\tensor{}M) \circ
(Y\tensor{}\mmm\tensor{}\cc\tensor{}M) \circ (Y\tensor{}\Delta')
\circ \rho^Y \\
   &=& (X\tensor{}M\tensor{}\varepsilon\tensor{}\cc\tensor{}M) \circ
(f\tensor{}M\tensor{}\cc\tensor{}\cc\tensor{}M) \circ
(Y\tensor{}\mmm\tensor{}\cc\tensor{}M) \circ
(Y\tensor{}\cc\tensor{}\mmm\tensor{}M) \\ &\,\,& \,\circ\,\,
(Y\tensor{}\cc\tensor{}\delta) \circ (Y\tensor{}\Delta\tensor{}M)
\circ \rho^Y  \\
   &\overset{\eqref{2-cell}}{=}&
(X\tensor{}M\tensor{}\varepsilon\tensor{}\cc\tensor{}M) \circ
(f\tensor{}M\tensor{}\cc\tensor{}\cc\tensor{}M) \circ
(Y\tensor{}\mmm\tensor{}\cc\tensor{}M) \circ
(Y\tensor{}\cc\tensor{}\mmm\tensor{}M) \\ &\,\,& \,\circ\,\,
(Y\tensor{}\Delta\tensor{}M\tensor{}M) \circ
(Y\tensor{}\delta) \circ \rho^Y  \\
& \overset{\eqref{1-cell}}{=} &
 (f\tensor{}M\tensor{}\cc\tensor{}M) \circ
(Y\tensor{}M\tensor{}\varepsilon\tensor{}\cc\tensor{}M) \circ
(Y\tensor{}M\tensor{}\Delta\tensor{}M) \circ
(Y\tensor{}\mmm\tensor{}M) \circ (Y\tensor{}\delta) \circ \rho^Y \\
   &=& (f\tensor{}M\tensor{}\cc\tensor{}M) \circ
(Y\tensor{}\mmm\tensor{}M) \circ (Y\tensor{}\delta) \circ \rho^Y
\end{eqnarray*}
Comparing the last two computation, we then get
$\varrho^{X\tensor{}M} \circ \tha{f}\,=\,
(\tha{f}\tensor{}\cc\tensor{}M) \circ \rho^Y$. Thus $\tha{f}$ is
right $\cc\tensor{}M$--colinear.

Let $g:(Y,\rho^Y) \to (X\tensor{}M,\varrho^{X\tensor{}M})$ be now
any left $A$--linear and right $\cc\tensor{}M$--colinear map. The
last condition means that $g$ satisfies
\begin{equation}\label{td-g}
(X\tensor{}\mmm\tensor{}M) \circ (X\tensor{}\delta) \circ
(\rho^X\tensor{}M) \circ g \,=\, (g\tensor{}\cc\tensor{}M) \circ
\rho^Y.
\end{equation}
Denote by $\td{g}= (X\tensor{}\varepsilon) \circ (X\tensor{}\xi)
\circ (\rho^X\tensor{}M) \circ g$ the image of $g$. It is clear
that $\td{g}$ is left $A$--linear. Let us check that it is right
$\cc$--colinearity. So, we have from one hand that
\begin{eqnarray*}
  \rho^X \circ \td{g} &=& \rho^X \circ (X\tensor{}\varepsilon) \circ (X\tensor{}\xi) \circ
(\rho^X\tensor{}M) \circ g \\
   &=& (X\tensor{}\cc\tensor{}\varepsilon) \circ (\rho^X\tensor{}\cc)
    \circ (X\tensor{}\xi)  \circ
(\rho^X\tensor{}M) \circ g \\
   &=& (X\tensor{}\cc\tensor{}\varepsilon) \circ (X\tensor{}\cc\tensor{}\xi)
    \circ (\rho^X\tensor{}\cc\tensor{}M)  \circ
(\rho^X\tensor{}M) \circ g \\
   &=& (X\tensor{}\cc\tensor{}\varepsilon) \circ (X\tensor{}\cc\tensor{}\xi)
    \circ (X\tensor{}\Delta\tensor{}M)  \circ
(\rho^X\tensor{}M) \circ g \\
   &=& (X\tensor{}\cc\tensor{}\varepsilon) \circ (X\tensor{}\Delta)
    \circ (X\tensor{}\xi)  \circ (\rho^X\tensor{}M) \circ g \\
   &=& (X\tensor{}\xi)  \circ (\rho^X\tensor{}M) \circ g,
\end{eqnarray*}
and from another hand, we have
\begin{eqnarray*}
   (\td{g}\tensor{}\cc) \circ (Y\tensor{}\xi) \circ \rho^Y &=& (X\tensor{}\varepsilon\tensor{}\cc) \circ
(X\tensor{}\xi\tensor{}\cc) \circ (\rho^X\tensor{}M\tensor{}\cc)
\circ (g\tensor{}\cc) \circ (Y\tensor{}\xi) \circ \rho^Y  \\
   &=& (X\tensor{}\varepsilon\tensor{}\cc) \circ
(X\tensor{}\xi\tensor{}\cc) \circ (\rho^X\tensor{}M\tensor{}\cc)
\circ (X\tensor{}M\tensor{}\xi) \circ (g\tensor{}\cc\tensor{}M) \circ \rho^Y  \\
  & \overset{\eqref{td-g}}{=} & (X\tensor{}\varepsilon\tensor{}\cc) \circ
(X\tensor{}\xi\tensor{}\cc) \circ (\rho^X\tensor{}M\tensor{}\cc)
\circ (X\tensor{}M\tensor{}\xi) \circ(X\tensor{}\mmm\tensor{}M)
\circ (X\tensor{}\delta) \\ &\,\,& \, \circ \,\, (\rho^X\tensor{}M) \circ g \\
   & \overset{\ref{Equivalent-Def}(ii)}{=} &
(X\tensor{}\varepsilon\tensor{}\cc) \circ
(X\tensor{}\xi\tensor{}\cc) \circ (\rho^X\tensor{}M\tensor{}\cc)
\circ (X\tensor{}\mmm) \circ (\rho^X\tensor{}M) \circ g \\
   &=& (X\tensor{}\varepsilon\tensor{}\cc) \circ
(X\tensor{}\xi\tensor{}\cc) \circ (X\tensor{}\mmm) \circ
(\rho^X\tensor{}\cc\tensor{}M) \circ (\rho^X\tensor{}M) \circ g  \\
   &=& (X\tensor{}\varepsilon\tensor{}\cc) \circ
(X\tensor{}\xi\tensor{}\cc) \circ (X\tensor{}\mmm) \circ
(X\tensor{}\Delta\tensor{}M) \circ (\rho^X\tensor{}M) \circ g  \\
   &\overset{\eqref{2-cell}}{=}& (X\tensor{}\varepsilon\tensor{}\cc)
\circ (X\tensor{}\Delta) \circ (X\tensor{}\xi) \circ
(\rho^X\tensor{}M) \circ g  \,\,=\,\, (X\tensor{}\xi) \circ
(\rho^X\tensor{}M) \circ g.
\end{eqnarray*}
These imply the equality $\rho^X \circ \td{g}\,=\,
(\td{g}\tensor{}\cc) \circ (Y\tensor{}\xi) \circ \rho^Y$, that is
$\td{g}$ is right $\cc$--colinear.

Next, we show that both maps $\td{-}$ and $\tha{-}$ are mutually
inverse. So let $f$ and $g$ two morphisms as above. By definitions,
we have
\begin{multline*}
\tha{\td{g}} \,=\, (X\tensor{}M\tensor{}\varepsilon) \circ
(\td{g}\tensor{}M\tensor{}\cc) \circ (Y\tensor{}\mmm) \circ \rho^Y
\\ \,=\, (X\tensor{}M\tensor{}\varepsilon) \circ
(X\tensor{}\varepsilon\tensor{}M\tensor{}\cc) \circ
(X\tensor{}\xi\tensor{}M\tensor{}\cc) \circ
(\rho^X\tensor{}M\tensor{}M\tensor{}\cc) \circ
(g\tensor{}M\tensor{}\cc) \circ (Y\tensor{}\mmm) \circ \rho^Y \\
\,=\, (X\tensor{}M\tensor{}\varepsilon) \circ
(X\tensor{}\varepsilon\tensor{}M\tensor{}\cc) \circ
(X\tensor{}\xi\tensor{}M\tensor{}\cc) \circ
(\rho^X\tensor{}M\tensor{}M\tensor{}\cc) \circ
(X\tensor{}M\tensor{}\mmm) \circ (g\tensor{}\cc\tensor{}M) \circ
\rho^Y \\ \overset{\eqref{td-g}}{=}
(X\tensor{}M\tensor{}\varepsilon) \circ
(X\tensor{}\varepsilon\tensor{}M\tensor{}\cc) \circ
(X\tensor{}\xi\tensor{}M\tensor{}\cc) \circ
(\rho^X\tensor{}M\tensor{}M\tensor{}\cc) \circ
(X\tensor{}M\tensor{}\mmm) \circ (X\tensor{}\mmm\tensor{}M) \circ
(X\tensor{}\delta) \circ (\rho^X\tensor{}M) \circ g \\ \,=\,
(X\tensor{}M\tensor{}\varepsilon) \circ
(X\tensor{}\varepsilon\tensor{}M\tensor{}\cc) \circ
(X\tensor{}\xi\tensor{}M\tensor{}\cc) \circ
(X\tensor{}\cc\tensor{}M\tensor{}\mmm) \circ
(\rho^X\tensor{}M\tensor{}\cc\tensor{}M) \circ
(X\tensor{}\mmm\tensor{}M) \circ (X\tensor{}\delta) \circ
(\rho^X\tensor{}M) \circ g \\ \,=\,
(X\tensor{}M\tensor{}\varepsilon) \circ
(X\tensor{}\varepsilon\tensor{}M\tensor{}\cc) \circ
(X\tensor{}\xi\tensor{}M\tensor{}\cc) \circ
(X\tensor{}\cc\tensor{}M\tensor{}\mmm) \circ (X
\tensor{}\cc\tensor{}\mmm\tensor{}M) \circ
(\rho^X\tensor{}\cc\tensor{}M\tensor{}M) \circ (X\tensor{}\delta)
\circ (\rho^X\tensor{}M) \circ g \\ \,=\,
(X\tensor{}M\tensor{}\varepsilon) \circ
(X\tensor{}\varepsilon\tensor{}M\tensor{}\cc) \circ
(X\tensor{}\xi\tensor{}M\tensor{}\cc) \circ
(X\tensor{}\cc\tensor{}M\tensor{}\mmm) \circ (X
\tensor{}\cc\tensor{}\mmm\tensor{}M) \circ
(X\tensor{}\cc\tensor{}\delta) \circ
(\rho^X\tensor{}\cc\tensor{}M) \circ (\rho^X\tensor{}M) \circ g \\
\,=\, (X\tensor{}M\tensor{}\varepsilon) \circ
(X\tensor{}\varepsilon\tensor{}M\tensor{}\cc) \circ
(X\tensor{}\xi\tensor{}M\tensor{}\cc) \circ
(X\tensor{}\cc\tensor{}M\tensor{}\mmm) \circ (X
\tensor{}\cc\tensor{}\mmm\tensor{}M) \circ
(X\tensor{}\cc\tensor{}\delta) \circ (X\tensor{}\Delta\tensor{}M)
\circ (\rho^X\tensor{}M) \circ g \\ \overset{\eqref{2-cell}}{=}
(X\tensor{}M\tensor{}\varepsilon) \circ
(X\tensor{}\varepsilon\tensor{}M\tensor{}\cc) \circ
(X\tensor{}\xi\tensor{}M\tensor{}\cc) \circ
(X\tensor{}\cc\tensor{}M\tensor{}\mmm) \circ (X
\tensor{}\cc\tensor{}\mmm\tensor{}M) \circ
(X\tensor{}\Delta\tensor{}M\tensor{}M) \circ (X\tensor{}\delta)
\circ (\rho^X\tensor{}M) \circ g \\ \,=\,
(X\tensor{}M\tensor{}\varepsilon) \circ
(X\tensor{}\varepsilon\tensor{}M\tensor{}\cc) \circ
(X\tensor{}\cc\tensor{}\mmm) \circ
(X\tensor{}\xi\tensor{}\cc\tensor{}M) \circ (X
\tensor{}\cc\tensor{}\mmm\tensor{}M) \circ
(X\tensor{}\Delta\tensor{}M\tensor{}M) \circ (X\tensor{}\delta)
\circ (\rho^X\tensor{}M) \circ g \\ \overset{\eqref{2-cell'}}{=}
(X\tensor{}M\tensor{}\varepsilon) \circ
(X\tensor{}\varepsilon\tensor{}M\tensor{}\cc) \circ
(X\tensor{}\cc\tensor{}\mmm) \circ (X\tensor{}\Delta\tensor{}M)
\circ (X \tensor{}\xi\tensor{}M)  \circ (X\tensor{}\delta) \circ
(\rho^X\tensor{}M) \circ g \\
\overset{\eqref{Equivalent-Def}(ii)}{=}
(X\tensor{}M\tensor{}\varepsilon) \circ (X\tensor{}\mmm) \circ
(X\tensor{}\varepsilon\tensor{}\cc\tensor{}M) \circ
(X\tensor{}\Delta\tensor{}M) \circ  (\rho^X\tensor{}M) \circ g \\
\overset{\eqref{1-cell'}}{=} (X\tensor{}\varepsilon\tensor{}M) \circ
(\rho^X\tensor{}M) \circ g \,\,=\,\, g.
\end{multline*}
and
\begin{eqnarray*}
  \td{\tha{f}} &=& (X\tensor{}\varepsilon) \circ (X\tensor{}\xi) \circ (\rho^X\tensor{}M) \circ \tha{f} \\
   &=& (X\tensor{}\varepsilon) \circ (X\tensor{}\xi) \circ (\rho^X\tensor{}M)
   \circ (X\tensor{}M\tensor{}\varepsilon) \circ (f\tensor{}M\tensor{}\cc)
\circ (Y\tensor{}\mmm) \circ \rho^Y\\
   &=& (X\tensor{}\varepsilon) \circ (X\tensor{}\xi) \circ (\rho^X\tensor{}M)
   \circ (f\tensor{}M) \circ (Y\tensor{}M\tensor{}\varepsilon)
\circ (Y\tensor{}\mmm) \circ \rho^Y \\
   &\overset{\eqref{1-cell'}}{=}&  (X\tensor{}\varepsilon) \circ (X\tensor{}\xi) \circ (f\tensor{}\cc\tensor{}M)
   \circ (Y\tensor{}\xi\tensor{}M) \circ (\rho^Y\tensor{}M)
\circ (Y\tensor{}\varepsilon\tensor{}M) \circ \rho^Y \\
   &=& (X\tensor{}\varepsilon) \circ (f\tensor{}\cc) \circ (Y\tensor{}\xi)
   \circ (Y\tensor{}\xi\tensor{}M) \circ (Y\tensor{}\cc\tensor{}M \tensor{}\varepsilon\tensor{}M)
\circ (\rho^Y\tensor{}\cc\tensor{}M) \circ \rho^Y  \\
   &=& f \circ (Y\tensor{}\varepsilon) \circ  (Y\tensor{}\xi)
   \circ (Y\tensor{}\xi\tensor{}M) \circ (Y\tensor{}\cc\tensor{}M \tensor{}\varepsilon\tensor{}M)
\circ (Y\tensor{}\Delta') \circ \rho^Y \\
   &=& f \circ (Y\tensor{}\varepsilon) \circ  (Y\tensor{}\xi)
   \circ (Y\tensor{}\xi\tensor{}M) \circ (Y\tensor{}\cc\tensor{}M \tensor{}\varepsilon\tensor{}M)
\circ (Y\tensor{}\cc\tensor{}\mmm\tensor{}M) \circ
(Y\tensor{}\cc\tensor{}\delta) \\ &\,\,& \,
\circ\, (Y\tensor{}\Delta\tensor{}M) \circ \rho^Y \\
   &\overset{\eqref{1-cell'}}{=}&
   f \circ (Y\tensor{}\varepsilon) \circ  (Y\tensor{}\xi)
   \circ (Y\tensor{}\xi\tensor{}M) \circ (Y\tensor{}\cc\tensor{}\varepsilon\tensor{}M\tensor{}M) \circ
(Y\tensor{}\cc\tensor{}\delta) \circ (Y\tensor{}\Delta\tensor{}M) \circ \rho^Y \\
&\overset{\eqref{2-cell}}{=}&
   f \circ (Y\tensor{}\varepsilon) \circ  (Y\tensor{}\xi)
   \circ (Y\tensor{}\xi\tensor{}M) \circ (Y\tensor{}\cc\tensor{}\varepsilon\tensor{}M\tensor{}M) \circ
(Y\tensor{}\Delta\tensor{}M\tensor{}M) \circ\, (Y\tensor{}\delta) \circ \rho^Y \\
   &=& f \circ (Y\tensor{}\varepsilon) \circ  (Y\tensor{}\xi)
   \circ (Y\tensor{}\xi\tensor{}M) \circ  (Y\tensor{}\delta)
\circ \rho^Y \\ &\overset{\ref{Equivalent-Def}(ii)}{=}& f \circ
(Y\tensor{}\varepsilon) \circ (Y\tensor{}\xi) \circ
\rho^Y\,\,=\,\, f \circ (Y\tensor{}\varepsilon') \circ
\rho^Y\,\,=\,\, f.
\end{eqnarray*}
\end{proof}

\begin{remark}
In general, for any $A$--corings morphism $\phi: \dd \to \cc$ , one
can associated to it the induction functor $(-)_{\phi}: \mM^{\dd}
\to \mM^{\cc}$. In this way $\dd$ becomes a $(\cc, \dd)$--bicomodule
and $(\dd,\cc)$--bicomodule. If the left module ${}_A\dd$ is assumed
to be flat or preserves certain equalizers. Then one can define,
using this bicomodule, the cotensor functor $-\cotensor{\cc}\dd:
\mM^{\cc} \to \mM^{\dd}$. In that case $(-)_{\phi}$ admits
$-\cotensor{\cc}\dd$ as right adjoint, see \cite{Gomez:2002} or
\cite{Brzezinski/Wisbauer:2003} for more details. Proposition
\ref{cow-adjoint} gives a right adjoint functor to $(-)_{\xi}$
without requiring any assumption on the cowreath product
$\cc\tensor{A}M$. Of course if we asume that $(\cc\tensor{A}M)$ is a
flat left $A$--module, then one can construct the cotensor product
functor $-\cotensor{\cc}(\cc\tensor{A}M): {}_A\mM^{\cc} \to
{}_A\mM^{\cc\tensor{A}M}$ associated to the coring morphisms
$\xi:\cc\tensor{A}M \to \cc$. In such case we have $(-)_{\xi}$ is
left adjoint to $-\cotensor{\cc}(\cc\tensor{A}M)$, and by
Proposition \ref{cow-adjoint}, we obtain a natural isomorphism
$-\tensor{A}M \cong -\cotensor{\cc}(\cc\tensor{A}M)$.
\end{remark}

\subsection{The functor $\oO: \wrcomod{(M,\,\mmm)}{\cc} \to
{}_A\rcomod{\cc\tensor{}M}$.}\label{oO} $ $

For any $\cc$--cowreath $(M,\mmm)$, we denote by
$\wrcomod{(M,\,\mmm)}{\cc}$ its category of all right
$(M,\mmm)$--comodules. Let $(X,\xxx)$ be an object of the category
$\wrcomod{(M,\,\mmm)}{\cc}$ with right coaction $\rho^{(X,\xxx)}:
\cc\tensor{A}X \to \cc\tensor{A}X\tensor{A}M$. Recall that
$\rho^{(X,\xxx)}$ is a morphism of $\cc$--bicomodules satisfying the
conditions of Proposition \ref{comod}(a)(ii). Define
$$\xymatrix@C=50pt{\varrho^{\cc\tensor{}X}: \cc\tensor{}X \ar@{->}^-{\Delta\tensor{}X}[r] &
\cc\tensor{}\cc\tensor{}X \ar@{->}^-{\cc\tensor{}\rho^{(X,\xxx)}}[r]
& \cc\tensor{}\cc\tensor{}X\tensor{}M
\ar@{->}^-{\cc\tensor{}\xxx\tensor{}M}[r] &
\cc\tensor{}X\tensor{}\cc\tensor{}M. }$$ We claim that
$(\cc\tensor{}X,\varrho^{\cc\tensor{}X})$ is an object of the
category ${}_A\mM^{\cc\tensor{}M}$. That is
$\varrho^{\cc\tensor{}X}$ is a left $A$--linear and right
$\cc\tensor{}M$--colinear map. The first property is clear, and the
second one is obtained as follows. We have
\begin{eqnarray*}
  (\cc\tensor{}X\tensor{}\varepsilon') \circ \varrho^{\cc\tensor{}X} &=&
  (\cc\tensor{}X\tensor{}\varepsilon) \circ (\cc\tensor{}X\tensor{}\xi)
  \circ (\cc\tensor{}\xxx\tensor{}M) \circ (\cc\tensor{}\rho^{(X,\xxx)}) \circ (\Delta\tensor{}X)  \\
   &\overset{\ref{comod}(a)(ii)}{=}& (\cc\tensor{}X\tensor{}\varepsilon) \circ
   (\cc\tensor{}\xxx) \circ (\Delta\tensor{}X)   \\
   &\overset{\eqref{1-cell'}}{=}&  (\cc\tensor{}\varepsilon\tensor{}X) \circ (\Delta\tensor{}X)
   \,\,=\,\, \cc\tensor{}X.
\end{eqnarray*}
The coassociativity property comes out as
\begin{multline*}
(\varrho^{\cc\tensor{}M} \tensor{}\cc\tensor{}M) \circ
\varrho^{\cc\tensor{}M} \,=\,
(\cc\tensor{}\xxx\tensor{}M\tensor{}\cc\tensor{}M) \circ
(\cc\tensor{}\rho^{(X,\xxx)}\tensor{}\cc\tensor{}M) \circ
(\Delta\tensor{}X\tensor{}\cc\tensor{}M) \circ
(\cc\tensor{}\xxx\tensor{}M) \circ (\cc\tensor{}\rho^{(X,\xxx)})
\circ (\Delta\tensor{}X) \\ \,=\,
(\cc\tensor{}\xxx\tensor{}M\tensor{}\cc\tensor{}M) \circ
(\cc\tensor{}\rho^{(X,\xxx)}\tensor{}\cc\tensor{}M) \circ
(\cc\tensor{}\cc\tensor{}\xxx\tensor{}M) \circ
(\Delta\tensor{}\cc\tensor{}X\tensor{}M) \circ
(\cc\tensor{}\rho^{(X,\xxx)}) \circ (\Delta\tensor{}X) \\
\overset{\eqref{2-cell}}{=}
(\cc\tensor{}\xxx\tensor{}M\tensor{}\cc\tensor{}M) \circ
(\cc\tensor{}\rho^{(X,\xxx)}\tensor{}\cc\tensor{}M) \circ
(\cc\tensor{}\cc\tensor{}\xxx\tensor{}M) \circ
(\Delta\tensor{}\cc\tensor{}X\tensor{}M) \circ
(\Delta\tensor{}X\tensor{}M) \circ \rho^{(X,\xxx)} \\ \,=\,
(\cc\tensor{}\xxx\tensor{}M\tensor{}\cc\tensor{}M) \circ
(\cc\tensor{}\rho^{(X,\xxx)}\tensor{}\cc\tensor{}M) \circ
(\cc\tensor{}\cc\tensor{}\xxx\tensor{}M) \circ
(\cc\tensor{}\Delta\tensor{}X\tensor{}M) \circ
(\Delta\tensor{}X\tensor{}M) \circ \rho^{(X,\xxx)} \\
\overset{\eqref{2-cell'}}{=}
(\cc\tensor{}\xxx\tensor{}M\tensor{}\cc\tensor{}M) \circ
(\cc\tensor{}\cc\tensor{}X\tensor{}\mmm\tensor{}M) \circ
(\cc\tensor{}\cc\tensor{}\xxx\tensor{}M\tensor{}M) \circ
(\cc\tensor{}\Delta\tensor{}X\tensor{}M\tensor{}M) \circ
(\cc\tensor{}\rho^{(X,\xxx)}\tensor{}M) \circ
(\Delta\tensor{}X\tensor{}M) \circ \rho^{(X,\xxx)} \\ \,=\,
(\cc\tensor{}X\tensor{}\cc\tensor{}\mmm\tensor{}M) \circ
(\cc\tensor{}\xxx\tensor{}\cc\tensor{}M\tensor{}M) \circ
(\cc\tensor{}\cc\tensor{}\xxx\tensor{}M\tensor{}M) \circ
(\cc\tensor{}\Delta\tensor{}X\tensor{}M\tensor{}M) \circ
(\cc\tensor{}\rho^{(X,\xxx)}\tensor{}M) \circ
(\Delta\tensor{}X\tensor{}M) \circ \rho^{(X,\xxx)} \\
\overset{\eqref{1-cell}}{=}
(\cc\tensor{}X\tensor{}\cc\tensor{}\mmm\tensor{}M) \circ
(\cc\tensor{}X\tensor{}\Delta\tensor{}M\tensor{}M) \circ
(\cc\tensor{}\xxx\tensor{}M\tensor{}M)  \circ
(\cc\tensor{}\rho^{(X,\xxx)}\tensor{}M) \circ
(\Delta\tensor{}X\tensor{}M) \circ \rho^{(X,\xxx)} \\ \,=\,
(\cc\tensor{}X\tensor{}\cc\tensor{}\mmm\tensor{}M) \circ
(\cc\tensor{}X\tensor{}\Delta\tensor{}M\tensor{}M) \circ
(\cc\tensor{}\xxx\tensor{}M\tensor{}M)  \circ
(\cc\tensor{}\rho^{(X,\xxx)}\tensor{}M) \circ
(\cc\tensor{}\rho^{(X,\xxx)}) \circ \rho^{(X,\xxx)} \\
\overset{\ref{comod}(a)(ii)}{=}
(\cc\tensor{}X\tensor{}\cc\tensor{}\mmm\tensor{}M) \circ
(\cc\tensor{}X\tensor{}\Delta\tensor{}M\tensor{}M) \circ
(\cc\tensor{}X\tensor{}\delta)  \circ (\cc\tensor{}\xxx\tensor{}M)
\circ (\cc\tensor{}\rho^{(X,\xxx)}) \circ \rho^{(X,\xxx)} \\
\overset{\eqref{2-cell}}{=}
(\cc\tensor{}X\tensor{}\cc\tensor{}\mmm\tensor{}M) \circ
(\cc\tensor{}X\tensor{}\cc\tensor{}\delta) \circ
(\cc\tensor{}X\tensor{}\Delta\tensor{}M)  \circ
\varrho^{\cc\tensor{}X} \,\, =\,\, (\cc\tensor{}X\tensor{}\Delta')
\circ \varrho^{\cc\tensor{}X},
\end{multline*}
which finishes the proof of the claim. Let $f: (X,\xxx) \to
(X',\xxx')$ be any morphism in category of right comodules
$\wrcomod{(M,\,\mmm)}{\cc}$. Then by Proposition \ref{comod}(b), $f$
satisfies $\rho^{(X',\xxx')} \circ f\,=\,(f\tensor{}M) \circ
\rho^{(X,\xxx)}$, using this equality, we get
\begin{eqnarray*}
  (f\tensor{}\cc\tensor{}M) \circ \varrho^{\cc\tensor{}X} &=&
  (f\tensor{}\cc\tensor{}M) \circ (\cc\tensor{}\xxx\tensor{}M) \circ
  (\cc\tensor{}\rho^{(X,\xxx)}) \circ (\Delta\tensor{}M)  \\
   &\overset{\eqref{2-cell}}{=}& (f\tensor{}\cc\tensor{}M) \circ (\cc\tensor{}\xxx\tensor{}M) \circ
  (\Delta\tensor{}X\tensor{}M) \circ \rho^{(X,\xxx)} \\
   &=& \lr{\lr{(f\tensor{}M) \circ (\cc\tensor{}\xxx) \circ (\Delta\tensor{}X)} \tensor{}M} \circ \rho^{(X,\xxx)} \\
   &\overset{\eqref{2-cell'}}{=}& \lr{\lr{(\cc\tensor{}\xxx') \circ (\Delta\tensor{}X')
   \circ f} \tensor{}M} \circ \rho^{(X,\xxx)}  \\
   &=& (\cc\tensor{}\xxx'\tensor{}M) \circ (\Delta\tensor{}X'\tensor{}M) \circ (f\tensor{}M) \circ \rho^{(X,\xxx)} \\
   &=& (\cc\tensor{}\xxx'\tensor{}M) \circ (\Delta\tensor{}X'\tensor{}M) \circ
   \rho^{(X',\xxx')} \circ f \\ &\overset{\eqref{2-cell}}{=}&
   (\cc\tensor{}\xxx'\tensor{}M) \circ
   (\cc\tensor{}\rho^{(X',\xxx')}) \circ (\Delta\tensor{}X') \circ f
   \\ &=& \varrho^{\cc\tensor{}X'} \circ f,
\end{eqnarray*}
which shows that $f$ is a morphism in the category
${}_A\mM^{\cc\tensor{A}M}$.

In conclusion we have establish a faithful functor $\oO:
\wrcomod{(M,\,\mmm)}{\cc} \to {}_A\mM^{\cc\tensor{A}M}$ acting by
$$\Oo\lr{\lr{(X,\xxx), \rho^{(X,\,\xxx)}}}\,=\,\lr{\cc\tensor{A}X,
\varrho^{\cc\tensor{A}X}:= (\cc\tensor{A}\xxx\tensor{A}M) \circ
(\cc\tensor{A}\rho^{(X,\xxx)}) \circ (\Delta\tensor{A}X)}, \qquad
\Oo(f) =f.$$

\subsection{The adjunction }\label{V-W} $\xymatrix{ \vV:
{}_A\mM^{\cc} \ar@<0,5ex>[r] & \ar@<0,5ex>[l] \rR_{(\cc:\, A)}: \wW.
}$

For any object $(X,\xxx)$ in the category $\rR_{(\cc:\,A)}$, we
had seen that $\lr{\cc\tensor{}M, (\cc\tensor{}\xxx) \circ
(\Delta\tensor{}X)}$ is an object in the category of bicomodules
${}_A\mM^{\cc}$. Of course any morphism $\Sf{f}: (X,\xxx) \to
(X',\xxx')$ in $\rR_{(\cc\,:A)}$ is by definition
$\Sf{f}:\cc\tensor{}X \to \cc\tensor{}X'$ an $A-\cc$--bicolinear
map. These in fact establish a functor which we denote by $\wW:
\rR_{(\cc\,:A)} \to {}_A\mM^{\cc}$. In other direction, for every
object $(Z,\rho^Z)$ in the category ${}_A\mM^{\cc}$, we have by
\cite[Lemma 2.1]{Kaoutit:2006c}, that $\lr{Z,\zzz= \rho^Z \circ
(\varepsilon\tensor{}Z)}$ is an object of the category
$\rR_{(\cc:\,A)}$. Given $g:(Z,\rho^Z) \to (Z',\rho^{Z'})$ any
morphisms in ${}_A\mM^{\cc}$, then clearly $\cc\tensor{}g:
\cc\tensor{}Z \to \cc\tensor{}Z'$ is a $\cc$--bicomodule morphism,
since the right $\cc$--coactions on $\cc\tensor{}Z$ and
$\cc\tensor{}Z'$ are given, respectively,  by
$\varrho^{\cc\tensor{}Z}\,=\, (\cc\tensor{}\zzz) \circ
(\Delta\tensor{}Z) \,=\, (\cc\tensor{}\rho^Z) \circ
(\cc\tensor{}\varepsilon\tensor{}Z) \circ (\Delta\tensor{}Z)\,=\,
\cc\tensor{}\rho^Z$ and $\varrho^{\cc\tensor{}Z'}\,=\,
\cc\tensor{}\rho^{Z'}$. That is $\cc\tensor{}g: (Z,\zzz) \to
(Z',\zzz')$ is a morphism in the category $\rR_{(\cc:\,A)}$. These
establish a functor which we denote by $\vV: {}_A\mM^{\cc}
\longrightarrow \rR_{(\cc:\,A)}$.

We claim that $\wW$ is left adjoint functor to $\vV$. Recall that,
for every pair of objects $(Y,\yyy) \in \rR_{(\cc:\,A)}$ and
$(Z,\rho^Z) \in {}_A\mM^{\cc}$, we have
\begin{eqnarray*}
  {\rm Hom}_{\rR_{(\cc:\,A)}}\lr{(Y,\yyy),\, \vV(Z,\rho^Z)} &=&
  {\rm Hom}_{\cc-\cc}\lr{(\cc\tensor{}Y),\,(\cc\tensor{}Z)} \\
  {\rm Hom}_{{}_A\mM^{\cc}}\lr{\wW(Y,\yyy),\, (Z,\rho^Z)} &=& {\rm Hom}_{A-\cc}\lr{(\cc\tensor{}Y),\,
  Z},
\end{eqnarray*}
where $\lambda^{\cc\tensor{}Y}\,=\, \Delta\tensor{}Y$,
$\rho^{\cc\tensor{}Y}\,=\, (\cc\tensor{}\yyy) \circ
(\Delta\tensor{}Y)$, and $\lambda^{\cc\tensor{}Z}\,=\,
\Delta\tensor{}Z$, $\rho^{\cc\tensor{}Z}\,=\, \cc\tensor{}\rho^Z$.
The natural isomorphism which gives the claimed adjunction is given
as follows. To every $g \in {\rm Hom}_{A-\cc}\lr{(\cc\tensor{}Y),\,
Z}$, we set $\tha{g}\,=\,(\cc\tensor{}g) \circ (\Delta\tensor{}Y):
\cc\tensor{}Y \to \cc\tensor{}Z $. While, to every $f \in {\rm
Hom}_{\cc-\cc}\lr{(\cc\tensor{}Y),\,(\cc\tensor{}Z)}$, we set
$\td{f}\,=\, (\varepsilon\tensor{}Z) \circ f: \cc\tensor{}Y \to Z$.
Keeping these notations, we have
\begin{eqnarray*}
  \varrho^{\cc\tensor{}Z} \circ \tha{g} &=& (\cc\tensor{}\rho^Z) \circ (\cc\tensor{}g) \circ (\Delta\tensor{}Y)  \\
   &=& (\cc\tensor{}g\tensor{}\cc) \circ (\cc\tensor{}\varrho^{\cc\tensor{}Y}) \circ (\Delta\tensor{}Y) \\
   &=& (\cc\tensor{}g\tensor{}\cc) \circ (\cc\tensor{}\cc\tensor{}\yyy)
   \circ (\cc\tensor{}\Delta\tensor{}Y) \circ (\Delta\tensor{}Y) \\
   &=& (\cc\tensor{}g\tensor{}\cc) \circ (\cc\tensor{}\cc\tensor{}\yyy)
   \circ (\Delta\tensor{}\cc\tensor{}Y) \circ (\Delta\tensor{}Y) \\
   &=& (\cc\tensor{}g\tensor{}\cc) \circ (\Delta\tensor{}\cc\tensor{}Y)
   \circ (\cc\tensor{}\yyy) \circ (\Delta\tensor{}Y)\\ &=&
   \lr{\lr{(\cc\tensor{}g) \circ (\Delta\tensor{}Y)}\tensor{}\cc}
   \circ \varrho^{\cc\tensor{}Y}\,=\, (\tha{g}\tensor{}\cc) \circ
   \varrho^{\cc\tensor{}Y},
\end{eqnarray*}
which implies that $\tha{g}$ is right $\cc$--colinear map. Hence
$\tha{g} \in {\rm
Hom}_{\cc-\cc}\lr{(\cc\tensor{}Y),\,(\cc\tensor{}Z)}$, since its
already a left $\cc$--colinear map.  This defines a map $\tha{-}:
{\rm Hom}_{{}_A\mM^{\cc}}\lr{\wW(Y,\yyy),\, (Z,\rho^Z)} \to {\rm
Hom}_{\rR_{(\cc:\,A)}}\lr{(Y,\yyy),\, \vV(Z,\rho^Z)}$. Now, we have
\begin{eqnarray*}
  \rho^Z \circ \td{f} &=& \rho^Z \circ (\varepsilon\tensor{}Z) \circ f \\
  &=& (\varepsilon\tensor{}Z\tensor{}\cc) \circ (\cc\tensor{}\rho^Z) \circ f\\
   &=& (\varepsilon\tensor{}Z\tensor{}\cc) \circ (f\tensor{}\cc) \circ (\cc\tensor{}\yyy) \circ (\Delta\tensor{}Y) \\
   &=& (\td{f}\tensor{}\cc) \circ \varrho^{\cc\tensor{}Y}.
\end{eqnarray*}
Thus $\td{f}$ is right $\cc$--colinear map. Whence $\td{f} \in
{\rm Hom}_{A-\cc}\lr{(\cc\tensor{}Y),\,Z}$, as $\td{f}$ is by
definition left $A$--linear. This establishes a map $\td{-}:{\rm
Hom}_{\rR_{(\cc:\,A)}}\lr{(Y,\yyy),\, \vV(Z,\rho^Z)} \to {\rm
Hom}_{{}_A\mM^{\cc}}\lr{\wW(Y,\yyy),\, (Z,\rho^Z)}$. One can
easily check that $\td{\tha{g}}\,=\, g$ and $\tha{\td{f}}\,=\, f$,
for every $f$ and $g$ as above. Therefore, there is a natural
isomorphism
$$ {\rm Hom}_{\rR_{(\cc:\,A)}}\lr{(Y,\yyy),\, \vV(Z,\rho^Z)} \,\,\cong \,\, {\rm
Hom}_{{}_A\mM^{\cc}}\lr{\wW(Y,\yyy),\, (Z,\rho^Z)},$$ for every pair
of objects $(Y,\yyy) \in \rR_{(\cc:\,A)}$ and $(Z,\rho^Z) \in
{}_A\mM^{\cc}$, and the desired claim is proved.

We now summarize the situation by giving the following non
commutative diagram:
\begin{equation}\label{digrama}
\xymatrix@C=40pt@R=40pt{\mM^{(M,\mmm)}_{\cc} \ar@{->}^-{\oO}[rr]
\ar@<0,5ex>^-{\uU}[dd] & &
{}_A\mM^{\cc\tensor{}M} \ar@<0,5ex>^-{(-)_{\xi}}[dd] \\ & & \\
\rR_{(\cc:\,A)} \ar@<0,5ex>^-{\wW}[rr]
\ar@<0,5ex>^-{-\Rtensor{C}(M,\mmm)}[uu] & &
\ar@<0,5ex>^-{-\tensor{}M}[uu] \ar@<0,5ex>^-{\vV}[ll]
{}_A\mM^{\cc},}
\end{equation}
where the adjunction $-\Rtensor{C}(M,\mmm) \dashv \uU$ is the usual
one associated to the category of comodules over a coalgebra defined
in monoidal category, and where the rest of pairwise arrows indicate
the adjunction $\wW \dashv \vV$ given in \ref{V-W} and $(-)_{\xi}
\dashv -\tensor{A}M$ stated in Proposition \ref{cow-adjoint}.

\section{The dual notions: Wreath over rings
extension.}\label{Sect-Dual}

In this section we give without proofs the "dual" of the majority of
results stated in the previous sections. Notice that the notion
"dual" is not at all perfect since there are several duales in the
present context. This is due probably to the fact that any
bicategory admits three kind of dualisation: by reversing $1$-cells,
by reversing $2$-cells, or by reversing both them.

The notion of coring is dual to that of ring. That is given any
ground base ring $A$, and consider its category of $A$-bimodules
$\Bimod{A}{A}$ as monoidal category with multiplication the tensor
product over $A$; an $A$--coring is then a coalgebra in
$\Bimod{A}{A}$, while an $A$--ring is an algebra in $\Bimod{A}{A}$.
In this way an $A$--ring is just a unital ring extension $\iota: A
\to T$ (i.e. a unital morphism of rings). The tensor product
$-\tensor{A}-$ will be denoted in some occasions by $-\tensor{}-$.

From now on, let us fix a rings extension $\iota: A \to T$, which we
expresse by $(A:T)$. The multiplication of $T$ will be denoted by
$\mu$ (or $\mu_T$) and it unit by $1$ (or $1_T$). Associated to
$(A:T)$, as in Section \ref{Sec-1} there is an additive monoidal
category $\rR_{(A:T)}$, defined by  the following data
\begin{enumerate}[$\bullet$]
\item \emph{Objects}. Are pairs $(P,\ppp)$ consisting of an
$A$--bimodule $P$ and an $A$--bilinear map $\ppp: T\tensor{A}P \to
P\tensor{A} T$ satisfying
\begin{eqnarray}
  (P\tensor{}\mu) \circ (\ppp\tensor{}T) \circ (T\tensor{}\ppp) &=& \ppp \circ (\mu\tensor{}P) \label{01}\\
  \ppp \circ (1\tensor{}P) &=& P\tensor{}1 \label{02}
\end{eqnarray}
Given any object $(P,\ppp)$ and any left $T$--module $X$ with action
$\Sf{l}_X: T\tensor{A}X \to X$. Then one can easily check that
$P\tensor{A}X$ inherits an structure of left $T$--module given by
the action $\Sf{l}_{P\tensor{A}X}\,=\, (P\tensor{A}\Sf{l}_X) \circ
(\ppp\tensor{A}T)$. Of course if $X$ is assumed to be a
$T$--bimodule, then $P\tensor{A}X$ becomes also a $T$--bimodule. In
this way for each object $(P,\ppp)$, the $A$--bimodule
$P\tensor{A}T$ will be always considered as a $T$--bimodule.

\item \emph{Morphisms}. For any two objects $(P,\ppp)$ and $(P',\ppp')$, the
$\KK$--module of morphisms is defined by
$$ {\rm Hom}_{\rR_{(A:\,T)}}\lr{ (P,\ppp),\, (P',\ppp')} : =\, {\rm
Hom}_{T-T}\lr{P\tensor{A}T,\, P'\tensor{A}T}.$$
\end{enumerate}

The category $\rR_{(A:\,T)}$ is monoidal with horizontal
multiplication given by
$$ (P,\ppp) \tensorbajo{(A:\,T)} (P',\ppp')\,=\, \lr{ P\tensor{A}P',\,
(P\tensor{A}\ppp') \circ (\ppp\tensor{A}P)} $$ Now, for any pair of
morphisms $\alpha: (P,\ppp) \to (Q,\qqq)$ and $\beta: (P',\ppp') \to
(Q',\qqq')$, the vertical multiplication
$\alpha\tensorbajo{(A:\,T)}\beta$ is defined as the composition map
\[
\xy *+{P\tensor{}P'\tensor{}T}="p",
p+<12cm,0pt>*+{Q\tensor{}Q'\tensor{}T}="1",
p+<0pt,-2cm>*+{P\tensor{}T\tensor{}P'\tensor{}T}="2",
p+<12cm,-2cm>*+{Q\tensor{}Q'\tensor{}T\tensor{}T}="3",
p+<3cm,-4cm>*+{Q\tensor{}T\tensor{}P'\tensor{}T}="4",
p+<9cm,-4cm>*+{Q\tensor{}P'\tensor{}T\tensor{}T}="5", {"p"
\ar@{-->}^-{\alpha \tensorbajo{(A:\,T)} \beta} "1"}, {"p"
\ar@{->}_-{P\tensor{}1\tensor{}P'\tensor{}T} "2"}, {"2"
\ar@{->}_-{\alpha\tensor{}P' \tensor{}T} "4"}, {"4"
\ar@{->}_-{Q\tensor{}\ppp'\tensor{}T} "5"}, {"5"
\ar@{->}_-{Q\tensor{}\beta\tensor{}T} "3"}, {"3"
\ar@{->}_-{Q\tensor{}Q'\tensor{}\mu} "1"}
\endxy
\]
or equivalently, $\alpha \tensorbajo{(A:\,T)} \beta\,=\,
(Q\tensor{}Q'\tensor{}\mu) \circ (Q\tensor{}\qqq'\tensor{}T) \circ
(\alpha\tensor{}\beta) \circ (P\tensor{}1\tensor{}P'\tensor{}T)$.
The identity object of this multiplication is proportioned by the
pair $(T\tensor{}A\cong A\tensor{}T,A)$.

\begin{remark}
Let $(A:T)$ be any rings extension and $P$ an $A$--bimodule. Then
one can easily check that the following statements are equivalent.
\begin{enumerate}[(i)]
\item $P\tensor{A}T$ is a $T$-bimodule with right $T$-action
$\Sf{r}_{P\tensor{}T}\,=\, P\tensor{}\mu$;

\item there is an $A$-bilinear map $\ppp: T\tensor{}P \to
P\tensor{}T$ such that $(P,\ppp)$ is an object of the category
$\rR_{(A:\,T)}$.
\end{enumerate}
\end{remark}

\begin{lemma}\label{FT}
Let $(P,\ppp)$ and $(Q,\qqq)$ two objects of the category
$\rR_{(A:\,T)}$, and $f:P\to Q$ an $A$--bilinear map. the following
conditions are equivalent
\begin{enumerate}[(a)]
\item $f\tensor{A}T: (P,\ppp) \to (Q,\qqq)$ is a morphism in
$\rR_{(A:\,T)}$;

\item $f$ satisfies the equality $$\qqq \circ (T\tensor{}f)
\,\,=\,\, (f\tensor{}T) \circ \ppp.$$
\end{enumerate}
\end{lemma}
\begin{proof}
Straightforward.
\end{proof}

\begin{definition}\label{Ring-w}
Let $(A:T)$ be a rings extension, and consider its associated
monoidal category $\rR_{(A:\,T)}$. A \emph{wreath} over $T$ or
$T$--\emph{wreath} is a an algebra in $\rR_{(A:T)}$, and
\emph{cowreath} or $T$--\emph{cowreath} is a coalgebra in
$\rR_{(A:T)}$.
\end{definition}

\begin{remark}\label{LT}
As in the case of corings, in Definition \ref{Ring-w} we are
defining in fact a right wreath and right cowreath. The left
versions of those definitions are given in the left monoidal
category $\lL_{(A:T)}$ whose object are pairs $(\uuu, U)$
consisting of an $A$--bimodule $U$ and $A$--bilinear map $\uuu:
U\tensor{A}T \to T\tensor{A}U$ satisfying the equalities
\begin{eqnarray}
  \uuu \circ (U\tensor{}1_T) &=& 1_T\tensor{}U \\
  \uuu \circ (U\tensor{}\mu) &=& (\mu\tensor{}U) \circ
  (T\tensor{}\uuu) \circ (\uuu\tensor{}T).
\end{eqnarray}
The $\KK$--modules morphisms are ${\rm
Hom}_{\lL_{(A:\,T)}}\lr{(\uuu,U),\, (\uuu',U')}\,=\, {\rm
Hom}_{T-T}\lr{T\tensor{A}U,\, T\tensor{A}U'}$, where
$T\tensor{A}U$ and $T\tensor{A}U'$ are endowed with a structure of
$T$--bimodule with right $T$--action given, respectively, by
$\Sf{r}_{T\tensor{A}U}\,=\,(\mu\tensor{A}U) \circ
(T\tensor{A}\uuu)$ and
$\Sf{r}_{T\tensor{A}U'}\,=\,(\mu\tensor{A}U') \circ
(T\tensor{A}\uuu')$. The multiplications in this monoidal category
are given as follows: Given $\alpha: (\uuu,U) \to (\vvv,V)$ and
$\beta: (\uuu',U') \to (\vvv',V')$ two morphisms in
$\lL_{(A:\,T)}$, the horizontal multiplication is defined by $$
(\uuu,U) \tensorbajo{(A:\,T)} (\uuu',U') \,=\,
\lr{(\uuu\tensor{A}U') \circ (U\tensor{A}\uuu'),\,U\tensor{A}U'}$$
and the vertical multiplication is defined by
$$\alpha\tensorbajo{(A:\,T)}\beta\,=\,
(\mu\tensor{A}V\tensor{A}V') \circ (T\tensor{A}\alpha\tensor{A}V')
\circ (T\tensor{A}\uuu\tensor{A}V') \circ
(T\tensor{A}U\tensor{A}\beta) \circ
(T\tensor{A}U\tensor{A}1\tensor{A}U').$$
\end{remark}

The dual version of Proposition \ref{Equivalent-Def} is the
following
\begin{proposition}\label{Ring-Equ-Def}
Let $(A:T)$ be a rings extension, and $(R,\rrr)$ an object of the
category $\rR_{(A:\,T)}$. The following statements are equivalent
\begin{enumerate}[(i)]
\item $(R,\rrr)$ is a $T$-wreath.

\item There is a $T$-bilinear maps $\bd{\eta}: T\to R\tensor{A}T$
and $\bd{\mu}: R\tensor{A}R\tensor{A}T \to R\tensor{A}T$ rendering
commutative the following diagrams
$$
\xy *+{R\tensor{}R\tensor{}T}="p", p+<3cm,0pt>*+{R\tensor{}T}="1",
p+<0pt,-2cm>*+{R\tensor{}T}="2",{"2" \ar@{->}^-{R\tensor{}\bd{\eta}}
"p"}, {"p" \ar@{->}^-{\bd{\mu}} "1"}, {"1" \ar@{=} "2"}
\endxy \qquad \xy *+{R\tensor{}R\tensor{}T}="p",
p+<4cm,0pt>*+{R\tensor{}T}="1", p+<4cm,-2cm>*+{T\tensor{}R}="2",
p+<0pt,-2cm>*+{R\tensor{}T\tensor{}R}="3", {"3"
\ar@{->}^-{R\tensor{}\rrr} "p"}, {"p" \ar@{->}^-{\bd{\mu}} "1"},
{"2" \ar@{->}_-{\rrr} "1"}, {"2" \ar@{->}^-{\bd{\eta}\tensor{}R}
"3"}
\endxy $$
$$
\xy *+{R\tensor{}R\tensor{}T\tensor{}R}="p",
p+<4cm,0pt>*+{R\tensor{}T\tensor{}R}="1",
p+<8cm,0pt>*+{R\tensor{}R\tensor{}T}="2",
p+<0pt,-2cm>*+{R\tensor{}R\tensor{}R\tensor{}T}="3",
p+<4cm,-2cm>*+{R\tensor{}R\tensor{}T}="4",
p+<8cm,-2cm>*+{R\tensor{}T}="5",  {"p"
\ar@{->}^-{\bd{\mu}\tensor{}R} "1"}, {"p"
\ar@{->}_-{R\tensor{}R\tensor{}\rrr} "3"}, {"1"
\ar@{->}^-{R\tensor{}\rrr} "2"}, {"2" \ar@{->}^{\bd{\mu}} "5"}, {"3"
\ar@{->}_-{R\tensor{}\bd{\mu}} "4"}, {"4" \ar@{->}_-{\bd{\mu}} "5"}
\endxy
$$
\end{enumerate}
\end{proposition}

The object and morphisms in the category of right modules over a
wreath are expressed in the following proposition

\begin{proposition}
Let $(R,\rrr)$ be a $T$--wreath with structure maps $\bd{\eta}: T
\to R\tensor{A}T$ and $\bd{\mu}: R\tensor{A}R\tensor{A}T \to
R\tensor{A}T$.
\begin{enumerate}[(a)]
\item Consider $(Y,\yyy)$ an object of $\rR_{(A:\,T)}$. The
following conditions are equivalent
\begin{enumerate}[(i)]
\item $(Y,\yyy)$ is right $(R,\rrr)$--module;

\item There is a $T$--bilinear map $\Sf{r}_{(Y,\,\yyy)}:
Y\tensor{A}R\tensor{A}T \to Y \tensor{A}T$ such that
$$\xymatrix@C=40pt@R=30pt{Y\tensor{}R\tensor{}T
\ar@{->}^-{\Sf{r}_{(Y,\,\yyy)}}[r]
 & Y\tensor{}T \\ Y\tensor{}T \ar@{=}[ur] \ar@{->}^-{Y\tensor{}\bd{\eta}}[u]} \,\,\,
\xymatrix@C=40pt@R=30pt{Y\tensor{}R\tensor{}T\tensor{}R
\ar@{->}^-{Y\tensor{}R\tensor{}\rrr}[r]
\ar@{->}_-{\Sf{r}_{(Y,\,\yyy)}\tensor{}R}[d]  &
Y\tensor{}R\tensor{}R\tensor{}T \ar@{->}^-{Y\tensor{}\bd{\mu}}[r]
& Y\tensor{}R\tensor{}T \ar@{->}^-{\Sf{r}_{(Y,\,\yyy)}}[d]   \\
Y\tensor{}T\tensor{}R \ar@{->}^-{Y\tensor{}\rrr}[r] &
Y\tensor{}R\tensor{}T \ar@{->}^-{\Sf{r}_{(Y,\,\yyy)}}[r]&
Y\tensor{}T  }$$ are commutative diagrams.
\end{enumerate}

\item Given two right $(R,\rrr)$--modules $(Y,\yyy)$ and
$(Y',\yyy')$. A morphism $\Sf{g}: (Y,\yyy) \to (Y',\yyy')$ in
$\rR_{(A:\,T)}$ is a morphism of right $(R,\rrr)$--modules if and
only if $\,\Sf{g}$ turns commutative the following diagram
$$\xymatrix@C=50pt@R=20pt{ & Y\tensor{}T\tensor{}R \ar@{->}^-{\Sf{g}\tensor{}R}[r]
\ar@{->}_-{Y\tensor{}\rrr}[dl] &
Y'\tensor{}T\tensor{}R \ar@{->}^-{Y'\tensor{}\rrr}[dr] & \\
Y\tensor{}R\tensor{}T \ar@{->}_-{\Sf{r}_{(Y,\,\yyy)}}[dr] & & &
Y'\tensor{}R\tensor{}T
\ar@{->}^-{\Sf{r}_{(Y',\,\yyy')}}[dl] \\
& Y\tensor{}T \ar@{->}^-{\Sf{g}}[r]
 & Y'\tensor{}T & }$$
\end{enumerate}
\end{proposition}

Analogously we expresse the objects and morphisms in the category of
left modules over a wreath

\begin{proposition}\label{leftW}
Let $(R, \rrr)$ be a $T$--wreath with structure maps $\bd{\eta}: T
\to R\tensor{A}T$ and $\bd{\mu}: R\tensor{A}R\tensor{A}T \to
R\tensor{A}T$.
\begin{enumerate}[(a)]
\item Consider $(Y,\yyy)$ an object of $\rR_{(A:\,T)}$. The
following are equivalent
\begin{enumerate}[(i)]
\item $(Y,\yyy)$ is left $(R,\rrr)$--module;

\item There is a $T$--bilinear map $\Sf{l}_{(Y,\yyy)}:
R\tensor{A}Y\tensor{A}T \to Y \tensor{A}T$ such that
$$\xymatrix@R=30pt{T\tensor{}Y \ar@{->}^-{\bd{\eta}\tensor{}Y}[r]
\ar@{->}_-{\yyy}[rrd] & R\tensor{}T\tensor{}Y
\ar@{->}^-{R\tensor{}\yyy}[r] & R\tensor{}Y\tensor{}T
\ar@{->}^-{\Sf{l}_{(Y,\,\yyy)}}[d] \\ & & Y\tensor{}T,} \,\,\,
\xymatrix@C=40pt@R=30pt{R\tensor{}R\tensor{}Y\tensor{}T
\ar@{->}^-{R\tensor{}\Sf{l}_{(Y,\,\yyy)}}[r]  &
R\tensor{}Y\tensor{}T \ar@{->}^-{\Sf{l}_{(Y,\,\yyy)}}[r] &
Y\tensor{}T   \\
R\tensor{}R\tensor{}T\tensor{}Y \ar@{->}^-{\bd{\mu}\tensor{}Y}[r]
\ar@{->}^-{R\tensor{}R\tensor{}\yyy}[u] & R\tensor{}T\tensor{}Y
\ar@{->}^-{R\tensor{}\yyy}[r]& R\tensor{}Y\tensor{}T
\ar@{->}_-{\Sf{l}_{(Y,\,\yyy)}}[u] }$$ are commutative diagrams.
\end{enumerate}

\item Given two left $(R,\rrr)$--modules $(Y,\yyy)$ and
$(Y',\yyy')$. A morphism $\Sf{g}: (Y,\yyy) \to (Y',\yyy')$ in
$\rR_{(A:\,T)}$ is a morphism of right $(R,\rrr)$--modules if and
only if $\,\Sf{g}$ turns commutative the following diagram
$$\xymatrix@C=50pt@R=30pt{ R\tensor{}Y\tensor{}T \ar@{->}_-{\Sf{l}_{(Y,\,\yyy)}}[d]
\ar@{->}^-{R\tensor{}\Sf{g}}[r]  & R\tensor{}Y'\tensor{}T \ar@{->}^-{\Sf{l}_{(Y',\,\yyy')}}[d] \\
Y\tensor{}T \ar@{->}^-{\Sf{g}}[r]
 & Y'\tensor{}T}$$
\end{enumerate}
\end{proposition}

Let $(R,\rrr)$ be a $T$-wreath and $(Y,\yyy)$ a right
$(R,\rrr)$-module and left $(R,\rrr)$-module with actions maps
$\Sf{l}_{(Y,\,\yyy)}: R\tensor{A}Y\tensor{A}T \to Y\tensor{A}T$
and $\Sf{r}_{(Y,\,\yyy)}: Y\tensor{A}R\tensor{A}T \to
Y\tensor{A}T$. Then $(Y,\yyy)$ is a $(R,\rrr)$--bimodule if and
only if
\begin{multline}
\Sf{l}_{(Y,\,\yyy)} \circ (R\tensor{}Y\tensor{}\mu) \circ
(R\tensor{}\Sf{r}_{(Y,\,\yyy)}\tensor{}T) \circ
(R\tensor{}Y\tensor{}R\tensor{}1\tensor{}T) \\ \,\,=\,\,
\Sf{r}_{(Y,\,\yyy)} \circ (Y\tensor{}R\tensor{}\mu) \circ
(Y\tensor{}\rrr\tensor{}T) \circ
(\Sf{l}_{(Y,\,\yyy)}\tensor{}R\tensor{}T) \circ
(R\tensor{}Y\tensor{}1\tensor{}R\tensor{}T).
\end{multline}

The wreath product is described by the following proposition

\begin{proposition}\label{w-prodct}
Let $(A:T)$ be any ring extension and $(R,\rrr)$ a $T$--wreath
with structure maps $\bd{\eta}:T \to R\tensor{A}T$ and $\bd{\mu}:
R\tensor{A}R\tensor{A}T \to R\tensor{A}T$. Then $(A:R\tensor{A}T)$
is a ring extension whose multiplication and unit are defined by
$$ \bd{\mu}' \,\,=\,\, (R\tensor{A}\mu)\circ (\bd{\mu}\tensor{A}T) \circ
(R\tensor{A}\rrr\tensor{A}T), \qquad \bd{\eta}' \,\,=\,\,
\bd{\eta} \circ 1.$$ Moreover, $\bd{\eta}: T \to R\tensor{A}T$ is
a morphism of $A$--rings.
\end{proposition}

\subsection{Functors connecting modules categories}
Given $(R,\rrr)$ a $T$--wreath, we denote by ${}_{(R,\,\rrr)}\mM$
its category of all right $(R,\rrr)$-module. As  in the case of
corings, we have a noncommutative diagram
\begin{equation}\label{digrama-d}
\xymatrix@C=40pt@R=40pt{{}_{(R,\rrr)}\mM \ar@{->}^-{\oO'}[rr]
\ar@<0,5ex>^-{\uU'}[dd] & &
{}_{R\tensor{}T}\mM_{T} \ar@<0,5ex>^-{(-)_{\eta}}[dd] \\ & & \\
\rR_{(A:\,T)} \ar@<0,5ex>^-{\wW'}[rr]
\ar@<0,5ex>^-{(R,\rrr)\tensorbajo{(A:\,T}-}[uu] & &
\ar@<0,5ex>^-{R\tensor{}-}[uu] \ar@<0,5ex>^-{\vV'}[ll] {}_T\mM_A,}
\end{equation}
where
\begin{enumerate}[$\bullet$]
\item the functor $\Oo': {}_{(R,\,\rrr)}\mM \to
{}_{R\tensor{}T}\mM_T$ is defined by $$\Oo'\lr{\lr{(Y,\yyy),
\Sf{l}_{(Y,\,\yyy)}}}\,=\,\lr{Y\tensor{A}T,
\Sf{l}_{Y\tensor{A}T}:= (Y\tensor{A}\mu) \circ
(\Sf{l}_{(Y,\,\yyy)}\tensor{A}T) \circ
(R\tensor{A}\yyy\tensor{A}X)}, \qquad \Oo'(f) =f;$$

\item the functors of the adjunction $\vV' \dashv \wW'$:
$\xymatrix{ \vV': {}_T\mM_A \ar@<0,5ex>[r] & \ar@<0,5ex>[l] \rR_{(
A:\,T)}: \wW'}$ are defined by
$$ \wW'\lr{(P,\ppp)}\,=\, P\tensor{A}T, \qquad \wW'(f)\,=\, f$$ and
$$\vV'(X,\Sf{l}_X)\,=\, \lr{X,\xxx:= (X\tensor{A}1) \circ
\Sf{l}_X}, \qquad \vV'(f)\,=\,f\tensor{A}T;$$

\item the adjunction $R\tensor{A}- \dashv (-)_{\eta}$ is the usual
adjunction associated to the $A$-rings morphisms of Proposition
\ref{w-prodct};

\item the adjunction $(R,\,\rrr)\tensorbajo{(A:\,T)}-\dashv \uU'$
is the usual adjunction for the category of right modules over an
algebra defined in monoidal category.
\end{enumerate}

\subsection{Twisted tensor product algebras are wreath
products}\label{TWALG} Let $T$ and $R$ two $A$--rings with
multiplications and units $\mu_T$, $1_T$ and $\mu_R$, $1_R$. Assume
that there is an $A$--bilinear map $\rrr: T\tensor{A}R  \to
R\tensor{A}T$ satisfying
\begin{eqnarray}
\rrr \circ (1_T\tensor{}R) &=& R\tensor{}1_T \label{TTP-1}\\
  \rrr \circ (\mu_T\tensor{}R) &=& (R\tensor{}\mu_T) \circ
  (\rrr\tensor{}T) \circ (T\tensor{}\rrr) \label{TTP-2}\\
  \rrr \circ (T\tensor{}1_R) &=& 1_R\tensor{}T \label{TTP-3} \\
  \rrr \circ (T\tensor{}\mu_R) &=& (\mu_R\tensor{}T) \circ (R\tensor{}\rrr) \circ (\rrr\tensor{}R). \label{TTP-4}
\end{eqnarray}
Equations \eqref{TTP-1} and \eqref{TTP-2} say that $(R,\rrr)$ is
an object of the category $\rR_{(A:\,T)}$. While equations
\eqref{TTP-3} and \eqref{TTP-4} say that $(\rrr,T)$ is an object
of the category $\lL_{(A:\,R)}$.\\ Taking the maps
$\bd{\eta}:=1_R\tensor{A}T: T \to R\tensor{A}T$ and
$\bd{\mu}:=\mu_R\tensor{A}T: R\tensor{A}R\tensor{A}T \to
R\tensor{A}T$, we can easily check using associativity and unitary
properties of $\mu_R$ and $1_R$, that $\bd{\eta}$ and $\bd{\mu}$
satisfy the commutativity of the diagrams stated in Proposition
\ref{Ring-Equ-Def}(ii). Equations \eqref{TTP-3} and \eqref{TTP-4}
show that both $\bd{\eta}$ and $\bd{\mu}$ are $T$--bilinear maps.
That is $(R,\rrr)$ a $T$-wreath with structure maps $\bd{\eta}$
and $\bd{\mu}$. The wreath product $R\tensor{A}T$ is by
Proposition \ref{Ring-Equ-Def} an $A$-ring extension of $T$ with
multiplication and unit
$$\bd{\mu}'\,\,=\,\, (\mu_R\tensor{A}\mu_T) \circ
(R\tensor{A}\rrr\tensor{A}T), \qquad 1_{R\tensor{A}T}\,\,=\,\,
1_R\tensor{A}1_T.$$ Of course $(\rrr,T) \in \lL_{(A:\,R)}$ can be
also considered as an $R$-wreath with structure maps
$R\tensor{A}1_T$ and $R\tensor{A}\mu_T$, and will leads to the
wreath product $R\tensor{A}T$ but this time a ring extension of
$R$.

In the commutative case (i.e. $A=\KK$), these algebras were
refereed in the literature by \emph{twisted tensor product
algebras}, and were intensively studied by several Mathematicians,
see \cite{Tambara:1990, Van Daele/Van Keer:1994,
Cap/Schichl/Vanzura:1995, Caenepeel/Ion/Militaru/Zhu:2000,
LopezPena:2006} and references cited there.

\subsection{Example: Ore extensions are wreath products}
(\cite[Example 2.11]{Caenepeel/Ion/Militaru/Zhu:2000}). We prove in
this subsection that the classical Ore extension
\cite{Goodearl/Warfield:1989, McConnell/Robson:1987} constructed by
using left skew derivations are in fact a product of wreath defined
over a commutative polynomials rings with one variable. So consider
such polynomials ring $T=\KK[X]$, and let $B$ be any ring (i.e.
$\KK$--algebra). Assume that $\delta$ is a left $\sigma$-derivation
of $B$, where $\sigma$ is a endomorphism of rings of $B$. That is
$\delta$ is a $\KK$--linear maps obeying the rule $\delta(bb')\,=\,
\delta(b)b' + \sigma(b)\delta(b')$, for every $b,b' \in B$. We
define by induction the following map $\bbb: T\tensor{\KK}B \to
B\tensor{\KK}T$,
\begin{eqnarray*}
  \bbb(1\tensor{}b) &=& b\tensor{}1 \\
  \bbb(X\tensor{}b) &=& \sigma(b)\tensor{}X + \delta(b)\tensor{}1,
\end{eqnarray*}
and if $\bbb(X^n\tensor{}b)\,=\, \sum_{1 \leq i \leq n}
b_i\tensor{}X^i$, for some $n\geq 1$, then $$\bbb(X^{n+1}
\tensor{}b) \,\,=\,\, \sum_{i=1}^{n} \lr{
\sigma(b_i)\tensor{}X^{i+1} + \delta(b_i)\tensor{}X^i}$$ Following
to the proof of \cite[Proposition 1.10]{Goodearl/Warfield:1989},
consider the $\KK$--linear endomorphisms ring $E={\rm
End}_{\KK}\lr{B[t]}$ of a commutative polynomial ring over $B$.
Clearly the ring $B$ is identified with its image in $E$ by using
left multiplications. The maps $\sigma$ and $\delta$ are extended to
maps in $E$ where $\sigma(bt^i)\,=\,\sigma(b)t^i$ and
$\delta(bt^i)\,=\,\delta(b)t^i$, for all $b \in B$ and
$i=0,1,\cdots$. Denote by $Y$ the element of $E$ defined by
$Y(f)\,=\, \sigma(f)t + \delta(f)$, for all $f \in B[t]$, and
construct a map $$\xymatrix@R=0pt{ B\tensor{\KK}\KK[X]
\ar@{->}^-{\tau}[r] &  {\rm End}_{\KK}\lr{B[t]} \\ b \tensor{} X^n
\ar@{|->}[r] &  bY^n}$$ $\tau$ is in fact injective since the left
$B$--submodule generated by the set $\{Y^n\}_{n=0,1,\cdots}$ is a
free. By the associativity of $E$, we have $\tau \circ \bbb \circ
(\mu\tensor{}B) \,=\, \tau \circ (B\tensor{}\mu) \circ
(\bbb\tensor{}T) \circ (T\tensor{}\bbb)$. Whence $\bbb \circ
(\mu\tensor{}B) \,=\, (B\tensor{}\mu) \circ (\bbb\tensor{}T) \circ
(T\tensor{}\bbb)$, thus $(\bbb,B)$ is an object of the category
$\rR_{(\KK:\,T)}$. Using again the injectivity of the map $\tau$, we
can prove that $(\bbb,B)$ is a $T$--wreath with structure maps
$\bd{\eta}: T \to B\tensor{\KK}T$ sending $X^n \mapsto
1\tensor{}X^n$  and $\bd{\mu}: B\tensor{\KK}B\tensor{\KK}T \to
B\tensor{\KK}T$ sending $b\tensor{}b'\tensor{}X^n \mapsto
bb'\tensor{}X^n$. In this way the product of this wreath is a
$\KK$--algebra with underlying $\KK$--module $B\tensor{\KK}\KK[X]$
and multiplication $$ (b\tensor{}X^n) (b'\tensor{}X^m)\,\,=\,\,
\sum_{i=1}^n bb_i'\tensor{}X^{i+m}$$ where
$\bbb(X^n\tensor{}b')\,=\, \sum_{1\leq i \leq n}b_i' \tensor{}X^i$.
This algebra is in fact an extension of $B$ via the map
$-\tensor{}1: B \to B\tensor{\KK}\KK[X]$. Now given any rings
extension $\phi: B \to S$ and assume that there exists an element $Z
\in S$ such that $Zb\,=\, \sigma(b)Z+ \delta(b)$, for all $b \in B$.
This condition leads to construct a ring extension $\bara{\phi}:
B\tensor{\KK}\KK[X] \to S$ sending $b\tensor{}X^n \mapsto bZ^n$. Its
clear now that $\phi\,=\, \bara{\phi} \circ (-\tensor{}1)$. In
conclusion the $\KK$--algebra $B\tensor{\KK}\KK[X]$ satisfies the
universal condition on Ore extension, and thus $B\tensor{\KK}\KK[X]
\,=\, B[Y;\sigma;\delta]$. Notice that $B\tensor{\KK}\KK[X]$ is
isomorphic to the subalgebra $\sum_{i=0,1,\cdots} BY^i$ of $E$.\\
Here in fact we have constructed a wreath product with commutative
base ring. An example of wreath with non-commutative base ring can
be constructed as above using an iterated Ore extensions over a
non commutative ring $A$. That is one can prove that the iterated
Ore extension of type
$A[X_1;\sigma_1;\delta_1][X_2;\sigma_2;\delta_2]$ is a wreath
product with base ring $A$.

\subsection{\v{C}ap et al construction of twisted modules} In this
subsection we review \v{C}ap et al \cite[Section
3]{Cap/Schichl/Vanzura:1995} constructions of what they call a
twisted modules over a twisted tensor algebra. The problem concerned
in \cite[Section 3]{Cap/Schichl/Vanzura:1995} is rephrased in the
non commutative case as follows: Given $A$ a non commutative base
ring, $T$, $R$ and $\rrr:T\tensor{A}R \to R\tensor{A}T$ as in
subsection \ref{TWALG}, a $(T, A)$-bimodule $X$, an $(R,
A)$-bimodule $Y$, can we make $X\tensor{A}Y$ into a left
$R\tensor{\rrr}T$--module in a way which is compatible with the
inclusion of $R$, i.e. such that $(r\tensor{}1_T) .
(x\tensor{}y)\,=\, (rx)\tensor{}y$, for every $r \in R$, $x \in X$
and $y \in Y$? where $R\tensor{\rrr}T:=R\tensor{A}T$ denotes the
wreath product of subsection \ref{TWALG}. The left action which the
authors proposed is the following one
$$\xymatrix@C=50pt{ \Sf{l}_{X\tensor{}Y}: R\tensor{}T
\tensor{}X\tensor{}Y \ar@{->}^-{R\tensor{}\xxx\tensor{}Y}[r] &
R\tensor{}X\tensor{}T\tensor{}Y
\ar@{->}^-{\Sf{l}_X\tensor{}\Sf{l}_Y}[r] & X\tensor{}Y}$$ where
$\Sf{l}_X: T\tensor{A}X \to X$ and $\Sf{l}_Y:R\tensor{A}Y \to Y$
are, respectively, the $A$--bilinear left action map of $X$ and $Y$,
and $\xxx: T\tensor{A}X \to X\tensor{A}T$ is some $A$-bilinear map.
A suffisent condition which should satisfies $\xxx$ in order to
answer positively to the above question using the map
$\Sf{l}_{X\tensor{}Y}$, was given in the commutative case in
\cite[3.6 and 3.7]{Cap/Schichl/Vanzura:1995} and says the following:
$\xxx$ is said to be a left module twisting map if and only if
\begin{eqnarray}
  \xxx \circ (1_T \tensor{}X) &=& X\tensor{}1_R \label{Cap-1} \\
  \xxx \circ (\mu_T\tensor{}X) &=& (X\tensor{}\mu_T) \circ (\xxx\tensor{}T) \circ (T\tensor{}\xxx) \label{Cap-2} \\
  \xxx \circ (T\tensor{}\Sf{l}_X) &=& (\Sf{l}_X\tensor{T}) \circ
  (R\tensor{}\xxx) \circ (\rrr\tensor{}X). \label{Cap-3}
\end{eqnarray}
\cite[Theorem 3.8]{Cap/Schichl/Vanzura:1995} says: If $\xxx$ is a
left module twisting map, the $\Sf{l}_{X\tensor{}Y}$ gives the
answer to the above question. A reciprocate implication was also
given in that Theorem: If $X$ is projective and for one left
faithful module $Y$ the map $\Sf{l}_{X\tensor{}Y}$ defines a left
action which compatible with the inclusion of $R$ then $\xxx$ is a
left module twisting map.

Let us traduce the above constructions in our context. First of
all, it is obvious that equations \eqref{Cap-1} and \eqref{Cap-2}
say that $(X,\xxx)$ is actually an object of the category
$\rR_{(A:\,T)}$. Take the $A$-bilinear map $\Sf{l}_{(X,\,\xxx)}: =
\Sf{l}_X\tensor{A}T: R\tensor{A}X\tensor{A}T \to X\tensor{A}T$. It
is easily seen that this map satisfies the commutativity  of the
diagrams sated in Proposition \ref{leftW}(ii). By Lemma \ref{FT},
$\Sf{l}_{(X,\,\xxx)}$ is a $T$--bilinear map if and only if
equation \eqref{Cap-3} is fulfilled. Therefore, $\xxx$ is a left
module twisting map if and only if $(X,\xxx)$ is a left
$(R,\rrr)$-module with action $\Sf{l}_{(X,\,\xxx)} =
\Sf{l}_X\tensor{A}T$. Thus what \v{C}ap et al were constructing is
just an induced left module over the wreath $(R,\rrr)$. Next, we
formulate the non commutative version  of \cite[Theorem
3.8]{Cap/Schichl/Vanzura:1995}.
\begin{proposition}\label{Cap-5}
Let $(R,\rrr)$ be the $T$-wreath of subsection \ref{TWALG} and $X$
an $(R,A)$--bimodule with left $R$-action $\Sf{l}_X$. Assume that
$(X,\xxx)$ is also a left $(R,\rrr)$-module with action
$\Sf{l}_{(X,\,\xxx)}: R\tensor{A}X\tensor{A}T \to X\tensor{A}T$ and
that $$\Sf{l}_{(X,\,\xxx)} \circ (R\tensor{}X\tensor{}1_T) \,\,
=\,\, (X\tensor{}1_T) \circ \Sf{l}_X .$$ Then, for every
$(T,A)$--bimodule $Y$, the $A$-bilinear map
$$\xymatrix@C=60pt{ R\tensor{}T\tensor{}X\tensor{}Y \ar@{->}^-{R\tensor{}\xxx\tensor{}T}[r]
\ar@{-->}_-{\Sf{l}_{X\tensor{}Y}}[drr] &
R\tensor{}X\tensor{}T\tensor{}Y
\ar@{->}^-{\Sf{l}_{(X,\,\xxx)}\tensor{}Y}[r] & X\tensor{}T\tensor{}Y
\ar@{->}^-{X\tensor{}\Sf{l}_Y}[d] \\  & & X\tensor{}Y}$$ define a
left $(R\tensor{\rrr}T)$-action which is compatible with the
inclusion of $R$.

Conversely, If $X\tensor{A}\mu_T$ preserves equalizers,
$-\tensor{A}T$ is faithful functor, and $\Sf{l}_{X\tensor{A}T}:=
(\Sf{l}_X\tensor{A}\mu_T) \circ (R\tensor{A}\xxx\tensor{A}T)$
define a left $(R\tensor{\rrr}T)$-action on $X\tensor{A}T$ which
is compatible with the inclusion of $R$ for some $A$--bilinear map
$\xxx:T\tensor{A}X \to X\tensor{A}T$ then $(X,\xxx)$ is a left
$(R,\rrr)$--module with action $\Sf{l}_{(X,\,\xxx)}=\Sf{l}_X
\tensor{A}T$. In particular $\xxx$ is a left module twisting map.
\end{proposition}
\begin{proof}
The unitary property of $\Sf{l}_{X\tensor{}Y}$ comes from that of
$\Sf{l}_Y$ and of $\Sf{l}_{(X,\,\xxx)}$( i.e. the first diagram in
Proposition \ref{leftW}(ii)). The associativity of $\Sf{l}_Y$ and of
$\Sf{l}_{(X,\,\xxx)}$(i.e. the second diagram in Proposition
\ref{leftW}(ii)) lead to that of $\Sf{l}_{X\tensor{}Y}$, taking into
the account that the structure maps of the wreath $(R,\rrr)$ are
$1_R\tensor{A}T$ and $\mu_R\tensor{A}T$. The proof of the
reciprocate implication is lifted to the reader.
\end{proof}

The right version of Proposition \ref{Cap-5} is expressed in the
monoidal category $\lL_{(A:\,R)}$ using the $R$-wreath $(\rrr,T)$.
If we combine both versions then we get a criterion on
$R\tensor{\rrr}T$-bimodules of the form $X\tensor{A}Y$ as in
\cite[Proposition 3.13]{Cap/Schichl/Vanzura:1995}.

\begin{proposition}
Let $(R,\rrr)$ and $(\rrr,T)$, respectively, the $T$-wreath and
$R$-wreath of subsection \ref{TWALG}. Given $X$ an $R$-bimodule,
and $V$ a $T$-bimodule with actions $\Sf{l}_X$, $\Sf{r}_X$ and
$\Sf{l}_V$, $\Sf{r}_V$. Assume that $(X,\xxx)$ is also a left
$(R,\rrr)$-module and $(\vvv,V)$ a right $(\rrr,T)$-module with
actions, respectively, $\Sf{l}_{(X,\xxx)}$ and
$\Sf{l}_{(\vvv,V)}$, and consider the maps
$$\Sf{l}_{X\tensor{A}V}\,=\, (M\tensor{A}\Sf{l}_V) \circ
(\Sf{l}_{(X,\,\xxx)}\tensor{A}V) \circ
(R\tensor{A}\xxx\tensor{A}V), \quad \Sf{r}_{X\tensor{A}V}\,=\,
(\Sf{r}_X\tensor{A}V) \circ (X\tensor{A}\Sf{r}_{(\vvv,V)}) \circ
(X\tensor{A}\vvv\tensor{A}T)$$ If $X\tensor{A}V$ is
$(R\tensor{A}T)$-bimodule with action $\Sf{l}_{X\tensor{A}V}$ and
$\Sf{r}_{X\tensor{A}V}$ then the following diagram is commutative
\begin{equation}\label{Cap-6}
\xymatrix@C=50pt@R=30pt{T\tensor{}X\tensor{}V\tensor{}R
\ar@{->}^-{\xxx\tensor{}V\tensor{}R}[r]
\ar@{->}_-{T\tensor{}X\tensor{}\vvv}[d] &
X\tensor{}T\tensor{}V\tensor{}R \ar@{->}^-{X\tensor{}\Sf{l}_V}[r]
& X\tensor{}V\tensor{}R \ar@{->}^-{X\tensor{}\vvv}[d] \\
T\tensor{}X\tensor{}R\tensor{}V \ar@{->}_-{T\tensor{}\Sf{r}_X\tensor{}V}[d] & &
X\tensor{}R\tensor{}V \ar@{->}^-{\Sf{r}_X\tensor{}V}[d] \\
T\tensor{}X\tensor{}V \ar@{->}_-{\xxx\tensor{}V}[r] &
X\tensor{}T\tensor{}V \ar@{->}_-{X\tensor{}\Sf{l}_V}[r] &
X\tensor{}V}
\end{equation}
This $(R\tensor{A}T)$-biaction on $X\tensor{A}V$ is left
compatible with the inclusion of $R$ and right compatible with the
inclusion of $T$ whenever $\Sf{l}_{(X,\,\xxx)}$ and
$\Sf{r}_{(\vvv,V)}$ satisfy the equalities $\Sf{l}_{(X,\,\xxx)}
\circ (R\tensor{A}X\tensor{A}1_T) \, =\, (X\tensor{A}1_T) \circ
\Sf{l}_X$ and $\Sf{r}_{(\vvv,V)} \circ (1_R\tensor{A}X\tensor{A}T)
\,\, =\,\, (1_R\tensor{A}V) \circ \Sf{r}_Y$.

Conversely, if $\Sf{l}_{(X,\,\xxx)}\,=\, \Sf{l}_X\tensor{A}T$,
$\Sf{r}_{(\vvv,V)}\,=\, R\tensor{A}\Sf{r}_V$, and the diagram of
equation \eqref{Cap-6} commutes then $X\tensor{A}V$ is an
$(R\tensor{A}T)$-bimodule.
\end{proposition}
\begin{proof}
Analogue to that of \cite[Proposition
3.13]{Cap/Schichl/Vanzura:1995}.
\end{proof}

\providecommand{\bysame}{\leavevmode\hbox
to3em{\hrulefill}\thinspace}
\providecommand{\MR}{\relax\ifhmode\unskip\space\fi MR }
\providecommand{\MRhref}[2]{
} \providecommand{\href}[2]{#2}

\end{document}